\documentclass[11pt]{amsart}
\usepackage{amssymb}    
\usepackage{amsmath}    
\usepackage{amsthm}     
\usepackage{epsfig}             
\begin{document}

\newcommand{\M}{\widetilde{M}}
\newcommand{\minus}{\setminus}
\newcommand{\cross}{\times}
\newcommand{\bndry}{\partial{\M}}
\newcommand{\st}{\mid}
\newcommand{\comp}{\circ}
\newcommand{\norm}{\Arrowvert}
\newcommand{\ra}{\rightarrow}
\newcommand{\R}{\mathbb{R}}
\newcommand{\Q}{\mathbb{Q}}
\newcommand{\C}{\mathbb{C}}
\newcommand{\Hy}{\mathbb{H}}
\newcommand{\Ot}{\widetilde{\Omega}}
\newcommand{\Et}{\widetilde{E}}
\newcommand{\At}{\widetilde{A}}
\newcommand{\Bt}{\widetilde{B}}
\newcommand{\Ob}{\overline{\Omega}}
\renewcommand{\epsilon}{\varepsilon}
\newcommand{\nin}{\not\in}
\renewcommand{\qed}{\square}
\renewcommand{\phi}{\varphi}
\newcommand{\ri}{\underline{\rho}}
\newcommand{\rs}{\overline{\rho}}
\newcommand{\te}{t \epsilon}
\newcommand{\cc}{d_{cc}}
\newcommand{\fcc}{f_{cc}}
\newcommand{\db}{\overline{d}}
\newcommand{\pt}{\|_t}
\newcommand{\Obar}{\overline{\Omega}}
\newcommand{\Hams}{Hamenst$\ddot{\text{a}}$dt's }
\newcommand{\hr}{$h$-rank}
\newcommand{\Ham}{Hamenst$\ddot{\text{a}}$dt }
\renewcommand{\phi}{\varphi}
\newcommand{\bmv}{\partial{\M} \setminus \{v_{-\infty}\}}
\newcommand{\bmvw}{\partial{\M} \setminus \{v_{-\infty},w_{-\infty}\}}
\newcommand{\ghlim}{\text{GH-lim}}
\newcommand{\wlim}{\text{$\omega$-lim}}
\newcommand{\Ha}{\mathcal{H}}
\newcommand{\ip}{<\cdot,\cdot>}
\newcommand{\pf}{\noindent {\em Proof: }}
\newcommand{\df}{\underset{def}{\equiv}}
\newtheorem{Pro}{Proposition}[section]
\newtheorem{Lem}[Pro]{Lemma}
\newtheorem{Example}[Pro]{Example}
\newtheorem{Sub}[Pro]{Sublemma}
\newtheorem{Thm}[Pro]{Theorem}
\newtheorem{MThm}{Theorem}
\renewcommand{\theMThm}{\Alph{MThm}}
\newtheorem{Def}[Pro]{Definition}
\newtheorem{Not}{Notation}
\newtheorem{Claim}{Claim}
\newtheorem{Cor}[Pro]{Corollary}
\newtheorem{Ithm}{Theorem}
\renewcommand{\theIthm}{\arabic{chapter}.\arabic{Ithm}}
\newtheorem{IndThm}{Theorem}
\newtheorem{Idef}{Definition}
\renewcommand{\theIdef}{\arabic{chapter}.\arabic{Idef}}
\newcommand{\hook}{\lfloor}
\newcommand{\bG}{\boldsymbol{G}}
\newcommand{\nuX}{\boldsymbol{\nu}_H}
\newcommand{\nuXp}{\boldsymbol{\nu}_H^\perp}
\newcommand{\Up}{\boldsymbol{\mathcal Y}_H}
\newcommand{\n}{\boldsymbol \nu}
\newcommand{\sigmau}{\boldsymbol{\sigma}^u_H}
\newcommand{\di}{\delta_{X,i}}
\newcommand{\del}{\delta_H}
\newcommand{\nui}{\nu_{H,i}}
\newcommand{\nuj}{\nu_{H,j}}
\newcommand{\dej}{\delta_{H,j}}
\newcommand{\til}{\tilde{\boldsymbol{\nu}}_H}
\newcommand{\tilp}{\tilde{\boldsymbol{\nu}}^\perp_H}
\newcommand{\CO}{C^\infty_0( \Omega)}
\newcommand{\Rn}{\mathbb R^n}
\newcommand{\Rm}{\mathbb R^m}
\newcommand{\Om}{\Omega}
\newcommand{\Hn}{\mathbb H^n}
\newcommand{\p}{\partial}
\newcommand{\bg}{\mathfrak g}
\newcommand{\bz}{\mathfrak z}
\newcommand{\bv}{\mathfrak v}
\newcommand{\e}{\epsilon}
\newcommand{\al}{\alpha}
\newcommand{\pb}{\overline{p}}
\newcommand{\qb}{\overline{q}}
\newcommand{\Zb}{\overline{Z}}
\newcommand{\gm}{\mathbb{G}}
\newcommand{\tilL}{\tilde{\mathcal{L}}_z(r)}
\newcommand{\tilLg}{\tilde{\mathcal{L}}_{\gamma(s)}(r)}
\newcommand{\tilLp}{\tilde{\mathcal{L}}_z'(r)}
\newcommand{\tilLzero}{\tilde{\mathcal{L}}_z(0)}

\theoremstyle{plain}
\newtheorem{Rmk}[Pro]{Remark}
\newtheorem*{Hyp}{HYPOTHESIS}
\numberwithin{equation}{section}

\title[The Bernstein problem in the Heisenberg group]{The Bernstein problem in
the Heisenberg group}

\author{Nicola Garofalo}
\address{Purdue University, West Lafayette, IN 47906}
\email{garofalo@math.purdue.edu}
\thanks{The first author was supported by NSF grants No. DMS-0070492 and No. DMS-0300477}

\author{Scott D. Pauls}
\address{Dartmouth College, Hanover, NH 03755}
\email{scott.pauls@dartmouth.edu}
\thanks{The second author is partially supported by NSF grants
  No. DMS-9971563 and No. DMS-0306752}


\begin{abstract}
We establish the following theorem of Bernstein type for the first
Heisenberg group $\mathbb H^1$: Let $S$ be a $C^2$ connected
$H$-minimal surface which is a graph over some plane $P$, then $S$ is
either a non-characteristic vertical plane, or its generalized seed
curve satisfies a type of constant curvature condition.
\end{abstract}

\maketitle

\tableofcontents

\section{\textbf{Introduction}}\label{S:intro}

\vskip 0.2in

The development of last century's mathematics has been enormously influenced by the desire to understand minimal surfaces. In such development the \emph{Bernstein problem} has played
 a central role. The latter states that a smooth complete minimal graph in $\mathbb R^{N+1}$, i.e., a $C^2$ entire solution $u:\mathbb R^N \to \R$ of the minimal surface equation
\begin{equation}\label{ms}
div\ \left(\frac{Du}{\sqrt{1 + |Du|^2}}\right) \ =\ 0\ ,
\end{equation}
must necessarily be an affine hyperplane. Bernstein himself \cite{Be} established this property when $N=2$. Almost fifty years later Fleming \cite{Fle} gave a completely
 different proof still for the case $N = 2$, and also conjectured that Bernstein result should be valid in any dimension.
 His approach inspired a new major development in the subject. In \cite{DG} De Giorgi was able to extend Bernstein's theorem to the case $N = 3$. Almgren \cite{Al} subsequently succeeded in treating the case $N = 4$, and Simons \cite{Sim} finally obtained a solution for $N\leq 7$. In their celebrated work \cite{BDG} Bombieri, De Giorgi and Giusti established the surprising result that the Bernstein property fails if $N\geq 8$, thus putting an end to the problem in Euclidean space. They proved in fact that: \emph{if $N\geq 8$ there exist complete minimal graphs in $\mathbb R^{N+1}$ which are not hyperplanes}.  For additional references and developments the reader should also consult \cite{O}, \cite{N1}, \cite{G}, \cite{MM}, \cite{Bo}, \cite{N2}, \cite{Si}, \cite{GMS}, \cite{CM}.

The main objective of this paper is investigating the Bernstein
problem in the first Heisenberg group $\mathbb H^1$. This is the
simplest and most important non-Abelian model of a class of
graded, nilpotent Lie groups known as Carnot groups. Such groups
arise naturally as tangent spaces to sub-Riemannian manifolds,
also known as Carnot-Carath\'eodory spaces, see \cite{F},
\cite{Stein}, \cite{Bel}, \cite{Gro1}, \cite{Gro2}, \cite{Mon},
\cite{G}. For $n\in \mathbb N$, the underlying manifold of the
Heisenberg group $\Hn$ is $\mathbb C^n \times \R$, with group law
\begin{equation}\label{gl}
g\ \circ\ g'\ =\ (z,t)\ \circ\ (z',t')\ =\ \left(z + z', t + t' - \frac{1}{2} \Im\ <z ,\overline{z'}>\right)\ .
\end{equation}

Here, we have let $z = x + i y, z' = x' + i y' \in \mathbb C^n$, $t , t' \in \R$, $<z , \overline{z'}> = \sum_{j=1}^n z_j \overline{z_j'}$. Denoting by $L_g(g') = g \circ g'$ the left-translation associated with \eqref{gl}, and by $(L_g)_*$ its differential, one readily recognizes that the generators of the (real) Heisenberg algebra $\mathfrak h_n$ are the left-invariant vector fields
\begin{equation}\label{vf}
X_j(g)\ =\ (L_g)_*(\frac{\partial}{ \partial x_j})\ =\ \frac{\partial}{\partial x_j} - \frac{y_j}{2}\ \frac{\partial}{\partial t}\ ,\quad\quad X_{n+j}(g)\ =\ (L_g)_*(\frac{\partial}{ \partial y_j})\ =\ \frac{\partial}{\partial y_j} + \frac{x_j}{2}\ \frac{\partial}{\partial t}\ ,
\end{equation}
$j = 1,..., n$. The sub-Laplacian on $\Hn$ is the second-order partial differential operator given by
 $\mathcal L = \sum_{j=1}^{2n} X_j^2$. Such operator is the real part of the Kohn sub-Laplacian in
  $\mathbb C^{n+1}$, see \cite{Stein}. One has
\begin{equation}\label{comm}
[X_i,X_{n+j}]\ =\  T\ \delta_{ij}\ ,\quad\quad\quad\quad i , j = 1, ... , n\ ,
\end{equation}
where $T = \partial/\partial t$ represents the characteristic direction, all other commutators being trivial. One can thus decompose the Heisenberg algebra as follows $\mathfrak h_n = V_1 \oplus V_2$, with $V_1 = \mathbb C^n \times \{0\}$, $V_2 = \{0\} \times \R$, and since \eqref{comm} implies $[V_1,V_1] = V_2$, we see that H\"ormander's finite rank condition \cite{H}
\begin{equation}\label{frc}
rank\ Lie\{X_1,...,X_{2n}\}\ \equiv\ dim\ \Hn\ =\ 2n + 1\ ,
\end{equation}
 is fulfilled at step $r = 2$. The natural non-isotropic dilations on $\Hn$ are $\delta_\lambda (z,t) = (\lambda z, \lambda^2 t)$, with corresponding homogeneous dimension  $Q = 2n + 2$. Noting that Lebesgue measure $dg = dz dt$ is a bi-invariant Haar measure on $\Hn$, one has in fact $d \delta_\lambda(g) = \lambda^Q dg$.

Minimal surfaces are surfaces of least perimeter, or area, among all those with the same boundary.
In a Carnot-Carath\'eodory space an appropriate sub-Riemannian version of perimeter generalizing the classical one due to De Giorgi was introduced in \cite{CDG},
and further developed in \cite{GN}, see also the papers \cite{BM} and \cite{FSS1} where two equivalent definitions were independently set forth. To provide the reader with
a broader perspective of the results in this paper, it may be useful to recall the relevant definitions in the setting of a Carnot group $\bG$
with a given subbundle $H\bG \subset T\bG$, and a distribution of smooth vector fields $X = \{X_1,...,X_m\}$ which is bracket-generating for $T\bG$.
Given an open set $\Om\subset \bG$, we let
\[
\mathcal F(\Om)\ =\ \{\zeta \in C^1_o(\Om,H\bG)\ \mid\ |\zeta|_\infty\ =\ \sup_{\Om}\ |\zeta|\ \leq 1\}\ .
\]

For a  function $u\in L^1_{loc}(\Om)$, the $H$-variation of $u$ with respect to $\Om$ is defined by
\[
Var_H(u;\Om)\ =\ \underset{\zeta\in \mathcal F(\Om)}{\sup}\ \int_{\bG} u\ \sum_{i=1}^m X_i \zeta_i\ dg\ .
\]

A function $u\in L^1(\Om)$ is called of bounded $H$-variation in $\Om$ if $Var_H(u;\Om) <\infty$. The space $BV_H(\Om)$ of functions with bounded $H$-variation in $\Om$, endowed with the norm
\[
||u||_{BV_H(\Om)}\ =\ ||u||_{L^1(\Om)}\ +\ Var_H(u;\Om)\ ,
\]
is a Banach space.
Let $E\subset \bG$ be a measurable set, $\Om$ be an open set. The $H$-perimeter of $E$ with respect to $\Om$ is defined by
\begin{equation}\label{Xper}
P_H(E;\Om)\ =\ Var_H(\chi_E;\Om)\ ,
\end{equation}
where $\chi_E$ denotes the indicator function of $E$.
An $H$-minimal surface in an open set $\Om \subset \bG$ was defined in \cite{GN} as the boundary of a set of least $H$-perimeter, among all those with the same boundaries outside $\Om$.
The existence of such surfaces and a measure theoretic solution of the corresponding Plateau problem were established in \cite{GN}.

In this paper, we do not study $H$-minimal surfaces from the geometric
measure theoretic viewpoint, but rather follow the classical
development of the subject. For us, an $H$-minimal surface will be a
$C^2$ hypersurface whose horizontal mean curvature $\mathcal H$, see Definition
\ref{D:Xmeancurv}, vanishes everywhere. This notion of a minimal surface was independently introduced in \cite{DGN} and \cite{Pauls:minimal}, and the present paper was inspired by these works.  Of course, the concept of horizontal mean curvature plays a central role here.
Just as the Gauss map guides the theory of
Euclidean minimal surfaces, the {\em horizontal Gauss map} dominates the study of $H$-minimal
surfaces.  To define this notion, we suppose that $\mathbb H^1$ has been endowed with a left-invariant Riemannian metric with respect to which $\{X_1,X_2,T\}$ constitute an orthonormal basis. Here, the vector fields $\{X_1,X_2\}$, defined in \eqref{vf}, span the horizontal subbundle $H\mathbb H^1\subset T\mathbb H^1$.
We consider an oriented surface $S\subset \mathbb H^1$ and let $\n$ be the Riemannian unit
normal to the surface.
We define the horizontal normal $\Up : S \rightarrow H\Hy^1$, relative to the basis $\{X_1,X_2\}$ of $H\mathbb H^1$,
by the formula
\[
\Up\ =\ <\n,X_1> X_1\ +\ <\n,X_2> X_2\ .
\]
The \emph{horizontal Gauss map} $\nuX : S \rightarrow \mathbb S^{1}$ is defined by
\begin{equation}\label{hgm0}
\nuX\ =\ \frac{\Up}{|\Up|}\ ,
\end{equation}
at every point of the set $\mathcal S_H = \{g\in S \mid |\Up(g)| \not= 0\}$. The following definition was introduced in \cite{DGN}.

\medskip

\begin{Def} \label{D:Xmeancurv}
The \emph{$H$-mean curvature} (or horizontal mean curvature) of $S$ at points of $\mathcal S_H$ is defined by
\[
\mathcal H\ =\ \sum_{j=1}^2\ \delta_{H,j}\ \nu_{H,j}\ ,
\]
where $\delta_{H,j}$, $j=1,2,$ indicate the tangential horizontal derivatives, see Definition \ref{D:delta}.
\end{Def}

\medskip

A $C^2$ surface $S\subset \mathbb H^1$ is called $H$-\emph{minimal} if it has $H$-mean curvature equal to zero.
It may be helpful for the reader to know that there exists a close connection between such notion of minimal surface and the
 geometric measure theoretic one introduced in \cite{GN}. One has in fact the following first variation formula (see Theorem 10.2 in \cite{DGN}).
Let $S \subset \mathbb H^1$ be a $C^2$ surface with
characteristic set $\Sigma$, and suppose that $h , k\in
C^\infty_o(S\setminus \Sigma)$. Consider the family of
surfaces $J_\lambda(S)$,
where for small $\lambda \in \mathbb R$ we have let
\begin{equation}\label{def} J_\lambda(g)\ =\ g\ +\ \lambda\ \bigg(h(g) \nuX +
k(g) T\bigg)\ ,\quad\quad\quad g\in S\ .
\end{equation}
 The first variation of the
$H$-perimeter with respect to the deformation \eqref{def} is given
by
\begin{equation}\label{fvH}
\frac{d}{d\lambda} P_H(J_\lambda(S))\Bigl|_{\lambda
= 0}\ =\ \int_{S} \mathcal H\ \left(h\ +\ \frac{T\phi}{W}\
k\right)\ d\sigma_H\ ,
\end{equation}
where as in \eqref{Sphi} we have denoted with $\phi$ a defining function for $S$, $W$ is the angle function defined in \eqref{isothermal}, and $d\sigma_H$ denotes the Borel
measure induced on $S$ by the $H$-perimeter defined in \eqref{Xper}.
In particular, we conclude from \eqref{fvH} that:  $S$ \emph{is stationary with respect to
\eqref{def} if and only if it is $H$-minimal}.

In \cite{Pauls:minimal}, the second author used the approximation
of $(\mathbb H^1,d)$ by Riemannian manifolds (in the
Gromov-Hausdorff topology) to characterize and investigate the
solutions of the minimal surface and Plateau problems in $(\mathbb
H^1,d)$ via approximations by solutions to the analogous problems
in the approximating Riemannian manifolds. In \cite{DGN}, D.
Danielli, the first named author, and D.M. Nhieu developed a
sub-Riemannian calculus for hypersurfaces in Carnot groups. In
\cite{DGN1} they also determined, under some conditions, the
isoperimetric sets in the Heisenberg group, i.e.,
 those bounded sets which minimize the $H$-perimeter under a volume costraint, and they proved that their boundaries are $C^2$ (but not $C^3$!) hypersurfaces with positive constant $H$-mean curvature.
In the same paper the following sub-Riemannian Bernstein problem was formulated: consider a complete $H$-minimal surface $S\subset \Hn$. Under which assumptions on $\Hn$ and $S$ is
 it true that $S$ must be a vertical hyperplane, i.e., there exist $(a,b) = (a_1,...,a_n, b_1,...,b_n) \in \mathbb R^{2n}\setminus \{0\}$, and $\gamma \in \R$, such that
\begin{equation}\label{conj}
S\ =\ \{(x,y,t)\in \Hn \mid <a,x> + <b,y>\ =\ \gamma \}\ ?
\end{equation}

One easily recognizes that the vertical hyperplanes in \eqref{conj} are $H$-minimal (however, \emph{any} hyperplane in $\Hn$ is $H$-minimal. In particular, such is the characteristic
 hyperplane $\{(x,y,t)\in \Hn \mid t = 0\}$). The following conjecture was proposed.

\vskip 0.2in

\noindent \textbf{Conjecture:} \emph{In the Heisenberg group $\Hn$, at least in low dimension, the Bernstein property should hold provided that the surface $S$ is an entire graph, and has empty characteristic locus. In particular, suppose that $S$ is a $C^2$ entire graph in the first Heisenberg group $\mathbb H^1$, and assume that $S$ has empty characteristic locus, then $S$ must be a vertical plane. I.e., there exist $a,b, \gamma \in \R$, with $a^2 + b^2 \not= 0$, such that}
\[
S\ =\ \{g = (x,y,t) \in \mathbb H^1 \mid a x\ +\ b y\ =\ \gamma \}\ .
\]

\vskip 0.2in

We recall that, given a $C^1$ hypersurface $S\subset \Hn$, a point $g_o\in S$ is
called characteristic if the vector fields \eqref{vf} which generate the
horizontal subbundle $H\Hn$ become tangent to $S$ at $g_o$, i.e., $H_{g_o}\Hn \subset T_{g_o}S$. The collection $\Sigma = \Sigma_S$ of all the characteristic points
 of $S$ is called \emph{the characteristic locus} of $S$. We note explicitly that the set $S_H$ introduced above is given by $S_H = S \setminus \Sigma$.

Some comments concerning the above conjecture are necessary. Substantial evidence seems in favor of it. On one hand, there is the close relation between
the Bernstein property and the classical Liouville theorem for
harmonic functions. Such connection continues to hold in the
sub-Riemannian setting. Now, the Liouville property in $\Hn$ presents a striking new phenomenon with
respect to the classical setting, namely that if $u$ is a bounded
entire solution of the sub-Laplacian in $\Hn$, then $u$ depends \emph{only} on the
horizontal variables $(x,y)\in \mathbb R^{2n}$. As a
consequence, such a function must, in fact, be an ordinary harmonic
function in $\mathbb R^{2n}$, and therefore by
the classical Liouville theorem it is constant.  Secondly, and perhaps more importantly, a basic result in \cite{FSS2} shows that when one adapts De Giorgi's method of the
blow-up to the sub-Riemannian setting of $\Hn$, one obtains in the
limit blow-up sets which are vertical (non-characteristic) planes as in \eqref{conj}. By imposing the non-characteristic assumption in the conjecture one rules thus out the undesired $H$-minimal characteristic hyperplanes such as $\{(x,y,t)\in \Hn \mid t = 0\}$.

In light of this evidence, our inital efforts went in the direction of proving the conjecture true. In the process of establishing its veracity, we have developed a basic representation result for a graph-like $H$-minimal surface which is based on the notions of seed curve and height function, see Theorem \ref{mrep}. While analyzing the various possibilities, however, we have made the striking discovery that the above conjecture is in fact not true.

\medskip

\noindent \textbf{Counterexample to the conjecture:} \emph{The real analytic surface
\begin{equation}\label{ce}
S\ =\ \{(x,y,t)\in \mathbb H^1 \mid y = - x \tan(\tanh(t))\ \}\ ,
\end{equation}
is an entire $H$-minimal graph, with empty characteristic locus, over the coordinate $(x,t)$-plane in $\mathbb H^1$.}

\medskip

This counterexample, which will be analyzed in detail in Sections \ref{S:examples} and \ref{S:return}, shows the
failure of the above formulated sub-Riemannian counterpart of the classical Bernstein
property. Nonetheless, in this paper we prove a result, Theorem \ref{maintheorem}, which we feel
is in fact closest in spirit to the classical theorem of Bernstein.  To precisely state this theorem, we first introduce some definitions, and present
some results of independent interest which are fundamental in the
proof of the main theorem.

To study $H$-minimal surfaces, we form a set of coordinates which
are adapted to the horizontal Gauss map on portions of the surface
which are graphs over a domain in the $xy$-plane.  In this sense,
they are analogous to isothermal coordinates on minimal surfaces
in $\R^3$.  Roughly, for an $H$-minimal surface $S$ written as a
graph $t=h(x,y)$ for a $C^2$ function $h$ over a domain $\Om$ of
the $xy$-plane where $S$ has empty characteristic locus, we
consider the unit horizontal Gauss map \eqref{hgm0} of this
surface. Setting $p=X_1  (t - h)$ and $q= X_2  (t - h)$, we have
\[
\nuX\ =\  \frac{p}{\sqrt{p^2+q^2}} \; X_1\ +\ \frac{q}{\sqrt{p^2+q^2}}
\; X_2\ .
\]

As the $H$-minimal surface in question is a graph over the
$xy$-plane, we have $p = p(x,y,h(x,y))$, $q=q(x,y,h(x,y))$, and
thus the horizontal Gauss map $\nuX(x,y,t)$ on $S$ depends in fact
only on the variables $(x,y)\in \Om$. By abuse of notation we
write $\nuX(x,y) = \nuX(x,y,h(x,y))$, and identify this function
$\nuX : \Om \to H\mathbb H^1$ with the unit vector field over
$\Om$
\begin{equation}\label{identification}
 (x,y)\ \longrightarrow\ \til(x,y)\ \overset{def}{=}\ \left(\frac{p(x,y)}{\sqrt{p(x,y)^2 + q(x,y)^2}}\ ,\
\frac{q(x,y)}{\sqrt{p(x,y)^2 + q(x,y)^2}}\right)\ .
\end{equation}

Henceforth in this paper, we routinely identify $z = x + i y$ with
the point $(x,y)$. Also, for $\zeta = (\zeta_1, \zeta_2)\in
\mathbb R^2$, we let $\zeta^\perp = (\zeta_2 , - \zeta_1)$.

\medskip

\begin{Def}
Given a point $z\in \Om\subset \mathbb R^2$, a \emph{seed curve
based at $z$} of the $H$-minimal surface $S$  is defined to be the
integral curve of the vector field $\til$ with initial point $z$.
Denoting such a seed curve by
  $\gamma_{z}(s)$, we then have
  \begin{equation}\label{seeddef0}
\gamma'_z(s)\ =\
  \til(\gamma_z(s))\ ,\quad\quad\quad \gamma_z(0)\ =\ z\ .
  \end{equation}
  If the base point $z$ is understood or irrelevant, we simply denote the seed curve by $\gamma(s)$.
We will indicate by $\tilL$ the integral curve of $\til^\perp$
starting at the point $z$.
\end{Def}

\medskip

We note explicitly that, thanks to $|\til| = 1$, a seed curve is
always parameterized by arc-length. Using
$\{\tilde{\mathcal{L}}_z, \gamma_z\}$ as our coordinate curves, we
obtain a new local parameterization of the $xy$-plane $F : \R^2
\ra \R^2$, given by
\begin{equation}\label{newpar}
(s,r)\ \ra\ (x(s,r), y(s,r))\ \overset{def}{=}\ F(s,r)\ =\
\gamma(s)\ +\ r\;\tilp(\gamma(s))\ .
\end{equation}

Keeping in mind that, from \eqref{seeddef0}, $\tilp(\gamma(s))=
\gamma'(s)^\perp =  (\gamma_2'(s), - \gamma_1'(s))$, we have
\begin{equation}\label{F0}
F(s,r)\ =\ (\gamma_1(s)\ +\ r\ \gamma_2'(s)\ ,\ \gamma_2(s)\ -\ r\ \gamma_1'(s))\ .
\end{equation}

As we will see in Section \ref{S:representation}, $F(s,r)$ defines a local diffeomorphism over a region of the $(s,r)$-plane, up to a certain curve $\mathcal C_\gamma$, see Lemma \ref{L:sl}.

Our first result is a basic representation theorem for $H$-minimal
surfaces with empty characteristic locus which are graphs  over a
portion of the $xy$-plane.

\medskip

\begin{MThm}\label{mrep}
A patch of a $C^k$ surface $S \subset \Hy^1$ of the type
\[
S\ =\ \{(x,y,t)\in \mathbb H^1 \mid (x,y)\in \Om\ ,\ t = h(x,y)\}\ ,
\]
where $h : \Om \to \R$ is a $C^k$ function over a domain $\Om$ in
the $xy$-plane, and with empty characteristic locus over $\Om$, is
an $H$-minimal surface if and only if for every $g = (z,t) \in S$,
there exists an open neighborhood of $g$ on $S$ that can be
parameterized by
\begin{equation}\label{para}
(s,r)\ \to\ (\gamma_1(s)+r\gamma_2'(s),\gamma_2(s)-r\gamma_1'(s),
h(s,r))\ ,
\end{equation}
where $h(s,r)$ is given by
\begin{equation}\label{para2}
h(s,r)\ =\ h_0(s)\ -\ \frac{r}{2} <\gamma(s), \gamma'(s)>\
\end{equation}
with
\[\gamma(s) \in C^{k+1},\;  h_0(s) \in C^k\ .
\]
Thus, to specify such a patch of smooth $H$-minimal surface, one must specify a single curve in $\Hy^1$ determined by a seed curve $\gamma(s)$, parameterized by arc-length, and an initial height function $h_0(s)$.
\end{MThm}

\medskip

One consequence of the representation in Theorem A is that $H$-minimal surfaces are, in fact, ruled surfaces.  In particular, for fixed $s$, the straight line
 $(\gamma_1(s)+r\gamma_2'(s),\gamma_2(s)-r\gamma_1'(s))$ in the $(s,r)$-plane lifts to a straight line on $S$ which is a geodesic in $\mathbb{H}^1$, see Corollaries \ref{C:geo} and \ref{C:geo2}.
  We stress that Theorem A is useful in both the
study of known examples as well as in the construction of new $H$-minimal surfaces.  Indeed, in Sections \ref{S:return} and \ref{S:characterization} we demonstrate the construction of new
types of surfaces, see Examples \ref{E:optreg}, \ref{E:optreg2}, \ref{E:cyl}, and Example \ref{E:gencurve}. We also stress that the
counterexample given above was discovered in a similar way.

Seed curves associated to $H$-minimal surfaces are our fundamental
objects of study.  With this in mind, we introduce the following definition.

\medskip

\begin{Def}
If no portion of a $C^2$, complete, connected $H$-minimal surface
  can be written as a graph over the $xy$-plane, we say that $S$ has
 \emph{trivial seed curve}.  Otherwise, $S$ has a \emph{non-trivial seed curve}.
\end{Def}

\medskip

In this definition and subsequent theorems, we will use the term {\em complete surface} to denote a metrically complete open surface in the Heisenberg group.

We show in Lemma \ref{trivial} that if $S$ has trivial seed curve,
then $S$ must be a vertical plane as in \eqref{conj}. However, to
completely describe an $H$-minimal surface we need to suitably
extend to notion of seed curve. A \emph{generalized seed curve} is
a collection of seed curves and associated height functions
together with patching data which define a single curve in
$\Hy^1$, see Definition \ref{D:genseed}.

The next theorem shows that generalized seed curves completely
determine $H$-minimal surfaces.

\medskip

\begin{MThm}\label{cleanup}Let $S\subset \mathbb H^1$ be a $C^2$, complete, connected $H$-minimal surface.
  Then, either $S$ is a vertical plane, or $S$ is determined by a generalized
  seed curve.
\end{MThm}

\medskip

In the following definition the signed curvature $\kappa(s)$ of a
seed curve is that given in \eqref{kappa} of Definition
\ref{D:seed}.

\medskip

\begin{Def}\label{D:introcc}
An $H$-minimal surface is said to have \emph{constant curvature} if either it
  has trivial seed curve (in which case $S$ is a vertical plane), or
  if the signed curvature $\kappa(s)$ of each seed curve which is part
  of the generalized seed curve defining the $H$-minimal surface is constant.
\end{Def}

\medskip
We emphasize two important points concerning this definition.
First, that the assumption of constant curvature does not imply
that the curve in $\mathbb{H}^1$ defined by the generalized seed
curve, $\Gamma_i = \{(\gamma^i(s),h_0^i(s))\}$, has constant
curvature, merely that the seed curves, $\gamma^i$, do.  Indeed,
Example \ref{E:gencurve} explicitly shows this:  each of its seed
curves are straight lines (constant curvature zero) but the lifted
curves, $(\gamma^i,h_0^i)$ do not have constant curvature. Second,
we point out that the individual seed curves, $\gamma^i$, for an
$H$-minimal surface of constant curvature may have different
signed curvatures for different $i$.  A method to construct
examples of such behavior is described after Definition
\ref{D:constcurv}.
\medskip

We can now state our main result.

\medskip

\begin{MThm}[\textbf{of Bernstein type}]\label{maintheorem}  Let $S$ be a $C^2$ connected $H$-minimal surface
  which is a graph over some plane $P$, then $S$ has constant
  curvature according to Definition \ref{D:introcc}.
\end{MThm}

\medskip

When the complete $H$-minimal surface $S$ fails to be a graph, then it
need not have constant curvature. An example is given by the
sub-Riemannian catenoids in Examples \eqref{E:catenoid},
\eqref{E:catenoid2}. We stress that in the Theorem C we have made no
assumption concerning the characteristic locus of $S$. Theorem C
points out a rigidity of the seed curve under the assumption that the
$H$-minimal surface is a graph:  it must be composed of circles or
lines.  In the case where all the seed curves in the generalized seed
curve are circles, $S$ may or may not have empty characteristic locus.
However, in the case where at least one seed curve in the generalized
seed curve is a straight line, we see that by (2) in Remark
\ref{R:rems}, $S$ always has non-empty characteristic locus. We stress,
however, that there are many
graph-like $H$-minimal surfaces given simply by specifying different initial
height functions. For instance, the $H$-minimal surface \eqref{ce} in
the counterexample above, is a graph over the $xt$-plane, it has empty
characteristic locus, see Example \ref{E:ce}, and as we show in
Example \ref{E:ceseed} its generalized seed curve consists of a single seed curve which is a circle. The plane $S =
\{(x,y,t)\in \mathbb H^1 \mid t = 0 \}$ is a graph over the
$xy$-plane, it has non-empty characteristic locus, and again its single seed
curve is a circle, see Example \ref{E:comp}. On the other hand, we
will construct a surface which is a graph
over the $xy$-plane, has non-empty characteristic locus and seed curve
a straight line, see Examples \ref{E:hyp} and \ref{E:hyp2}. In connection with this latter example,
if we further restrict our attention to graphs over
the $xy$-plane, then we can completely classify the possible $H$-minimal surfaces.

\medskip

\begin{MThm}  Suppose $S$ is a $C^2$ connected $H$-minimal
graph over the entire $xy$-plane, then:
\begin{enumerate}
\item Either $S$ has seed curve a circle, and $S$ is a plane of the form $ax+by + ct = d$ for some real numbers
  $a,b, c, d$, with $c\not= 0$, (with characteristic locus $\Sigma = \{(-2b/c,2a/c,d/c))\}$).
\item Or, $S$ has seed curve a straight line, and $S$ can be written as
\[\left (x+x_0,y+y_0,t_0-\frac{1}{2}ab(x^2-y^2)-\frac{1}{2}(b^2-a^2)xy+h_0(ax+by)+\frac{1}{2}x_0y-\frac{1}{2}xy_0
\right)\ ,
\] where $(x_0,y_0,t_0)$ is a point in $\R^3$ and $a,b$ are constants
so that $a^2+b^2=1$.
\end{enumerate}
\end{MThm}

\medskip
We note that, if we left translate so that the point $(x_0,y_0,t_0)$
is moved to the origin and further compose with an appropriate
rotation about the $t$-axis, we can write these surfaces as:
\[ \left (s+r,s-r,h_0(s)-\frac{sr}{2}\right)\]
In addition, we note that the both of these classes of surfaces were found in
\cite{Pauls:minimal}, and, after the appropriate changes to reflect
the difference in the representation of the Heiseneberg group, the second set of examples can be written as
\[
t\ =\ \alpha x^2-\frac{xy}{2}+f(x-\alpha y)
\]
for some real number $\alpha$, and some function $f$
(which is of course equivalent to the choice of $h_0(s)$).

We mention in closing that Theorem D has been first established in
the recent paper \cite{CHMY}, but with a completely different
approach from ours. Besides Theorem D, the paper \cite{CHMY}
contains several other interesting results and also deals with the
Heisenberg group from the broader perspective of three-dimensional
CR manifolds. The analysis in \cite{CHMY} is based on a detailed
study of the characteristic locus of $S$, whereas our approach is
centered on the concepts of seed curve and height function.  It is
perhaps worth mentioning that we received from the authors the
preprint \cite{CHMY} on June 23, 2003, when the present paper was
already undergoing an extensive revision with respect to its first
version, in which we attempted to establish the above mentioned
conjecture. At that time, in some email exchanges the authors of
\cite{CHMY} expressed to us their belief that our paper contained
a flaw. This was caused by a misstatement of the main result in
the abstract, which talked about minimal surfaces rather than
minimal graphs. This unfortunate oversight led the authors of
\cite{CHMY} to incorrectly believe that we had missed the catenoid
type example found by one of us (the latter was first described in
\cite{Pauls:minimal} and is reviewed in Example \ref{E:catenoid}).
As the above discussion of the counterexample to the conjecture
clarifies, the true reason for which our original attempt was
flawed was in fact completely different, and it was brought to
light by our discovery of the counterexample.

\medskip

\textbf{Acknowledgement:} We would like to thank the anonymous referee for his/her constructive criticism and for several comments which contributed to improve the presentation of the paper.

\vskip 0.6in


\section{\textbf{The horizontal Gauss map and the minimal surface equation}}\label{S:GM}

\vskip 0.2in

In this section, we describe two notions of fundamental importance
in the study of minimal surfaces in $\mathbb H^1$, the horizontal
Gauss map, and the $H$-mean curvature. To set the stage, we will
suppose $\Hy^1$ endowed with a Riemannian metric which makes
$\{X_1,X_2,T\}$ an orthonormal basis, where the vector fields
$X_1, X_2$ are given by \eqref{vf}, and as in \eqref{comm} we have
$[X_1,X_2] = T$, all other commutators being zero. We consider an
oriented hypersurface $S\subset \Hy^1$, and indicate with $\n$ the
Riemannian unit normal pointing outward. We recall $X =
\{X_1,X_2\}$ is an orthonormal basis of the horizontal subbundle
$H\Hy^1$ of the tangent bundle $T\Hy^1$. We recall that the
characteristic locus of $S$ is given by
\begin{equation}\label{char}
\Sigma\ =\ \{g\in S \mid H_g\Hy^1 \subset T_g S\}\ =\ \{g  \in S
\mid X_i(g) \in T_g S , i = 1 , 2\}\ .
\end{equation}

Given a function $f$ on $\Hy^1$, its \emph{horizontal gradient}
$\nabla_H f$ is the projection onto the subbundle $H\Hy^1$ of its
Riemannian gradient $\nabla f$ with respect to the orthonormal
basis $\{X_1,X_2,T\}$, i.e.,
\[
\nabla_H\ f\ =\ <\nabla f , X_1> X_1\ +\ <\nabla f , X_2> X_2\ =\
X_1 f\ X_1\ +\ X_2 f\ X_2\ .
\]

\medskip

\begin{Def}\label{D:HGauss}
We define the \emph{horizontal normal} $\Up : S \rightarrow H\Hy^1$, relative to the basis $X_1, X_2$ of the subbundle $H\mathbb H^1$,
by the formula
\[
\Up\ =\ <\n,X_1> X_1\ +\ <\n,X_2> X_2\ .
\]
The \emph{horizontal Gauss map} $\nuX : S \rightarrow \mathbb S^{1}$ is defined by
\[
\nuX\ =\ \frac{\Up}{|\Up|}\ ,
\]
whenever
$|\Up| \not= 0$.
\end{Def}

\medskip

Throughout the following discussion we will denote by $\mathcal S_H$ the natural domain of $\nuX$, i.e., the open set on $S$
\[
\mathcal S_H\ =\ \{g\in S\ \mid\ |\Up(g)| \not= 0\}\ .
\]

We note explicitly that $\Up$ is the projection $pr_H(\boldsymbol \nu)$ of the normal $\boldsymbol \n$ on the horizontal subbundle $H\Hy^1 \subset T \Hy^1$, whereas $\nuX$ is the normalized projection.
An obvious, yet important consequence of the definition is

\begin{equation}\label{normone}
|\nuX|^2\ \equiv \ 1\ ,\quad\quad\quad\quad \text{in}\quad \mathcal S_H\ .
\end{equation}

We next consider a function $u \in C^1(\mathcal O)$, where $\mathcal O$ is an open subset of $\mathbb H^1$ containing $S$.
The following notion, introduced in \cite{DGN}, is central to the development of sub-Riemannian calculus on a hypersurface.

\medskip

\begin{Def}\label{D:delta}
At every point $g\in \mathcal S_H$ the \emph{tangential horizontal gradient} of $u$ on $S$ with respect to the subbundle $H\mathbb H^1$ is defined as follows
\[
\delta_H u\ \overset{def}{=}\ \nabla_H u\ -\ <\nabla_H u,\nuX>\
\nuX\ .
\]
\end{Def}

In the next definition we introduce a sub-Riemannian analog of the classical notion of mean curvature.

\medskip

\begin{Def}\label{D:XMC}
The $H$-\emph{mean curvature} of $S$ at points of $\mathcal S_H$ is defined by
\[
\mathcal H\ =\ \sum_{j=1}^2\ \delta_{H,j}\ \nu_{H,j}\ .
\]
If $g_o\in \Sigma$ we let
\[
\mathcal H(g_o)\ =\ \underset{g\to g_o, g\in S_H}{\lim}\ \mathcal H(g)\ ,
\]
provided that such limit exists, finite or infinite. We do not define the $H$-mean curvature at those points $g_o\in \Sigma$ at which the limit does not exist.
\end{Def}

\medskip

\begin{Def}\label{D:MS}
A $C^2$ surface $S$ is called $H$-\emph{minimal} if its $H$-mean curvature $\mathcal H$ vanishes everywhere.
\end{Def}

\medskip

It is interesting to write the $H$-minimal surface equation in terms of a defining function for the hypersurface $S$. Suppose that
\begin{equation}\label{Sphi}
S\ =\ \left\{g \in \Hy^1 \mid \phi(g) = 0\right\}\ ,
\end{equation}
where $\phi : \Hy^1 \to \R$ is $C^2$, and there exists an open neighborhood $\mathcal O$ of $S$ such that
\[
|\nabla \phi(g)|\ \geq\ \alpha\ >\ 0\ \quad\quad\quad g \in \mathcal O\ .
\]

We will think of $S$ as the boundary of the open set $\mathcal U =
\{g\in \Hy^1 \mid \phi(g) < 0\}$, so that the Riemannian outer
unit normal
 to $S$ with respect to the orthonormal basis $\{X_1,X_2,T\}$
is given $\n = \nabla \phi/|\nabla \phi|$. One easily recognizes
that in such case the horizontal Gauss map is given by
\begin{equation}\label{GM}
\nuX\ =\ \frac{\nabla_H \phi}{|\nabla_H \phi|}\
,\quad\quad\quad\quad \text{in}\quad\mathcal S_H\ .
\end{equation}

Definition \ref{D:XMC}, combined with \eqref{normone}, gives in such case
\begin{equation}\label{MSE}
\mathcal  H\ =\ \sum_{j=1}^2\ X_j\ \nu_{H,j}\ =\ \sum_{j=1}^2\
X_j\ \left\{\frac{X_j\phi}{|\nabla_H \phi|}\right\}\ .
\end{equation}

Thereby, the notion of $H$-minimal surface introduced in \cite{DGN} coincides with that set forth in \cite{Pauls:minimal}, and the corresponding \emph{minimal surface equation} is given by
\begin{equation}\label{MSE2}
\sum_{j=1}^2\ X_j\ \left\{\frac{X_j\phi}{|\nabla_H \phi|}\right\}\
=\ 0\ .
\end{equation}

For a function $u : \Hy^1 \to \mathbb R$ we define the \emph{symmetrized horizontal Hessian} of $u$ at $g\in \Hy^1$ as the $2\times 2$ matrix with entries
\begin{equation}\label{hessian}
u_{,ij}\ \overset{def}{=}\ \frac{1}{2}\ \bigg\{X_iX_j u\ +\ X_j X_i u\bigg\}\ ,\quad\quad\quad\quad i , j = 1 , 2\ .
\end{equation}

Setting $\nabla^2_H u = (u_{,ij})$, we introduce the matrix
\begin{equation}\label{H}
\mathcal A\ =\ \begin{pmatrix}
\mathcal A_{11} & \mathcal A_{12}
\\
\mathcal A_{21} & \mathcal A_{22}
\end{pmatrix}
\ \overset{def}{=}\ \frac{1}{|\nabla_H \phi|^3}\ \left\{|\nabla_H
\phi|^2 \nabla^2_H \phi\ -\ \nabla^2_H \phi(\nabla_H \phi) \otimes
\nabla_H \phi\right\}\ .
\end{equation}

The next proposition follows by a direct computation based on \eqref{MSE} and we leave the details to the reader.

\medskip

\begin{Pro}\label{P:XMC}

At every point of $\mathcal S_H$ one has
\[
\mathcal H\ =\ tr\ \mathcal A\ .
\]
In particular, the $H$-minimal surface equation for $S$ can be re-written as
\begin{equation}\label{XMC2}
tr\ \mathcal A\ =\  0\ .
\end{equation}
\end{Pro}

\medskip

It is standard to verify that the eigenvalues of $\mathcal A$ are given by the equations
\begin{equation}\label{ev}
\lambda_1\ +\ \lambda_2\ =\ tr\ \mathcal A\ ,\quad\quad\quad \lambda_1\ \lambda_2\ =\ det\ \mathcal A\ .
\end{equation}

At non-characteristic points of $S$ we have $\nabla_H\phi \not=
0$. If we consider the $i$-th entry of the vector $\mathcal A
(\nabla_H \phi)$, from the second equality in \eqref{H} we obtain
\begin{align*}
& \frac{1}{|\nabla_H \phi|^3}\ \big\{|\nabla_H \phi|^2\
\phi_{,ij}\ X_j\phi\ -\ \phi_{,ik}\ X_k\phi\ X_j\phi\
X_j\phi\big\}
\\
& =\ \frac{1}{|\nabla_H \phi|^3}\ \big\{|\nabla_H \phi|^2\
\phi_{,ij}\ X_j\phi\ -\ |\nabla_H \phi|^2\ \phi_{,ik}\
X_k\phi\big\} = \ 0\ .
\end{align*}

This proves that at points where $\nabla_H \phi \not= 0$, we have
$\mathcal A (\nabla_H \phi) = 0$, so that $\nabla_H \phi$ is an
eigenvector of $\mathcal A$ with corresponding eigenvalue
$\lambda_1 = 0$.  From \eqref{ev} we conclude that the second
eigenvalue of $\mathcal A$ is given by $\lambda_2 = tr\ \mathcal
A$. One easily finds that an eigenvector for $\lambda_2$ is given
by $(\mathcal A_{12}, \mathcal A_{22})^T$. Thanks to $det\
\mathcal A = 0$, another eigenvector for $\lambda_2$ is given by
$(\mathcal A_{11}, \mathcal A_{21})^T$.

With an homage to classical notation, we now define
\begin{equation}\label{isothermal}
p\ =\ X_1 \phi\ ,\quad\quad q\ =\ X_2 \phi\ , \quad\quad  W\ =\
\sqrt{p^2+q^2}\ =\ |\nabla_H \phi|\ .
\end{equation}

Moreover, for convenience we let
\begin{equation}\label{isobar}
\pb\ =\ \frac{p}{W}\ ,\quad\quad\quad\quad \qb\ =\ \frac{q}{W}\ ,
\end{equation}
so that the minimal surface equation \eqref{MSE2} becomes $X_1 \pb  + X_2 \qb = 0.$
We also note the following useful equivalent form in terms of the functions $p$ and $q$
\begin{equation}\label{mse2}
X_1\bar p + X_2 \bar q\ =\ \frac{1}{W^3}\ \left\{q^2\ X_1p\ +\ p^2\ X_2q\ -\ p q\ (X_1q + X_2 p)\right\}\ =\ 0\ .
\end{equation}

We stress that if \eqref{MSE2} is satisfied, i.e., if $S$ is $H$-minimal, then from \eqref{XMC2} and \eqref{ev}, we conclude that also $\lambda_2 = 0$, so that $\mathcal A$
 has a double eigenvalue equal to zero with a single
eigenvector $\nabla_H \phi = p X_1 + q X_2$.

The expression of $\mathcal A$ with respect to the variables $p, q$ is

\begin{equation}\label{Hpq}
\mathcal A\ =\ \frac{1}{W^3}\ \begin{pmatrix}
q^2\ X_1p - p q  X_1q &  p^2 X_1 q - pq X_1 p
\\
q^2 X_2p\ -\ p q  X_2q &  p^2 X_2 q\ - pq X_2 p
\end{pmatrix}\ .
\end{equation}
With respect to these variables, the eigenvectors are $(p,q)^T$
(corresponding to eigenvalue $\lambda_1=0$) and $(q,-p)^T$
(corresponding to eigenvalue $\lambda_2=\mathcal{H}$).

Returning to Definition \ref{D:HGauss} we see that with respect to the variables $\bar p$ and $\bar q$ the horizontal Gauss map is given by
\begin{equation}\label{horgaussiso}
\nuX\ =\  \bar p\ X_1\ +\ \bar q\ X_2\ .
\end{equation}

We end this section with a useful observation.

\medskip

\begin{Lem}
If $S$ in an $H$-minimal surface in $\mathbb{H}^1$ then any left-translation of $S$ is an $H$-minimal surface.  Moreover, any rotation of
$S$ about the t-axis is also an $H$-minimal surface.
\end{Lem}

The proof of this is a computation of the effect of left-translation
and rotation about the $t$-axis on the horizontal Gauss map coupled with
the fact that the horizontal subbundle is preserved by both actions.
The details are left to the reader.

\vskip 0.6in


\section{\textbf{Some examples of $H$-minimal surfaces}}\label{S:examples}

\vskip 0.2in

In this section we analyze some interesting examples of
$H$-minimal surfaces. Most of the surfaces discussed in this
section, with the notable exception of Example \ref{E:ce}, are not
new, as they have already appeared in \cite{Pauls:minimal}.
However, the presentation in this section differs in spirit from
that in \cite{Pauls:minimal} since it is organized having in mind
the new ideas to be introduced in Section \ref{S:representation}.
In this sense, it should be viewed as preparatory to Section
\ref{S:return}, where such new ideas will be illustrated in
detail.

The reader should be aware that, throughout this section, and also
in the rest of the paper, when a surface $S$ is given as a graph
over a portion of the $xy$-plane, as in \eqref{S1}, we will
identify the horizontal Gauss map, which is given by
\eqref{horgaussiso}, with the two dimensional vector field $\til =
(\bar p, \bar q)$ defined in \eqref{identification}. Therefore, if
sometimes we say that $\nuX$ is equal to $ (\bar p, \bar q)$, we
mean that such equation has to be interpreted through the
identification $\nuX \cong \til$. As we explain in the opening of
Section \ref{S:representation} such identification is justified
and also significantly aids in the computations.

\medskip

\begin{Example}\label{E:charplane}
\emph{Consider the characteristic plane $\Pi = \{(x,y,t)\in \mathbb H^1 \mid  t = 0\}$. Using the defining function $\phi(x,y,t) = t$, and the equation \eqref{GM},
a simple computation shows that
\[
\til\ =\ \left(- \frac{y}{|z|}, \frac{x}{|z|}\right)\
,\quad\quad\quad\quad (x,y,t) \not= (0,0,0)\ .
\]
The horizontal Gauss map is not defined at the characteristic point $0 = (0,0,0)$ of $\Pi$. Nonetheless, one easily sees that
\[
\mathcal H\ =\ \partial_x \left(- \frac{y}{|z|}\right)\ +\ \partial_y
\left(\frac{x}{|z|}\right)\ \equiv\ 0\ ,\quad\quad\quad\text{for}\quad (x,y,t) \not= 0\ ,
\]
and therefore it is possible to define $\mathcal H(0) = 0$ as a limit. This example is a special case of the following more general fact. Let
\begin{equation}\label{charplane}
S\ =\ \{(x,y,t)\in \mathbb H^1 \mid a x + b y + ct = d \}\ ,
\end{equation}
be an arbitrary plane in $\mathbb H^1$. When $c = 0$ and $a^2 + b^2 \not= 0$, we obtain a vertical plane. We emphasize that such vertical planes have empty characteristic locus, with constant
 horizontal Gauss map coinciding with the Riemannian unit normal.
 Recalling \eqref{identification} we thus obtain
\[
\til\ =\ \left(\frac{a}{\sqrt{a^2 + b^2}}, \frac{b}{\sqrt{a^2 +
b^2}}\right)\ .
\]
When $c\not= 0$, the characteristic locus is composed of one single point $\Sigma = \{(- 2b/c, 2a/c, d/c)\}$, and for $(x,y,t) \not\in \Sigma$ we have
\begin{equation}\label{nuXplanes}
\til\ =\  \left(\frac{a - \frac{c y}{2}}{\sqrt{a^2 + b^2 +
\frac{c^2}{4} |z|^2 + bcx - ac y}}, \frac{b + \frac{c
x}{2}}{\sqrt{a^2 + b^2 + \frac{c^2}{4} |z|^2 + bcx - ac
y}}\right)\ .
\end{equation}
In either case, it is readily recognized from \eqref{nuXplanes} that
\[
\mathcal H\ =\ div\ \til\ \equiv\ 0\ ,
\]
and therefore $S$ is a complete embedded $H$-minimal surface. Incidentally, we notice the case $c\not= 0$ also follows from the analysis of the characteristic plane $\Pi$, and
by the left-translation invariance of the notion of $H$-mean curvature. All one needs to observe is that by left-translating $\Pi$ with respect to the point $g_o = (- 2b/c, 2a/c,d/c)$,
one obtains the plane \eqref{charplane}.}
\end{Example}

\medskip

\begin{Example}\label{E:hyp}
\emph{We next analyze the surface given by
\begin{equation}\label{S}
S\ =\ \left\{(x,y,t) \in \mathbb H^1 \mid t = \frac{x y}{2} \right\}\ ,
\end{equation}
with defining function $\phi(x,y,t) = t - \frac{xy}{2}$. One has
\[
X_1\phi\ =\ -\ y\ ,\quad\quad\quad\quad X_2\phi\ =\ 0\ ,
\]
and therefore the characteristic locus of $S$ is given by the $x$-coordinate axis, i.e.,  $\Sigma = \{(x,y,t) \in \mathbb H^1 \mid y = t = 0\}$. The horizontal Gauss map is well-defined
 in either of the two regions $\mathcal S_H^+ = S \cap \{(x,y,t)\in \mathbb H^1 \mid y > 0\}$, $\mathcal S_H^- = S \cap \{(x,y,t)\in \mathbb H^1 \mid y <
 0\}$. Furthermore, $\til$ is constant on each region, and given by
\begin{equation}\label{nuXhyp}
\til\ =\ \left(- \frac{y}{|y|} , 0\right)\ =\ ( \mp 1 , 0)\ ,
\quad\quad \text{on} \quad \mathcal S_H^\pm\ .
\end{equation}
It is immediately seen from \eqref{nuXhyp} that
\[
\mathcal H\ =\ \partial_x \left(- \frac{y}{|y|}\right)\  \equiv\ 0\ ,\quad\quad\quad\text{for}\quad (x,y,t)\in \mathcal S_H = \mathcal S_H^+ \cup \mathcal S_H^- \ .
\]
At points of the line $\Sigma$, the $H$-mean curvature can be defined
as a limit. We thus conclude that $\mathcal H \equiv 0$ on $S$, and therefore
$S$ is an entire $H$-minimal graph over the $xy$-plane.  Figure
\ref{eg3.2} shows a picture of this surface.  The black line is the
characteristic locus.
}
\begin{figure}
\epsfig{figure=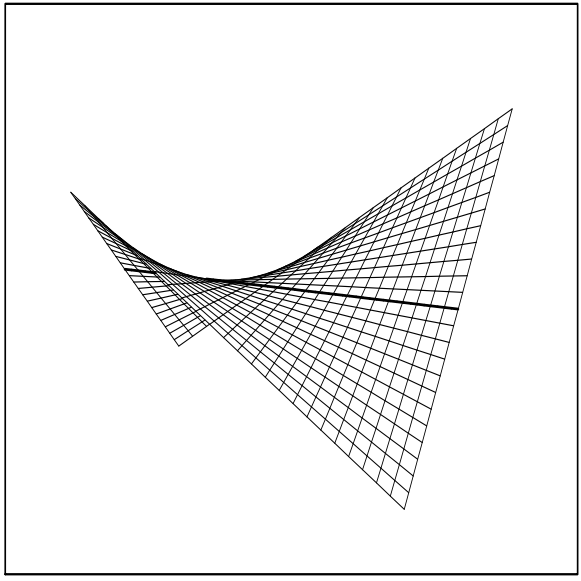,width=3in,height=3in,angle=270}
\caption{The $H$-minimal surface $\{t = \frac{xy}{2}\}$}\label{eg3.2}
\end{figure}

\end{Example}

\medskip

Before discussing the next example, we recall a useful result from \cite{DGN}.

\medskip

\begin{Pro}\label{P:tphi}
Let $S$ be a hypersurface in $\mathbb H^n$ which is the level set of a defining function of the type
\begin{equation}\label{specialphi}
\phi(z,t)\ =\ u\left(\frac{|z|^2}{4}\right)\ -\ t\ ,
\end{equation}
for some $C^2$ function $u:[0,\infty) \to \mathbb R$.
For every point point $g = (z,t)\in S$, such that $z\not= 0$, the $H$-mean curvature at $g$ is given by
\begin{equation}\label{specialHX}
\mathcal H\ =\ \frac{2\ s\ u''(s)\ +\ (Q -3)\ u'(s)\ (1 + u'(s)^2)}{2\ \sqrt s\ (1 + u'(s)^2)^{3/2}}\ ,\quad\quad s = \frac{|z|^2}{4}\ .
\end{equation}
\end{Pro}

\medskip

Using Proposition \ref{P:tphi} it is possible to compute an interesting family of complete (globally defined) rotationally invariant $H$-minimal surfaces
 in $\mathbb H^1$. Let us assume that $S$ be a $C^2$ hypersuface in $\Hn$ described by a defining function as in \eqref{specialphi}. Using
\eqref{specialHX} for its $H$-mean curvature, and imposing that $\mathcal H \equiv 0$, we find for $s>0$
\begin{equation}\label{cate1}
2\ s\ u''(s)\ +\ (Q -3)\ u'(s)\ (1 + u'(s)^2)\ =\ 0\ .
\end{equation}

To solve \eqref{cate1} we multiply it by $u'$, and set $y = 1 + (u')^2$, which gives $y' = 2 u' u''$, obtaining
\begin{equation*}
y'(s)\ -\ \frac{Q-3}{s}\ y(s)\ =\ -\ \frac{Q-3}{s}\ y^2(s)\ .
\end{equation*}

This is a Bernoulli equation which is solved by the substitution $v = y^{-1}$, which yields
\[
v(s)\ =\ \frac{v_o}{s^{Q-3}}\ -\ 1\ .
\]

Returning to the dependent variable $u$, we finally obtain
\begin{equation}\label{cate2}
u(s)\ =\ u_o\ \pm\ \int \frac{1}{\sqrt{a s^{Q-3} - 1}}\ ds\ ,
\end{equation}
where $a>0$ and $u_o\in \mathbb R$ are arbitrary. If we specialize \eqref{cate2} to the case of $\mathbb H^1$, for which $Q = 4$, we find
\begin{equation}\label{cate3}
u(s)\ =\ u_o\ \pm\ \frac{2}{a}\ \sqrt{a s - 1}\ .
\end{equation}

Substituting $s = |z|^2/4$ in \eqref{cate3}, we obtain the following interesting conclusion.

\medskip

\begin{Example}\label{E:catenoid}
\emph{The surfaces
\[
S\ =\ \left\{(x,y,t)\in \mathbb H^1 \mid (t - u_o)^2\ =\ \frac{4}{a^2}\left(\frac{a}{4} |z|^2 - 1\right)\right\}
\]
are real-analytic complete $H$-minimal surfaces in $\mathbb H^1$, with empty characteristic locus, which are not graphs on any of the three coordinate planes. They are the sub-Riemannian analogue of the classical catenoids.
Thanks to \eqref{specialHX}, we see that $S$ satisfies the equation $\mathcal H \equiv 0$ for every $(z,t)\in \mathbb H^1$ such that $z\not= 0$. Since $S$ does not intersect the $t$-axis, we conclude that the $H$-mean curvature vanishes everywhere on $S$.
The fact that $S$ has empty characteristic locus easily follows by the rotational symmetry of $S$ about the group center, the $t$-axis,
and by the fact that a surface with such property can only have characteristic points at the intersection with the $t$-axis. Since such
intersection is empty for $S$, the conclusion follows. These surfaces where first discovered by the second named author, see Section 4.1.2 in \cite{Pauls:minimal}.
The choice $u_o = 0$, and $a = 2$, yields the surface}
\[
S\ =\ \left\{(x,y,t)\in \mathbb H^1 \mid t^2\ =\ \frac{|z|^2}{2} - 1\right\}\ ,
\]
\emph{which is represented in figure \ref{eg3.4}.
Although we have already computed the $H$-mean curvature using Proposition \ref{P:tphi}, it is nonetheless useful, for later purposes, to compute the horizontal Gauss map for such $S$.
Restricting attention to that portion of $S$ that lies above the plane $\{(x,y,t) \in \mathbb H^1\mid t=0\}$, we can use the defining function
\[
\phi(x,y,t)\ =\ t\ -\ \sqrt{\frac{|z|^2}{2} - 1}\ ,
\]
obtaining for projected horizontal Gauss map}
\begin{equation}\label{nuXcat}
\til\ =\ -\ \sqrt{2}\ \left(\frac{x + y \sqrt{\frac{|z|^2}{2} -
1}}{|z|^2}\ ,\ \frac{y - x \sqrt{\frac{|z|^2}{2} -
1}}{|z|^2}\right)\ .
\end{equation}
\begin{figure}
\epsfig{figure=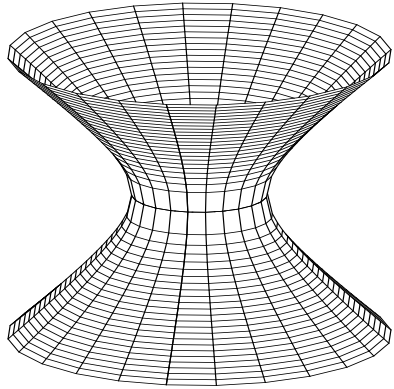,width=2in,height=2in}
\caption{The catenoid type surface $\{t^2=\frac{|z|^2}{2}-1\}$}\label{eg3.4}
\end{figure}
\end{Example}

\medskip

\begin{Example}\label{E:ce}
\emph{The following example is quite striking, since it provides a $C^\infty$ (in fact, real analytic) $H$-minimal surface which is an entire graph over the coordinate $xt$-plane, has empty characteristic locus, and which fails to be a vertical plane. Such example shows that the conjecture in the introduction fails to be true. Consider the surface
\begin{equation}\label{ce1}
S\ =\ \{(x,y,t)\in \mathbb H^1 \mid y = - x \tan(\tanh(t))\ ,\ (x,t)\in \mathbb R^2\}\ .
\end{equation}
Clearly, \eqref{ce1} gives $S$ as an entire real analytic graph over
the $xt$-plane.  Figure \ref{eg3.5} shows the surface as a graph over
the $xt$-plane.  Let us verify first that $S$ has empty characteristic locus $\Sigma$. Using the defining function
\[
\phi(x,y,t)\ =\ y - f(x,t)\ ,\quad\quad\quad \text{where}\quad f(x,t) \ =\ -\ x\ \tan(\tanh(t))\ ,
\]
we form
\begin{equation}\label{4}
X_1\phi\ =\ - f_x + \frac{y}{2} f_t\ ,\quad\quad \quad X_2\phi\ =\  1 - \frac{x}{2} f_t\ .
\end{equation}
From the second equation in \eqref{4} it is clear that if $x=0$,
then $X_2\phi = 1$ and so $\Sigma$ has empty intersection with the
plane $\{(x,y,t)\in \mathbb H^1 \mid x=0\}$. Suppose then that
$x\not= 0$. In such case we see that $X_2\phi = 0$ if and only if
\[
f_t\ =\ \frac{2}{x}\ .
\]
Now a simple computation gives
\[
f_t\ =\ -\ \frac{x}{\cosh(2t)}\ ,
\]
and therefore
\begin{equation}\label{5}
X_2\phi = 0\ \quad\quad \text{if and only if} \quad\quad -\ \frac{x^2}{2}\ =\  \cosh(2t)\ .
\end{equation}
Since the latter equation has no solution for $x\not=0$, it is clear from \eqref{5} that $X_2\phi \not=0$, and thus $\Sigma = \varnothing$. We next want to show that $S$ is $H$-minimal. In this respect, it is easier to look at those patches of $S$ which can be written as a graph over the $xy$-plane. One easily recognizes that, on either side of the hyperplane $\{(x,y,t) \in \mathbb H^1 \mid x = 0\}$, we can write the surface as follows
\begin{equation}\label{hor}
t\ =\ h(x,y)\ \overset{def}{=}\ -\  \tanh^{-1}\left(\tan^{-1}\left(\frac{y}{x}\right)\right)\ .
\end{equation}
In this form, we take as defining function
\[
\psi(x,y,t)\ =\ t\ -\ h(x,y)\ ,
\]
obtaining
\[
X_1\psi\ =\ -\ h_x\ -\ \frac{y}{2}\ ,\quad\quad\quad X_2\psi\ =\ -\ h_y\ +\ \frac{x}{2}\ ,
\]
\[
|\nabla_H \psi|^2\ =\ |\nabla h|^2\ +\ \frac{|z|^2}{4}\ +\ y\ h_x\
-\ x\ h_y\ =\  =\ |\nabla h|^2\ +\ \frac{r^2}{4}\ +\ <\nabla h,
z^\perp>\ ,
\]
where we have let $r = |z| = \sqrt{x^2 + y^2}$. From these
formulas we obtain
\begin{align}\label{HXce}
\mathcal H\ & =\ div\ \til\ =\ -\ div\ \left(|\nabla_H \psi|^{-1}\
\left( h_x + \frac{y}{2}, h_y - \frac{x}{2}\right)\right)
\\
& =\ |\nabla_H \psi|^{-3}\ \left\{|\nabla_H \psi|^2\ \Delta h\ -\
\frac{1}{2} <\nabla(|\nabla_H \psi|^2), \left( h_x + \frac{y}{2},
h_y - \frac{x}{2}\right)>\right\}\ . \notag
\end{align}
In the following computation it will be convenient to write $h = - \tanh^{-1}\ \alpha$, with $\alpha(x,y) = \tan^{-1}(y/x)$. We see that
\[
\nabla \alpha\ =\ -\ \frac{z^\perp}{r^2}\ ,\quad\quad |\nabla \alpha|^2\ =\ \frac{1}{r^2}\ ,\quad\quad \Delta \alpha\ =\ 0\ ,
\]
\[
\nabla h\ =\ \frac{z^\perp}{r^2 (1 - \alpha^2)}\ ,\quad\quad |\nabla h|^2\ =\ \frac{1}{r^2 (1 - \alpha^2)^2}\ ,\quad\quad \Delta h\ =\ -\ \frac{2 \alpha}{r^2 (1 - \alpha^2)^2}\ ,
\]
\[
|\nabla_H \psi|^2\ =\ \frac{1}{r^2 (1 - \alpha^2)^2}\ +\
\frac{r^2}{4}\ +\ \frac{1}{1 - \alpha^2}\ ,
\]
\[
\nabla(|\nabla_H \psi|^2)\ =\ -\ \frac{2}{r^4 (1 - \alpha^2)^2}\
z\ -\ \frac{4 \alpha}{r^4 (1 - \alpha^2)^3}\ z^\perp\ +\
\frac{z}{2}\ -\ \frac{2 \alpha}{r^2 (1 - \alpha^2)^2}\ z^\perp\ ,
\]
and therefore,
\[
<\nabla(|\nabla_H \psi|^2),z^\perp>\ =\ -\ \frac{4 \alpha}{r^2 (1
- \alpha^2)^3}\ -\ \frac{2 \alpha}{(1 - \alpha^2)^2}\ .
\]
From these formulas we conclude
\begin{align}\label{HXce2}
& |\nabla_H \psi|^2\ \Delta h\ -\ \frac{1}{2} <\nabla(|\nabla_H
\psi|^2), \left( h_x + \frac{y}{2}, h_y - \frac{x}{2}\right)>
\\
& =\ |\nabla_H \psi|^2\ \Delta h\ -\ \left[\frac{1}{2 r^2 (1 -
\alpha^2)}\ + \frac{1}{4}\right] <\nabla(|\nabla_H
\psi|^2),z^\perp>
\notag\\
 & =\ -\ \frac{2 \alpha}{r^2 (1 - \alpha^2)^2}\ \left[\frac{1}{r^2 (1 - \alpha^2)^2}\ +\ \frac{r^2}{4}\ +\ \frac{1}{1 - \alpha^2}\right]
\notag\\
& +\ \left[\frac{1}{2 r^2 (1 - \alpha^2)}\ + \frac{1}{4}\right]\ \left[\frac{4 \alpha}{r^2 (1 - \alpha^2)^3}\ +\ \frac{2 \alpha}{(1 - \alpha^2)^2}\right]
\notag\\
& =\ 0\ .
\notag
\end{align}
Inserting \eqref{HXce2} into \eqref{HXce}, we conclude that $\mathcal H \equiv 0$ on either of the two regions  $S\cap \{(x,y,t)\in \mathbb H^1 \mid x>0\}$,
$S\cap \{(x,y,t)\in \mathbb H^1 \mid x < 0\}$. Since, as we saw from the representation \eqref{ce}, the characteristic locus of $S$ is empty, and therefore
 the horizontal Gauss map is defined everywhere, we can extend $\mathcal H$ by continuity to all of $S$, thus concluding that such surface is globally $H$-minimal.
\\
 We conclude the analysis of this example by writing explicitly the projection onto the $xy$-plane of the horizontal Gauss map when $S$ is expressed in the form \eqref{hor}
\begin{equation}\label{nuXce}
\til\ =\ \left( -\ \frac{\left(\frac{1}{r^2 ( 1 - \alpha^2)} +
\frac{1}{2}\right)\ y}{\sqrt{\frac{1}{r^2 (1 - \alpha^2)^2} +
\frac{r^2}{4} + \frac{1}{1 - \alpha^2}}}\ ,\
\frac{\left(\frac{1}{r^2 ( 1 - \alpha^2)} + \frac{1}{2}\right)\
x}{\sqrt{\frac{1}{r^2 (1 - \alpha^2)^2} + \frac{r^2}{4} +
\frac{1}{1 - \alpha^2}}}\right)\ .
\end{equation}
}

\begin{figure}
\epsfig{figure=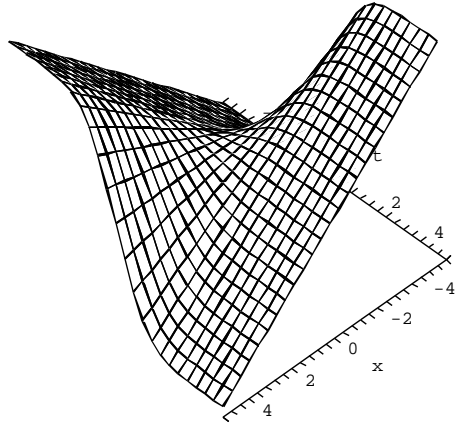,width=3in,height=3in}
\caption{The $H$-minimal surface $\{y=-x\tan(\tanh t)\}$}\label{eg3.5}
\end{figure}
\end{Example}

\medskip

\vskip 0.4in

We will return to the above examples in Section \ref{S:return}. To close this section we remark that, although all the above examples of $H$-minimal surfaces are not just $C^2$,
 but in fact real-analytic, due to the presence of characteristic points it is possible for a surface of constant $H$-mean curvature to have only a fixed amount of regularity.
As an instance of this negative phenomenon, we recall the following (partial) solution of the isoperimetric problem in $\Hn$ established in \cite{DGN}.
Consider the collection of sets $\mathcal E = \{E\subset \Hn\,|\, E \quad\text{satisfies $(i)-(ii)$}\}$,
where
\begin{itemize}
\item[(i)] $|E\cap \mathbb H^n_{+}|\ =\ |E \cap \mathbb H^n_{-}|$
; \item[(ii)] There exist $R>0$, and $C^1$ functions $u,v:
[0,R]\to \mathbb R$, satisfying $u(R) = v(R) = 0$, such that
 $\partial E \cap \mathbb H^n_{+}\ =\ \{(z,t)\,|\,t=u(|z|)\}$ and
 $\partial E \cap \mathbb H^n_{-}\ =\ \{(z,t)\,|\,t=v(|z|)\}$ .
\end{itemize}

\medskip

\begin{Thm}\label{isoprofile}
Given $V>0$, the variational problem
\[
min\{P_H(E;\Hn)\,|\, E \in \mathcal E,\ |E|=V\}
\]
has a unique solution $E_o\in\mathcal E$. Furthermore,
$\partial E_o$ is given explicitly as the graph
\begin{align}\label{iso}
t\ =\ & \pm\ \left\{\frac{1}{4}|z|\sqrt{R^2-|z|^2}\ -\ \frac{R^2}{4}\,\tan^{-1}\left(\frac{|z|}{\sqrt{R^2-|z|^2}}\right)\ +\ \frac{\pi R^2}{8}\right\}\ , \quad |z|\ \leq\ R\ .
\end{align}
The sign $\pm$ depends on whether one considers $\partial E_o
\cap \Hn_+$, or $\partial E_o\cap \Hn_-$. Finally, the set $E_o$ is only of class
$C^2$ near its two characteristic points $\left(0,\pm \frac{\pi R^2}{8}\right)$, it is $C^\infty$ away from them, and the $C^2$ compact hypersurface without boundary $S = \p E_o$ has positive
constant $H$-mean curvature given by
\[
\mathcal H\ =\  \frac{Q-2}{R}\ .
\]
\end{Thm}

The isoperimetric sets \eqref{iso} are surfaces of positive constant $H$-mean curvature. It is natural to ask whether it is possible to construct surfaces of zero $H$-mean curvature which
 are $C^2$, but not better than that. In Section \ref{S:return} we will show how to construct such $H$-minimal surfaces, see Examples \ref{E:optreg} and \ref{E:optreg2}.

\vskip 0.6in


\section{\textbf{A representation theorem for $H$-minimal graphs over the
  $xy$-plane}}\label{S:representation}

\vskip 0.2in

In this section, we consider an a priori special case, where we have an open set $\Om \subset \mathbb R^2$, and a $C^2$ function $h : \Om \to \R$ such that
\begin{equation}\label{S1}
S\ =\ \{(x,y,t)\in \mathbb H^1\mid (x,y)\in \Om\ ,\ t\ =\ h(x,y)\}\ .
\end{equation}

With $\phi(x,y,t)= t - h(x,y)$, we will assume that
$S = \p \mathcal U \cap \{g\in \mathbb H^1\mid \phi(g) = 0\}$,
where
$\mathcal U = \{g = (x,y,t) \in \mathbb H^1 \mid (x,y)\in \Om\ ,\ t < h(x,y)\}$,
so that the outward Riemannian unit normal to $S$ is given by
\begin{align}\label{nu}
\boldsymbol{\nu}\ =\ & -\ \frac{h_x + \frac{y}{2}}{\sqrt{1 +
    (h_x + \frac{y}{2})^2+(h_y-\frac{x}{2})^2}}\ {X_1}\ -\
\frac{h_y - \frac{x}{2}}{\sqrt{1 +(h_x + \frac{y}{2})^2+(h_y - \frac{x}{2})^2}}\ {X_2}
\\
& +\ \frac{1}{\sqrt{1 +(h_x +\frac{y}{2})^2+(h_y - \frac{x}{2})^2
}}\ T\ .\notag
\end{align}

We recall that $\mathbb H^1$ has been endowed with a
left-invariant Riemannian metric with respect to which
$\{{X_1},{X_2},T\}$ constitute an orthonormal basis. We note that
if we let $c(x,y)=(h_x+\frac{y}{2})^2+(h_y-\frac{x}{2})^2$, then
the characteristic locus of $S$ is given by $c^{-1}(0)$.  As $c$
is a continuous function, we note that the characteristic locus is
a closed set on the surface.

Re-writing the functions defined in \eqref{isothermal} for the present situation, we have
\begin{equation}\label{type1}
\begin{split}
p&\ =\  -\ \left(h_x\ +\ \frac{y}{2}\right)\ ,\quad\quad\quad\quad q\ =\ -\ \left(h_y\ -\ \frac{x}{2}\right)
\\
W \ &=\ \sqrt{p^2+q^2}\ =\ \sqrt{|\nabla h|^2 + \frac{|z|^2}{4} + <\nabla h, z^\perp>}\ ,\\
\end{split}
\end{equation}
where we have denoted $z^\perp = y - i x$.
Moreover, since from \eqref{type1} we see that $\bar p = p/W$, $\bar q = q/W$ do not depend on $t$, the
minimal surface equation \eqref{mse2} presently reads
\begin{equation}\label{msdiv}
div \; \til\ =\ 0\ ,
\end{equation}
where $div$ is the usual divergence operator on $\R^2$ and $\til$
is the projection of the horizontal Gauss map defined in
\eqref{identification}. Thus, minimal surfaces of this
non-parametric type correspond to divergence free unit vector
fields on $\R^2$. This observation leads us to introduce the
following definition.

\medskip

\begin{Def}\label{D:seedcurve}
Given a point $z\in \Om\subset \mathbb R^2$, a \emph{seed curve
based at $z$} of the $H$-minimal surface $S$ of the type
\eqref{S1} is defined to be the integral curve of the vector field
$\til$ with initial point $z$. Denoting such a seed curve by
  $\gamma_{z}(s)$, we then have
  \begin{equation}\label{seeddef}
\gamma_z'(s)\ =\
  \til(\gamma_z(s))\ ,\quad\quad\quad \gamma_z(0)\ =\ z\ .
  \end{equation}
  We note explicitly that, thanks to $|\til| = 1$, a seed curve is
always parameterized by arc-length, and so $|\gamma_z'(s)| \equiv
1$. If the base point $z$ is understood or irrelevant, we simply
denote the seed curve by $\gamma(s)$. We will indicate by $\tilL$
the integral curve of $\tilp$ starting at the point $z$, i.e.,
\begin{equation}\label{seedp} \tilLp\ =\ \tilp(\tilL)\
,\quad\quad\quad \tilLzero\ =\ z\ .
\end{equation}
\end{Def}

\medskip

We now use some elementary analysis to investigate the types of
planar vector fields such as $\til$ and $\tilp$.

\medskip

\begin{Lem}\label{L:straight}
If $V(x,y)=(v(x,y),w(x,y))$ is a $C^1$, unit, divergence free vector field on $\R^2$, then the integral curves of the perpendicular vector field $V^\perp = (w,-v)$ are straight lines.
\end{Lem}

\pf An integral curve $\tilde{\mathcal{L}}(r) = (x(r),y(r))$ of
$V^\perp$ satisfies the system
\begin{equation}\label{ic}
\begin{cases}
x'\ =\ w(x,y)\ , \quad\quad\quad y'\ =\ -\ v(x,y)\ ,
\\
x(0)\ =\ x\ ,\quad\quad\quad y(0)\ =\ y\ .
\end{cases}
\end{equation}

Letting $\psi(r) = V^\perp(\tilde{\mathcal{L}}(r)) =
(w(x(r),y(r)), - v(x(r),y(r)))$, we find from \eqref{ic}
\begin{equation}\label{phi'}
\psi'(r)\ =\ \bigg(<\nabla w(\tilde{\mathcal{L}}(r)),
V^\perp(\tilde{\mathcal{L}}(r))>, - <\nabla
v(\tilde{\mathcal{L}}(r)), V^\perp(\tilde{\mathcal{L}}(r))>\bigg)\
.
\end{equation}

Since $V$ is assumed to be unit, by differentiating
the equation $v^2+w^2=1$ with respect to $x$ and $y$ we obtain

\begin{equation}\label{impl}
v\ v_x\ +\ w\ w_x\ =\ 0\ ,\quad\quad\quad v\ v_y\ +\ w\ w_y\ =\ 0\ .
\end{equation}

Since $V$ is also divergence free, we have
\begin{equation}\label{divfree}
v_x\ +\ w_y\ =\ 0\ .
\end{equation}

Combining (\ref{impl}) with (\ref{divfree}) yields
\begin{equation}\label{dot}
<\nabla w, V^\perp>\ \equiv\ 0\ ,\quad\quad\quad\quad <\nabla v, V^\perp>\ \equiv\ 0\ .
\end{equation}

If we use \eqref{dot} in \eqref{phi'} we conclude that $\psi'
\equiv 0$, i.e., the vector field $V^\perp$ is constant along its
own integral curves, hence, in particular, $\tilde{\mathcal L}'(r)
= V^\perp(\tilde{\mathcal L}(r)) \equiv V^\perp(z)$. From the
latter equality we conclude that the integral curves of $V^\perp$
are given by
\begin{equation}\label{icL}
\tilde{\mathcal L}(r)\ =\ (w(z) r + x, - v(z) r + y)\ =\ z\ +\ r\
V^\perp(z)\ , \end{equation}
 hence they are straight lines in
$\R^2$ starting at $z\cong (x,y)$, with $|\tilde{\mathcal L}'(r)|
\equiv 1$.

$\qed$

\medskip

\begin{Cor}\label{C:ltilde}
Let $S$ be an $H$-minimal surface of the type \eqref{S1}. Given a
point $z\in \Om$, the integral curve of $\tilp$ starting at $z$ is
given by a straight line perpendicular to $\gamma_z$
\[
\tilL\ =\ z\ +\ r\ \tilp(z)\ .
\]
\end{Cor}

\pf The hypothesis that $S$ be $H$-minimal guarantees that the
vector field $z \to \til(z)$ is divergence free, see
\eqref{msdiv}, therefore applying Lemma \ref{L:straight} with $V =
\til$ we conclude from its proof that $\tilp$ is constant along
its integral curves, and that the latter are given by \eqref{icL}.

$\qed$

\medskip

According to Definition \ref{D:seedcurve} a seed curve starting at
$z$ is represented by the equation
\[
\gamma_z(s)\ =\ z\ +\ \int_0^s \til(\gamma_z(s))\ ds\ .
\]

From the latter we see that, if the point $g=(z,t)\in S$ is
non-characteristic, then the assumption that $S$ be $C^2$ implies
that the vector field $\zeta\to \til(\zeta)$ be $C^1$ in a
neighborhood of $z$, and thereby from the regularity theory for
solutions of ordinary differential equations we infer that $s\to
\gamma_z(s)$ is of class $C^2$ in a neighborhood of $s=0$. More in
general, is $S$ is of class $C^k$, then $s\to \gamma_z(s)$ if of
class $C^k$ in a neighborhood of $s=0$.

\medskip

\begin{Def}\label{D:seed}
We define the \emph{signed curvature} of $\gamma$ by
\begin{equation}\label{kappa}
\kappa(s)\ \overset{def}{=}\ <\gamma''(s) , \tilp(\gamma(s))>\ =\
<\gamma''(s),\gamma'(s)^\perp> \ =\ \gamma_1''(s) \gamma_2'(s)\ -\
\gamma_2''(s) \gamma_1'(s)\ .
\end{equation}
\end{Def}

\medskip

For later purposes we note that
\begin{equation}\label{kappa2}
\kappa(s)^2\ =\ \gamma_1''(s)^2\ +\ \gamma_2''(s)^2\ .
\end{equation}

One has in fact from \eqref{kappa} \begin{align*} \kappa(s)^2\ &
=\ \gamma_1''(s)^2 \gamma_2'(s)^2\ +\ \gamma_2''(s)^2
\gamma_1'(s)^2\ -\ 2\
\gamma_1'(s)\gamma_2'(s)\gamma_1''(s)\gamma_2''(s) \\
& =\  \gamma_1''(s)^2 (\gamma_1'(s)^2 + \gamma_2'(s)^2)\ +\
\gamma_2''(s)^2 (\gamma_1'(s)^2 + \gamma_2'(s)^2)\ -\ 2\
<\gamma'(s),\gamma''(s)>^2
\\
& =\ \gamma_1''(s)^2\ +\ \gamma_2''(s)^2\ ,
\end{align*}
where we have repeatedly used the identity $|\gamma'(s)|^2 = 1$.

Returning to Definition \ref{D:seedcurve}, if we use
$\{\tilde{\mathcal{L}}_z, \gamma_z\}$ as our coordinate curves, we
have a new local parameterization of the $xy$-plane $F : \R^2 \ra
\R^2$ given by
\begin{equation}\label{newpar}
(s,r)\ \ra\ (x(s,r), y(s,r))\ \overset{def}{=}\ F(s,r)\ =\
\gamma(s)\ +\ r\;\tilp(\gamma(s))\ .
\end{equation}

We note explicitly that $F$ maps the straight line $r=0$ into the
seed curve $\gamma(s)$, i.e., \[ F(s,0)\ =\ \gamma(s)\ .
\]

On the other hand, the straight line $s=0$ is mapped into the
straight line passing through the base point $z$ of the seed curve
and having direction vector $\tilp(z)$, i.e., \[ F(0,r)\ =\ z\ +\
r\ \tilp(z)\ . \]

Thanks to Corollary \ref{C:ltilde} we thus recognize that $F(0,r)
= \tilL$. In particular, when $z = \gamma(s)$ we obtain from this
identity
\begin{equation}\label{L}
\tilLg\ =\ \gamma(s)\ +\ r\ \tilp(\gamma(s))\ .
\end{equation}

Formula \eqref{L} will be useful later in the section. Keeping in
mind \eqref{seeddef} we see that along the seed curve we have
$\tilp(\gamma(s)) = \gamma'(s)^\perp = (\gamma_2'(s) , -
\gamma_1'(s))$. One thus has the following explicit expression for
$F(s,r)$
\begin{equation}\label{explicit}
F(s,r)\ =\ \gamma(s) + r \gamma'(s)^\perp\ =\ \bigg(\gamma_1(s) +
r \gamma_2'(s)\ ,\ \gamma_2(s) - r \gamma_1'(s)\bigg)\ .
\end{equation}

Since $\gamma = \gamma_z$, clearly $F(s,r)$ depends on the choice
of the base point $z\in \Om$ as well. However, since in our
arguments such point will be fixed, we will routinely omit such
dependence. As we observed above if $S$ is $C^2$, then $\gamma \in
C^2$, and therefore $F\in C^1$ in its domain. More in general, if
$S$ is of class $C^k$, then $F\in C^{k-1}$. We next point out that
$F$ fails to be a diffeomorphism everywhere. In fact,
\eqref{newpar} gives
\begin{align}\label{DF}
DF\ =\ \begin{pmatrix} \frac{\partial x}{\partial s} & \frac{\partial x}{\partial r} \\ \frac{\partial y}{\partial s} &
 \frac{\partial y}{\partial r}\end{pmatrix}\ =\ \begin{pmatrix} \gamma_1'(s) + r \gamma_2''(s) & \gamma_2'(s) \\ \gamma_2'(s) - r \gamma_1''(s) & - \gamma_1'(s)\end{pmatrix}\ .
\end{align}

With \eqref{kappa} in mind, we thus obtain from \eqref{DF}
\begin{equation}\label{detDF}
det \;DF\ =\  - \gamma_1'(s)^2 - \gamma_2'(s)^2 +
r\ (\gamma_2'(s)\gamma_1''(s)-\gamma_1'(s)\gamma_2''(s))\ =\ - 1 + r \kappa(s)\ .
\end{equation}

Thus, we have shown the following lemma.

\medskip

\begin{Lem}\label{L:sl}
The mapping $F$ defined in \eqref{newpar} ceases to be a local
$C^1$ parameterization of a surface, in the sense that $DF$ is not
invertible, along the curve in the $(s,r)$-plane given by
\[
\mathcal C_\gamma\ =\ \left\{(s,r)\in \mathbb R^2 \mid r = \frac{1}{\kappa(s) }\right\}\ ,
\]
where $\kappa(s)$ is the signed curvature of the seed curve $\gamma(s)$.
\end{Lem}

\medskip

 Next, we impose the new coordinates \eqref{newpar} on the
function $h$ in \eqref{nu}. Under the new parameterization, the
surface $S$ in \eqref{S1} is given by
\begin{equation}\label{np}
(s,r)\ \longrightarrow\ \bigg(\gamma_1(s)+r\gamma_2'(s),\gamma_2(s)-r\gamma_1'(s),
h\big(\gamma_1(s)+r\gamma_2'(s),\gamma_2(s)-r\gamma_1'(s)\big)\bigg)\ ,
\end{equation}
at least on that portion of the $sr$-plane where the map $F$
defines a diffeomorphism. We stress that, henceforth in the paper,
by abuse of notation we indicate
\[
h(s,r)\ =\ h(F(s,r))\ =\ h(\gamma_1(s) + r \gamma'_2(s) ,
\gamma_2(s) - r \gamma_1'(s))\ ,
\]
with similar understanding
when we express any other function, originally defined on an open set of the $xy$-plane, in terms of the coordinates $(s,r)$. Thus, for instance, $W(s,r)$ has to be interpreted analogously.
We are now ready to establish the main representation theorem for $H$-minimal surfaces of the type \eqref{S1}, see Theorem A in the introduction.

\medskip

\begin{Thm} \label{rep}
Let $k \ge 2$. A patch of a $C^k$ surface $S \subset \Hy^1$ of the
type \eqref{S1}, with empty characteristic locus over $\Om$, is an
$H$-minimal surface if and only if, for every $g=(z,t) \in S$,
there is a neighborhood of $g$ which can be parameterized by
\begin{equation}\label{para}
(\gamma_1(s)+r\gamma_2'(s),\gamma_2(s)-r\gamma_1'(s),
h(s,r))\ ,
\end{equation}
with $h(s,r)$ given in the form
\begin{equation}\label{para2}
h(s,r)\ =\ h_0(s)\ -\ \frac{r}{2} <\gamma(s), \gamma'(s)>\ .
\end{equation}
with
 \[\gamma \in C^{k+1},\quad\quad\quad\text{and}\quad  h_0 \in C^k\ .
\]
We conclude from this, and from \eqref{explicit}, that the map $F$
is not only $C^{k-1}$, but in fact $C^k$. Thus, to specify such a
patch of a smooth $H$-minimal surface of this type, one must
  specify a single curve in $\Hy^1$ determined by a seed curve $\gamma(s)$, parameterized by arc-length, and an
  initial height function $h_0(s)$.
\end{Thm}

\pf First, we prove the necessity of the representation
\eqref{para}, \eqref{para2}. Assuming $S$ is an $H$-minimal graph
over a domain $\Omega$ in the $xy$-plane, with empty
characteristic locus, pick a base point, $g = (z,t) \in S$ with $z
\in \Omega$ being its projection onto the $xy$-plane. Let
$\gamma=\gamma_{z}$ be the integral curve of $\til$ such that
$\gamma(0)= z$. Further, suppose that the domain of $\gamma$ is
the interval $(s_0,s_1)$. Now, for a fixed $s \in (s_0,s_1)$, let
\[
r_0(s)\ =\ \sup\{r>0 \mid F(s,-r) \subset \Omega\}\ ,
\quad\quad\quad r_1(s)\ =\ \sup\{r>0 \mid F(s,r) \subset \Omega\}\
,
\]
and let $O = \{ (s,r)\in \mathbb R^2 \mid s \in (s_0,s_1), r \in
(-r_0(s),r_1(s))\}$. By construction we have that, if
$F(O)=\Omega_0$, then $(s,r) \to (F(s,r),h(s,r))$ maps $O$ onto
$S_0$, where $S_0$ is the portion of $S$ which is a graph over
$\Omega_0$.  Now, for any fixed $s\in (s_0,s_1)$, consider
$h_{s}(r) \overset{def}{=} h(s,r)$, i.e., the function $h$
restricted to a curve parameterized by $r$. As the trace of such
curve lies on the surface, its tangent vector $\boldsymbol
\tau(r)$ must be perpendicular to the Riemannian normal to $S$
along the curve. For the subsequent calculations it might be
helpful for the reader to keep in mind the following formula which
gives the relation between the expression of a vector field
$(x,y,t) \mapsto (a(x,y,t), b(x,y,t), c(x,y,t))$ in the
rectangular coordinates of $\mathbb R^3$, and that with respect to
the orthonormal basis $\{X_1,X_2,T\}$,
\begin{equation}\label{coord}
(a,b,c)\ =\ a\ X_1\ +\ b\ X_2\ +\ \left (c + \frac{1}{2} a y -
  \frac{1}{2} b x \right)\ T\ .
\end{equation}

We first compute $\boldsymbol \tau(r)$ using the representation \eqref{np} of $S$
\begin{equation}\label{tau1}
\begin{split}
\boldsymbol \tau(r) &= \left(\gamma_2'(s), - \gamma_1'(s), h_x
\frac{\p x}{\p r}+h_y\frac{\p y}{\p r}\right)\; \; \text{(by the
definition \eqref{explicit} of
  $F$)}\\
&= (\gamma_2'(s), - \gamma_1'(s), h_x
\gamma_2'(s)-h_y\gamma_1'(s)) \\
&= \gamma_2'(s) X_1 - \gamma_1'(s) X_2 + \left(<\nabla h , \tilp>
+\
\frac{1}{2} <\gamma(s),\gamma'(s)> \right) T\ ,\\
\end{split}
\end{equation}
where we have used \eqref{coord} along with the fact that in the coordinates $(s,r)$ one has
\begin{equation}\label{xy}
\begin{split}
x\ &=\  x(s,r)\ =\ \gamma_1(s)\ +\ r \gamma_2'(s)
\\
y\ & =\ y(s,r)\ =\ \gamma_2(s)\ -\ r \gamma_1'(s)\ .
\end{split}
\end{equation}

The Riemannian unit normal on the surface is given by \eqref{nu}, which, according to \eqref{type1} and \eqref{isobar},  we can write as
\begin{align}\label{run}
\boldsymbol \nu\ & =\ \frac{W(s,r)}{\sqrt{1 + W(s,r)^2}}\ \left\{\overline p\ X_1\ +\ \overline q\ X_2\right\}\ +\ \frac{1}{\sqrt{1 + W(s,r)^2}}\ T
\\
& =\
\frac{W(s,r)}{\sqrt{1 + W(s,r)^2}}\ \nuX\ +\ \frac{1}{\sqrt{1 + W(s,r)^2}}\ T\ .
\notag
\end{align}

Since $\boldsymbol \tau$ and $\boldsymbol \nu$ must be orthogonal, we see that
\begin{align}
& 0\ =\ <\boldsymbol \nu, \boldsymbol \tau>
\\
&=\ \frac{1}{\sqrt{1 +
W(s,r)^2}}\bigg\{(\gamma_2'(s)\gamma_1'(s)-\gamma_1'(s)\gamma_2'(s))\
W(s,r) + \bigg(<\nabla h,\tilp> + \frac{1}{2} <\gamma(s),
\gamma'(s)>\bigg)\bigg\}
\notag\\
&=\ \frac{1}{\sqrt{1 + W(s,r)^2}} \bigg(<\nabla h,\tilp>\ +\
\frac{1}{2} <\gamma(s), \gamma'(s)>\bigg) \ . \notag
\end{align}

Thus, we have the equation
\begin{equation}\label{vanish}
\frac{\partial h}{\partial r}(s,r)\ =\ <\nabla h(s,r) ,
\tilp(\gamma(s))>\ =\ -\ \frac{1}{2} <\gamma(s), \gamma'(s)>\ ,
\end{equation}
which implies that
\[
h(s,r)\ =\ -\ \frac{r}{2} <\gamma(s), \gamma'(s)>\ +\ h(s,0)\ .
\]

Letting
\begin{equation}\label{hf}
h_0(s)\ \overset{def}{=}\ h(s,0)\ =\ h(\gamma(s))\ ,
\end{equation}
we have the desired result that $S$ is described as in \eqref{para}, with $h(s,r)$ given by
\[
h(s,r)\ =\ h_0(s)\ -\ \frac{r}{2} <\gamma(s),\gamma'(s)>\ ,
\]
which is \eqref{para2}.

We next want to show that $\gamma \in C^{k+1}$, and that $h_0 \in
C^k$. Before doing this, however, we prove that $F:O \ra \Omega_0$
is one-to-one, and hence is a diffeomorphism.  To show this, we
argue by contradiction, assuming that there exist two points
$(s_2,r_2) \not= (s_3,r_3)$ in $O$ such that
\begin{equation}\label{F=}
F(s_2,r_2)\ =\ F(s_3,r_3)\ . \end{equation}

 Since $S$ is a graph, \eqref{F=} implies
  $(F(s_2,r_2),h(s_2,r_2))=(F(s_3,r_3),h(s_3,r_3))$. Also, according to
\eqref{newpar}, \eqref{L}, \eqref{F=} is equivalent to saying that
$\mathcal L_{\gamma(s_2)}(r_2) =  \mathcal L_{\gamma(s_3)}(r_3)$.
This implies that $\mathcal{L}_{\gamma(s_2)} \cap
\mathcal{L}_{\gamma(s_3)} \cap \Omega_0 \neq \varnothing$. We
notice that the possibility $s_2 = s_3$ and $r_2 \not= r_3$ is
excluded, since in view of \eqref{L}, \eqref{F=} can be re-written
\begin{equation}\label{intersection}
\gamma(s_2)\ +\ r_2\ \tilp(\gamma(s_2))\ =\ \gamma(s_3)\ +\ r_3\
\tilp(\gamma(s_3))\ , \end{equation}
 which if $s_2 = s_3$ would reduce to $(r_2 - r_3)
\tilp(\gamma(s_2)) = 0$. Since $|\tilp(\gamma(s_2))| = 1$, the
latter equation is impossible, unless $r_2 = r_3$. We must thus
have $s_2 \not= s_3$. Assuming this, if $\gamma(s_2) =
\gamma(s_3)$, then we have also $\tilp(\gamma(s_2)) =
\tilp(\gamma(s_3))$. But then the two straight lines represented
by $\mathcal L_{\gamma(s_2)}(r)$ and $\mathcal L_{\gamma(s_3)}(r)$
are parallel and therefore, unless they coincide, they cannot
intersect in $\Om_0$ thus violating the assumption $\mathcal
L_{\gamma(s_2)}(r_2) = \mathcal L_{\gamma(s_3)}(r_3)$. However,
they cannot coincide unless $r_2 = r_3$. If this happens (i.e., if
$r_2 = r_3$ and the two straight lines coincide) then the
parametric curve $\gamma(s)$ would self-intersect ($\gamma(s_2) =
\gamma(s_3)$) with equal tangents at the intersection point. Since
this would contradict the fact that $\gamma(s)$ is a smooth
embedded curve, we conclude that we must have $\gamma(s_2) \not=
\gamma(s_3)$. Again, this forces
\begin{equation}\label{intersect}
\tilp(\gamma(s_2))\ \not=\ \tilp(\gamma(s_3))\ , \end{equation}
since otherwise the two lines $\tilde{\mathcal L}_{\gamma(s_2)}$,
$\tilde{\mathcal L}_{\gamma(s_3)}$ would be parallel and therefore
they could not possibly intersect. Finally, we observe that
\eqref{intersect} is equivalent to saying $\gamma'(s_2) \not=
\gamma'(s_3)$.

Keeping in mind that $\nuXp = \bar q X_1 - \bar p X_2$, we
conclude from \eqref{intersect} that we must have
\[
\lim_{r \ra r_2} \nuXp((F(s_2,r),h(s_2,r)))\ \neq \ \lim_{r \ra
r_3} \nuXp((F(s_3,r),h(s_3,r)))\ .
\]

 Thus, $\nuXp$ (and hence
$\nuX$) cannot be continuous at the point
$(F(s_2,r_2),h(s_2,r_2))=(F(s_3,r_3),h(s_3,r_3))\in S$ . However,
the non-characteristic assumption implies that $\nuX$ is
continuous on $S$, and we have thus reached a contradiction. This
proves that $F:O \ra \Omega_0$ is one-to-one and onto.  Since it
is a local diffeomorphism as well, we conclude $F|_O$ is a
diffeomorphism of $O$ onto $\Omega_0$.

We now turn to proving that $h_0 \in C^k$ and $\gamma \in
C^{k+1}$. To this effect, we analyze the curve
$(\gamma(s),h_0(s))$. Looking at the parameterization we have just
shown, we have $S$ given as
\[
\left (\gamma_1(s)+r
\gamma_2'(s),\gamma_2(s)-r\gamma_1'(s),h_0(s)-\frac{r}{2}<\gamma(s),\gamma'(s)> \right)\ .
\]

As observed following Definition \ref{D:seedcurve}, if $S$ is
$C^k$ then, away from characteristic points, $\gamma$ is at least
$C^k$. Therefore, since from \eqref{hf} we have $h_0(s) =
h(F(s,0)) = h(\gamma(s))$, from the chain rule and the fact that
$h\in C^k(\Om)$, we conclude that $h_0 \in C^k$ as well. We want
to show next that $\gamma \in C^{k+1}$.

\begin{Rmk}\label{R:intcurve}
Suppose $S$ is parameterized as in Theorem \ref{rep} by a single
seed curve and height function. For a fixed $r$, consider
\[\gamma_{r}(s)\ =\ (\gamma_1(s)+r\gamma_2'(s),\gamma_2(s)-r\gamma_1'(s))\ =\ \gamma(s)\ +\ r\ \tilp(\gamma(s))\ .\]
We note that this curve is a re-parametrization of the integral
curve of $\til$ passing through $\gamma_r(0)$. To see this, we
observe that, proceeding as in the proof of Lemma
\ref{L:straight}, one can show that $\til$ is constant along the
straight lines given by $r \to
(\gamma_1(s)+r\gamma_2'(s),\gamma_2(s)-r\gamma_1'(s))$, with $s$
fixed. Keeping in mind \eqref{explicit} we thus conclude that for
any fixed $s$, and for every $r$, one has
\[ \til(\gamma(s) + r \gamma'(s)^\perp)\ =\ \til(F(s,r))\ =\
\til(F(s,0))\
 =\ \til(\gamma(s))\ =\
\gamma'(s)\ . \]

 Since
$\gamma_r'(s)=(\gamma_1'(s)+r\gamma_2''(s),\gamma_2'(s)-r\gamma_1''(s))$,
we see that
\begin{equation}
\begin{split}
& <\gamma_r'(s), \tilp(\gamma(s))>\ =\ <\gamma_r'(s) ,
\gamma'(s)^\perp>
\\
& =\ <\gamma'(s),\gamma'(s)^\perp> + r <\gamma''(s)^\perp,
\gamma'(s)^\perp>\ =\ 0 \ .
\end{split}
\end{equation}
the last equality being true since $|\gamma'(s)| \equiv 1$. This
implies that $\gamma_r'(s) = \alpha(s) \gamma'(s)$, for some
function $\alpha(s)$. We also note that \eqref{kappa} and
\eqref{kappa2} show
\[
|\gamma_r'(s)|^2\ =\ 1\ +\ r^2 (\gamma_1''(s)^2 +
\gamma_2''(s)^2)\ -\ 2 r \kappa(s)\ =\ (1-r\kappa(s))^2\ .
\]
We thus conclude that $|\alpha(s)| = |1 - r \kappa(s)|$.
\end{Rmk}

For a fixed $r=r_0$, we see that the image of
\[
\gamma_{r_0}(s)\ =\ \left (\gamma_1(s)+r_0
\gamma_2'(s),\gamma_2(s)-r_0\gamma_1'(s) \right)
\]
coincides with the image integral curve of $\til$ through
$(\gamma_1(s)+r_0 \gamma_2'(s),\gamma_2(s)-r_0\gamma_1'(s))$ (see
Remark \ref{R:intcurve}).  Let $g_{r_0}(\sigma)$ be this integral
curve. Since
\[ \gamma_{r_0}'(s) =   \left (\gamma_1'(s)+r_0
\gamma_2''(s),\gamma_2'(s)-r_0\gamma_1''(s) \right) \] and \[
|\gamma_{r_0}'(s)| = |1 - r_0 \kappa(s)|\ .
\]

Let
\[
\kappa_0\ \overset{def}{=}\ \underset{(s_0,s_1)}{\sup}\
|\kappa(s)|\ ,
\]
then if $\kappa_0 >0$ we take $r_0 \leq \kappa_0$, whereas if
$\kappa_0 = 0$, then we choose $r_0$ to be any positive number.
With such choice of $r_0$ it is clear that
\[
|1\ -\ r_0\ \kappa(s)|\ =\ 1\ -\ r_0\ \kappa(s) \] for every $s\in
(s_0,s_1)$, and therefore $\gamma_{r_0}$ is simply a
re-parameterization of the integral curve, $g_{r_0}$. Again, as
the surface is $C^k$, we have that $g_{r_0}$ is at least $C^k$ as
well. To show the desired regularity of $\gamma$, we will consider
the defining function of the surface, $h$, restricted to
$g_{r_0}$. Again, since the surface is $C^k$, $h$ is $k$ times
differentiable along $g_{r_0}$. First, we differentiate $h$ with
respect to $\sigma$.
\[ \frac{\partial}{\partial \sigma} h(s,r) =
\frac{\partial h}{\partial
  s}\  \frac{d s}{d \sigma} \]
Fixing $r_0$ sufficiently close to zero (so that $|r_0 \kappa(s)| < 1$), we know that $|1-r_0
\kappa(s)|=1-r_0\kappa(s)$ and so $\frac{\partial s}{\partial
  \sigma} = \frac{1}{1-r_0 \kappa(s)}= 1+ r_0 \kappa(s) + r_0^2
\kappa(s)^2 + \dots$.  Thus,
\begin{equation*}
\begin{split}
\frac{\partial}{\partial \sigma} h(s,r) &=
\left(h_0'(s)-\frac{r_0}{2} - \frac{r_0}{2} <\gamma(s),
\gamma''(s)>\right )\left(  1+ r_0 \kappa(s) + r_0^2
\kappa(s)^2 + \dots \right) \\
&= h_0'(s) + r_0 \alpha_1+r_0^2 \alpha_2 + \dots
\end{split}
\end{equation*}
where, for example, \[\alpha_1=\left(-\frac{1}{2} - <\gamma(s),
\gamma''(s)> +
  \kappa(s)h_0'(s)\right )\]

We recall that if $f+g$ is differentiable and $f$ is differentiable
that $g$ must be differentiable as well.  Thus, as $h_0'$ is differentiable
we have that $r_0 \alpha_1+r_0^2 \alpha_2 + \dots$ is as well.
Moreover, as $r_0$ is independent of $\sigma$, we can divide by $r_0$
yielding that $\alpha_1 + r_0 \alpha_2 \dots$ is differentiable.  As
the function $k(s,r_0)=\sum_{i=1}^\infty \alpha_i r_0^i$ is infinitely
differentiable in $r_0$, for it to be differentiable in $s$ as well,
each $\alpha_i$ must be differentiable in $s$.  For example,
$\alpha_1$ must be differentiable in $s$ since
$k(s,r_0)-r_0k_{r_0}(s,r_0)=\alpha_1(s)$.  Considering $\alpha_1$, we
note that
\begin{equation*}
\begin{split}
\alpha_1(s)&= -\frac{1}{2} - <\gamma(s),\gamma''(s)> +
  \kappa(s)h_0'(s)\\
 &= -\frac{1}{2} -
  \gamma_1(s)\gamma_1''(s)-\gamma_2(s)\gamma_2''(s)+h_0'(s)(\gamma_1''(s)\gamma_2(s)-\gamma_1''(s)\gamma_1(s))\\
&=  -\frac{1}{2} +\gamma_1''(s)\beta_1(s)+\gamma_2''(s)\beta_2(s) \\
\end{split}
\end{equation*}
where $\beta_1(s)= h_0'(s)\gamma_2'(s)-\gamma_1(s)$ and
$\beta_2(s)=-h_0'(s)\gamma_1'(s)-\gamma_2(s)$.  We note that both
$\beta_1$ and $\beta_2$ are differentiable.  Moreover, since
$\gamma_1'(s)\gamma_1''(s) + \gamma_2'(s)\gamma_2''(s)=0$ we have
\begin{equation*}
\begin{split}
\gamma_1'(s)\gamma_1''(s)\beta_1(s)+\gamma_1'(s)\gamma_2''(s)\beta_2(s)
&= -\gamma_2'(s)\gamma_2''(s)
\beta_1(s)+\gamma_1'(s)\gamma_2''(s)\beta_2(s)\\
&= \gamma_2''(s) (-\gamma_2(s)\beta_1(s)+\gamma_1(s)\beta_2(s))
\end{split}
\end{equation*}

Now, putting these together we have that
\[ \frac{ \gamma_1'(s)\left ( \alpha_1(s) - \frac{1}{2}
  \right)}{-\gamma_2(s)\beta_1(s)+\gamma_1(s)\beta_2(s)}=
  \gamma_2''(s)\]
Since all of the functions on the left hand side of the equation are
  differentiable along $g_{r_0}$, we have that $\gamma_2''(s)$ is
  differentiable in $\sigma$ and hence in
  $s$.  As, $<\gamma', \gamma''> =0$, this shows that
  $\gamma_1''(s)$ is differentiable as well.  In other words, if
  the $H$-minimal surface is at least $C^2$, then $\gamma$ is at least
  $C^3$ and $h_0$ is at least $C^2$.

If, in addition, we assume that the surface is $C^k$, we have that
\[ \frac{ \gamma_1'(s)\left ( \alpha_1(s) - \frac{1}{2}
  \right)}{-\gamma_2(s)\beta_1(s)+\gamma_1(s)\beta_2(s)}=
  \gamma_2''(s)\]
is differentiable $k-1$ times in $\sigma$.  A straightforward (but
  somewhat tedious) calculation shows that this implies that $\gamma
  \in C^{k+1}$.


We now prove the reverse implication. We assume that the surface $S$ is parameterized by \eqref{para}, \eqref{para2}, i.e.,
\[
\bigg(F(s,r)\ ,\  h_0(s)\ -\ \frac{r}{2} <\gamma(s), \gamma'(s)>\bigg)\ ,
\]
where $F(s,r)$ is as in \eqref{newpar}, \eqref{explicit}, $\gamma(s)$ is a $C^{k+1}$ curve parameterized by arc-length and
$h_0(s) \in C^k(\R)$.  We then show that $S$ is a $C^k$ $H$-minimal
surface.  To do this we will compute the horizontal Gauss map $\nuX$ for
$S$ and show that it is divergence free. We will calculate the $s$ and
$r$ derivatives of the parameterization yielding two tangent vectors
to the surface from which to compute the Riemannian normal.

As above, we obtain from \eqref{newpar} for the $r$ derivative
\begin{align}\label{tau}
\boldsymbol \tau\ & =\ \bigg(F_r(s,r)\ ,\ -\ \frac{1}{2} <\gamma(s),\gamma'(s)>\bigg)\ =\ \bigg(\boldsymbol \nu_H^\perp(\gamma(s))\ ,\ -\ \frac{1}{2} <\gamma(s),\gamma'(s)>\bigg)
\\
& =\  \gamma_2'(s)\ X_1\ -\ \gamma_1'(s)\ X_2\ ,
\notag
\end{align}
where in the last equality we have used \eqref{coord}, and the equation $\boldsymbol \nu_H^\perp(\gamma(s)) = (\gamma_2'(s),- \gamma_1'(s))$.
We note that \eqref{tau} shows, in particular, that the vector $\boldsymbol \tau$ is horizontal.
We next compute the $s$ derivative of the new parameterization of $S$. Keeping in mind the convention $z^\perp = y - i x$, using \eqref{coord} we find
\begin{align}\label{sigma}
\boldsymbol \sigma\ & =\ \bigg(F_s(s,r)\ ,\ h_0'(s) - \frac{r}{2} |\gamma'(s)|^2 - \frac{r}{2} <\gamma(s),\gamma''(s)>\bigg)
\\
& =\ \bigg(\gamma'(s) + r\  (\gamma''(s))^\perp\ ,\ h_0'(s) - \frac{r}{2} - \frac{r}{2} <\gamma(s),\gamma''(s)>\bigg)
\notag\\
& =\  (\gamma_1'(s)+r\gamma_2''(s))\ X_1\ +\ (\gamma_2'(s)-r\gamma_1''(s))\ X_2
\notag\\
& +\ \left(h_0'(s) - r - \frac{1}{2} <\gamma'(s),\gamma(s)^\perp> + \frac{r^2}{2} \kappa(s) \right)\ T\ .
\notag
\end{align}

To complete the proof, we need to calculate the divergence of $\nuX$, and prove that the latter vanishes identically. In order to obtain a convenient expression of $\nuX$, we note that taking the cross-product, with respect to the orthonormal frame $\{X_1,X_2,T\}$, of the vectors $\sigma$ and $\tau$ tangent to $S$, we easily obtain
\begin{align}\label{wedge}
\boldsymbol \sigma \times \boldsymbol \tau\ & =\  \gamma'_1(s)\ \left(h_0'(s) - r - \frac{1}{2} <\gamma'(s),\gamma(s)^\perp> + \frac{r^2}{2} \kappa(s)\right)\ X_1
\\
& +\ \gamma'_2(s)\ \left(h_0'(s) - r - \frac{1}{2} <\gamma'(s),\gamma(s)^\perp> + \frac{r^2}{2} \kappa(s)\right)\ X_2\ +\ (- 1 + r \kappa(s))\ T\ .
\notag
\end{align}

The equation \eqref{wedge} yields the following claim.

\medskip

\begin{Claim} For $(s,r)$ such that
\[
h_0'(s) - r - \frac{1}{2} <\gamma'(s),\gamma(s)^\perp> + \frac{r^2}{2} \kappa(s)\ \neq\ 0\ ,
\]
we let
\[
\beta(s,r)\ =\
\frac{-1+r\kappa(s)}{h_0'(s) - r - \frac{1}{2} <\gamma'(s),\gamma(s)^\perp> + \frac{r^2}{2} \kappa(s)}\ .
\]
We thus see that
\[
\boldsymbol \eta\ =\ \gamma_1'(s)\ X_1\ +\ \gamma_2'(s)\ X_2\ +\ \beta(s,r)\ T
\]
points
in the same direction as the
Riemannian normal to $S$ and hence the horizontal Gauss map for $S$ is given by
\begin{equation}\label{goodnuX}
\nuX\ =\ \gamma_1'(s) X_1\ +\ \gamma_2'(s) X_2\ .
\end{equation}
\end{Claim}

Having established this, and recalling \eqref{xy}, a routine computation using the inverse function theorem and \eqref{DF}, \eqref{detDF}, yields that
\begin{equation}\label{ps_xy}
\begin{split}
\frac{\partial s}{\partial x}\ & =\ \frac{ \gamma_1'(s)}{1-r\kappa(s)}\ ,
\\
\frac{\partial s}{\partial y}\ & =\ \frac{\gamma_2'(s)}{1-r\kappa(s)}\ .
\end{split}
\end{equation}

We thus obtain from \eqref{goodnuX}, \eqref{ps_xy}
\begin{equation*}
\begin{split}
div\ \nuX\ =\ (\gamma_1'(s))_x+ (\gamma_2'(s))_y\ &=\ \gamma_1''(s) \frac{\partial s}{\partial x}\ +\
\gamma_2''(s) \frac{\partial s}{\partial y}\\
&=\  \frac{\gamma_1''(s)
  \gamma_1'(s)+\gamma_2''(s)\gamma_2'(s)}{1-r\kappa(s)}\\
&=\ 0\ .
\end{split}
\end{equation*}

The last equality stems from the assumption that $\gamma(s)$ is
parameterized by arc-length and hence $\gamma'(s)$ is orthogonal to
$\gamma''(s)$. In view of \eqref{msdiv} we conclude that $S$ is $H$-minimal, and this completes the proof.

 $\qed$

\medskip

\begin{Rmk}\label{R:char}
We note that, as we see from \eqref{sigma}, if
\begin{equation}\label{remchar}
h_0'(s) - r - \frac{1}{2} <\gamma'(s),\gamma(s)^\perp> + \frac{r^2}{2} \kappa(s)\ =\ 0\ ,
\end{equation}
for some $(s,r)$, then the vector $\boldsymbol \sigma$ is horizontal. Since the
other tangent vector to $S$, $\boldsymbol \tau$, is always horizontal, see
\eqref{tau}, we conclude that if the parameterization were extended to
include such a point, the point is a
characteristic point.  When examining individual H-minimal surfaces in
terms of seed curves and height functions,
as we will in section \ref{S:return}, this equation is useful for
determining the characteristic locus.
\end{Rmk}

\medskip

We end this section by noting that in a different case, opposite
to the case of graphs over the $xy$-plane, the $H$-minimal
surfaces are still ruled surfaces.

\medskip

\begin{Lem}\label{trivial}  If $S$ is a $C^2$ $H$-minimal surface so
  that no portion of $S$ can be written as a graph over the $xy$-plane, then $S$ is a subset of a single
  vertical plane.
\end{Lem}
\pf  Suppose $S$ is given locally by $\phi(x,y,t)=0$.  Then, since no portion
of $S$ can be written as a graph over the $xy$-plane, we must have
$\phi_t=0$.  Then, $\phi(x,y,t)=\phi_o(x,y)$
and $S$ is a ruled surface whose generators are straight lines
perpendicular (in $\R^3$) to a curve $\gamma \subset \R^2 \times \{0\}$.  Using
the implicit function theorem in $\R^2$, we can locally write $\gamma$ as either $y=f(x)$ or $x=g(y)$.  For $y=f(x)$, we have the
patch of surface defined by \[f(x)-y=0\]  In this case, $p=f'(x),
q=-1$ and
\[
\nuX\ =\ \frac{f'(x)}{\sqrt{1+f'(x)^2}}\ X_1\ -\ \frac{1}{\sqrt{1+f'(x)^2}}\ X_2\ .
\]

Using the assumption that $S$ is $H$-minimal, we take the
divergence of $\nuX$ to find that
\[
0\ =\ div\ \nuX\ =\ \frac{d}{d x} \left ( \frac{f'}{\sqrt{1+(f')^2}} \right) \ =\
\frac{f''}{(1+(f')^2)^{\frac{3}{2}}}\ ,
\]
which gives $f'' = 0$. This shows that the planar curve $\gamma$ must be a straight line and $\nuX$ is
constant. A similar computation yields the same conclusion in the case $x=g(y)$.
Thus, surfaces of this type must be composed of pieces of vertical
planes.  As $S$ is $C^2$, it must be a portion of a single vertical plane.

$\qed$

\vskip 0.6in


\section{\textbf{$H$-minimal surfaces as ruled surfaces}} \label{S:ruled}

\vskip 0.2in

In this section, based on the work in Section
\ref{S:representation}, we put forth a geometric interpretation of
$H$-minimal surfaces which can be written as a graph over a
portion of the $xy$-plane.  In what follows, we will indicate with
$d_{cc}$ the Carnot-Carath\'eodory distance in $\mathbb H^1$
generated by the vector fields $\{X_1,X_2\}$. In view of Theorem
\ref{rep} we see that if $S \subset \Hy^1$ is a $C^2$ $H$-minimal
surface of the type \eqref{S1}, with empty characteristic locus,
then there exist $\mathcal Q = \{(s,r) \in \mathbb R^2 \mid s_0 <
s <s_1\ ,\ r_0 < r < r_1\}$, such that letting
\[
d(s)\ =\ (\gamma_1(s),\gamma_2(s), h_0(s))\ ,
\]
and
\[
v(s)\ =\ (\gamma_2'(s),-\gamma_1'(s),0)\ ,
\]
we have
\begin{align}\label{repre}
 S\ & =\ \big\{d(s) \circ \delta_r v(s) \mid (s,r)\in \mathcal
Q\big\} \\
&  =\ \bigg\{\big(\gamma_1(s) + r \gamma_2'(s), \gamma_2(s) - r
\gamma_1'(s), h_0(s) - \frac{r}{2} <\gamma(s),\gamma'(s)>\big)
\mid  (s,r)\in \mathcal Q\bigg\}\ , \notag
\end{align}
 where $\delta_r(x,y,t) = (r x, r y, r^2 t)$ is the anisotropic
dilation of $\mathbb H^1$, and we have indicated with $\circ$ the
group law \eqref{gl}. We note explicitly that, since $v(s)$
belongs to the horizontal plane passing through the origin (the
$xy$-plane), we have
\begin{equation}\label{one}
\cc(v(s),0)\ =\ \sqrt{\gamma_1'(s)^2\ +\ \gamma_2'(s)^2}\ =\ 1\ .
\end{equation}

\medskip

\begin{Cor}\label{C:geo} Given a $C^2$ non-characteristic $H$-minimal
  graph $S$ over the $xy$-plane parameterized as in \eqref{repre} via Theorem \ref{rep}, for every fixed $s_0$ the curve
  \begin{equation}\label{rule}
  \mathcal{L}(r)\ =\ d(s_0)\ \circ\ \delta_r v(s_0)
  \end{equation}
   is a portion of a geodesic in $\Hy^1$.
\end{Cor}

\pf It suffices to observe the following elementary properties of the distance $\cc$
\begin{align*}
\cc(\mathcal L(r_1), \mathcal L(r_2))\ & =\ \cc(d(s)\circ \delta_{r_1}v(s),d(s)\circ \delta_{r_2}v(s))\ =\ \cc(\delta_{r_1}v(s),\delta_{r_2}v(s))
\\
& =\ \cc((\delta_{r_2}v)^{-1} \circ \delta_{r_1}v(s), 0)\ =\ \cc(\delta_{r_1-r_2}v(s),0)\ =\ (r_1-r_2)\ \cc(v(s),0)
\\
& =\ r_1 - r_2\ ,
\end{align*}
thanks to \eqref{one}.

 $\qed$

\medskip

\begin{Def} A \emph{ruled surface} in $(\mathbb{H}^1,\cc)$ is a surface
  $S$ which is foliated by geodesics of $(\mathbb{H}^1,\cc)$.  Such a
  geodesic is called a \emph{rule}. A curve
  $d(s)$ transverse to the foliation is called a \emph{directrix} of
  the ruled surface.
\end{Def}

\medskip

Corollary \ref{C:geo} and the representation \eqref{repre} yield a nice geometric characterization of
$H$-minimal surfaces.

\medskip

\begin{Cor} \label{C:geo2}
 $S$ is a portion of a $C^2$ non-characteristic $H$-minimal graph over the $xy$-plane in
  $(\mathbb{H}^1,\cc)$, parameterized as in Theorem \ref{rep}, if and only if it is a piece of a ruled surface where the rulings are all straight
  lines
  in $(\mathbb{H}^1,\cc)$.
\end{Cor}

\pf  The geodesics in the ruling are given by the curves
$\mathcal{L}(r)$ defined by \eqref{rule} in Corollary \ref{C:geo}.
Thanks to \eqref{repre} such geodesic is a straight line in
$\mathbb H^1$ (identified with $\mathbb R^3$). A directrix of the
ruling is parameterized by
\[
d(s)\ =\ (\gamma_1(s),\gamma_2(s), h_0(s))\ .
\] $\qed$

\medskip

\begin{Rmk}
Note that the directrix contains two pieces of information, the first
two coordinates compose a seed curve for $S$ and the
last coordinate yields the initial height function $h_0(s)$.  If we
were using completely classical nomenclature, we would call the curve $(\gamma(s),h_0(s))$ a line of striction since the curve is orthogonally transverse to the rules.
\end{Rmk}

\medskip

\begin{Rmk}  The results of this section can be extended in several directions.  First, similar results hold for the analogue of
non-parametric constant
mean curvature surfaces in the Heisenberg group - they are ruled
surfaces as well but are foliated by geodesics which are lifts of
circles in the plane.  Second, both of these results can be
generalized to some more general Carnot groups.  Third, some of
the regularity assumptions can be weakened yielding similar
representation results.  These results
and some of their applications will be addressed in a future paper.
\end{Rmk}

\vskip 0.3in


\section{\textbf{Further analysis of the examples}}\label{S:return}

\vskip 0.2in

At this point, it might be helpful for the reader to return to the examples of Section \ref{S:examples} with the purpose of explicitly computing their seed curves $\gamma(s)$, the diffeomorphism $F(s,r)$, the signed curvature $\kappa(s)$ and their singular locus $\mathcal C_\gamma$.

\medskip

\begin{Example}\label{E:comp}
 \emph{Consider the characteristic plane $\Pi$ in Example \ref{E:charplane}, whose horizontal Gauss map is given by
\[
\nuX\ =\ \left(- \frac{y}{|z|}, \frac{x}{|z|}\right)\ ,\quad\quad\quad\quad z = (x,y) \in \Om \overset{def}= \mathbb R^2 \setminus \{(0,0)\}\ .
\]
For every $z\in \Om$ the seed curve $\gamma(s) = \gamma_z(s)$ solves the system}
\[
\begin{cases}
\gamma_1'(s)\ =\ -\ \frac{\gamma_2(s)}{|\gamma(s)|}\ ,\quad\quad\quad \gamma_1(0)\ =\ x\ ,
\\
\gamma_2'(s)\ =\  \quad \frac{\gamma_1(s)}{|\gamma(s)|}\ ,\quad\quad\quad \gamma_2(0)\ =\ y\ .
\end{cases}
\]
\emph{It is easy to recognize that if $\gamma(s)$ solves this system, then it must be $|\gamma(s)| \equiv const. \equiv |z|$. Using this latter information we find that the seed curve is given by the circle
\begin{equation}\label{seedpi}
\gamma(s)\ =\ \left(x\ \cos\left(\frac{s}{|z|}\right) - y\ \sin\left(\frac{s}{|z|}\right)\ ,\ y\ \cos\left(\frac{s}{|z|}\right) + x\ \sin\left(\frac{s}{|z|}\right)\right)\ .
\end{equation}
Keeping in mind that $\nuX^\perp(\gamma(s)) = (\gamma_2'(s), - \gamma_1'(s))$, from \eqref{kappa} we obtain
\[
\kappa(s)\ \equiv\ -\ \frac{1}{|z|}\ .
\]
The singular locus in Lemma \ref{L:sl} is thus presently given by the straight line
\[
\mathcal C_\gamma\ =\ \left\{(s,r)\in \mathbb R^2 \mid r = - |z|\right\}\ ,
\]
and the map
\begin{align*}
F(s,r)\ & =\ \left(1 + \frac{r}{|z|}\right)\ \gamma(s)
\\
& =\ \left(1 + \frac{r}{|z|}\right)\ \left(x\  \cos\left(\frac{s}{|z|}\right) - y\sin\left(\frac{s}{|z|}\right)\ ,\ y\ \cos\left(\frac{s}{|z|}\right) + x\sin\left(\frac{s}{|z|}\right)\right)\ ,
\end{align*}
defines a local diffeomorphism on either side of such line. Of course, $F(s,r)$ fails to be globally one-to-one, however by suitably restricting the range of $(s,r)$ we obtain a global diffeomorphism.
For instance, if we let
\[
U\ =\ \{(s,r)\in \mathbb R^2 \mid 0 < s < 2 \pi |z|\ ,\ r > - |z|\}\ ,
\]
then $F : U \to \Om $ defines a global diffeomorphism onto. Notice that in the present example we have $h(x,y) \equiv 0$, and therefore the height function defined by \eqref{hf} is given by $h_0(s) = 0$. We next use the equation \eqref{remchar} in Remark \ref{R:char} to determine the characteristic locus of $\Pi$. Such equation presently reduces to
\[
r^2\ + 2 |z| r\ +\ |z|^2\ =\ 0\ ,
\]
which admits the double root $r = - |z|$. Since the parameterization of $\Pi$ with respect to the coordinates $(s,r)$ is given by
\[
(s,r)\ \to\ (F(s,r), 0)\ ,
\]
 we conclude that the characteristic locus $\Sigma$ is given by the
 image of the line $r = - |z|$ through such map. We readily see that
 $\Sigma = \{(0,0,0)\}$, in accordance with what we found in the
 discussion of Example \ref{E:charplane}.  Figure \ref{eg5.1} shows
 the vector field $\nuX$, the seed curve $\gamma(s)$ (in black) and
 some representative rules (in grey) to help illustrate this example.
 Note that the origin corresponds, as shown above, to both the
 characteristic locus and the image of the singular locus.
\begin{figure}
\epsfig{figure=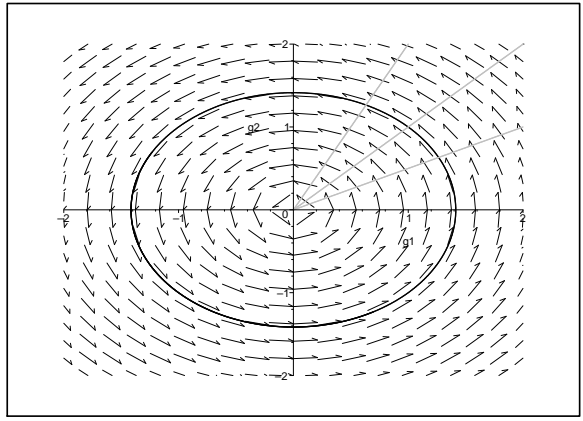,angle=-90,width=3in}
\caption{A seed curve and the projections of some rules to the
$xy$-plane for the surface $\{t=0\}$}\label{eg5.1}
\end{figure}
From the same discussion, it is easy to construct the seed curves
and the corresponding height functions for the more general
characteristic planes \eqref{charplane} in the case $c\not= 0$.
Recalling that such planes are obtained by left-translating $\Pi$
with respect to the point $g_o = (- 2b/c, 2a/c,d/c)$, one easily
recognizes from \eqref{nuXplanes} that, given $z\not= (- 2b/c,
2a/c)$, the seed curve with initial point $z$ is given by
\begin{equation}\label{seedcharplane}
\tilde \gamma_1(s)\ =\ \gamma_1(s)\ -\ \frac{2b}{c}\ ,\quad\quad\quad \tilde \gamma_2(s)\ =\ \gamma_2(s)\ +\ \frac{2a}{c}\ ,
\end{equation}
where $\gamma = (\gamma_1,\gamma_2)$ represents the seed curve in \eqref{seedpi} for the plane $\Pi$ with initial point $(x + 2b/c, y - 2a/c)$. It might be helpful to observe, in this respect, that
\[
|\gamma|\ =\ \frac{2}{c}\ \sqrt{a^2 + b^2 + \frac{c^2}{4} |\tilde \gamma|^2 + b c \tilde \gamma_1 - a c \tilde \gamma_2}\ .
\]
}
\end{Example}

\medskip

\begin{Example}\label{E:hyp2}
\emph{We next analyze the surface
\[
S\ =\ \left\{(x,y,t) \in \mathbb H^1 \mid t = \frac{x y}{2} \right\}\ ,
\]
 given in Example \ref{E:hyp}. Recalling the expression \eqref{nuXhyp} of the horizontal Gauss map, if for instance $z \in \Om^+ = \{(x,y)\in \mathbb R^2 \mid y > 0\}$, then one readily verifies that the seed curve $\gamma = \gamma_z$ is a straight line given by
\begin{equation}\label{E:hyp:g}
\gamma(s)\ =\ (x - s , y)\ ,
\end{equation}
and therefore
\[
\kappa(s)\ \equiv\ 0\ .
\]
From this we infer that the singular locus $\mathcal C_\gamma =
\varnothing$, and in fact
\begin{equation}\label{Fhyp}
F(s,r)\ =\ (x - s, y + r)\ ,
\end{equation}
defines a global diffeomorphism of $\mathbb R^2$ onto itself. In the coordinates $(s,r)$ the defining function $h(x,y) = xy/2$ of $S$ is expressed by
\[
h(s,r)\ =\ \frac{xy + xr - ys - sr}{2}\ ,
\]
and therefore we find in particular
\[
h_0(s)\ =\ h(s,0)\ =\ \frac{y}{2}\ (x - s)\ .
\]
Imposing equation \eqref{remchar} we see that in the plane $(s,r)$ the characteristic locus is given by the straight line $r = - y$. The image of this line on the surface $S$ through the diffeomorphism \eqref{Fhyp} is given by the equation
\[
\left(x -s , y + r, \frac{xy + xr - ys - sr}{2}\right)\ ,
\]
in which we need to set $r = - y$. We thus obtain the collection of points
\[
\Sigma\ =\ \{(x - s , 0 , 0)\in S \mid s\in \mathbb R\}\ =\ \{(s,0,0)\in S \mid s\in \mathbb R\}\ .
\]
This is in accordance with the characteristic locus that we found in Example \ref{E:hyp} in terms of the coordinates $(x,y,t)$.
\\
If, instead, $z = (x,y) \in \Om^- = \{(x,y)\in \mathbb R^2 \mid y < 0\}$, then the corresponding seed curve is the straight line given by
\[
\gamma(s)\ =\ (x + s , y)\ ,
\]
and the relative diffeomorphism and height function by
\[
F(s,r)\ =\ (x + s, y - r)\ ,
\]
\[
h_0(s)\ =\ \frac{y}{2}\ (x + s)\ .
\]
If we impose \eqref{remchar} with these data, we reach the same
conclusion for the characteristic locus.  In figure \ref{eg5.2}, we
show the vector field $\nuX$, a seed curve corresponding to
$(x,y)=(0,1)$ in equation \ref{E:hyp:g} (in black) and some
representative rules (in grey).  We
can see from the picture that the rules are all parallel, which
provides visual confirmation of the fact that the singular locus is
empty.
}
\begin{figure}
\epsfig{figure=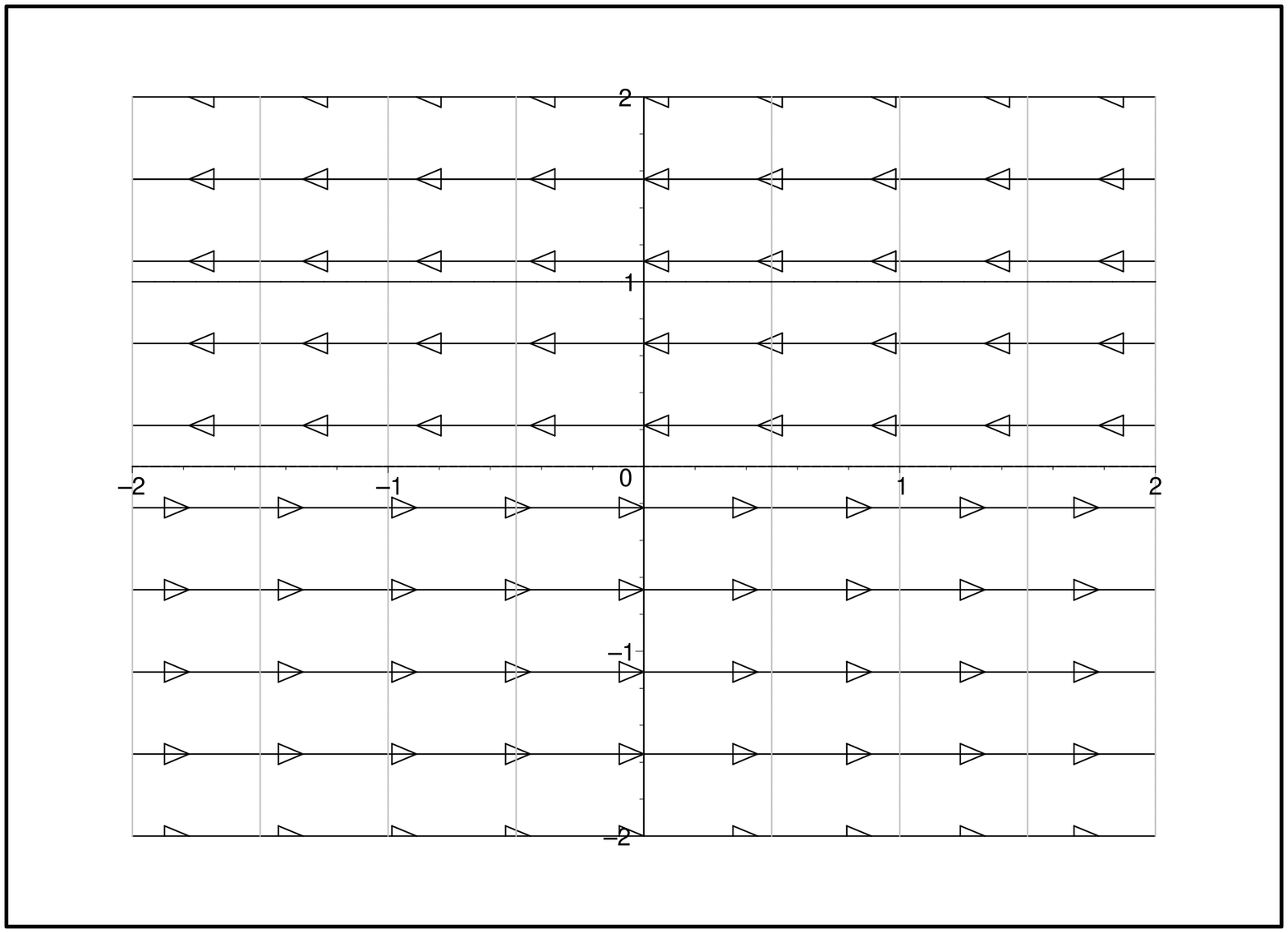,angle=-90,width=3in}
\caption{A seed curve and projections of rules for $\{t = \frac{xy}{2}\}$}\label{eg5.2}
\end{figure}

\end{Example}

\medskip

\begin{Example}\label{E:catenoid2}
\emph{We consider here the sub-Riemannian catenoids in Example \ref{E:catenoid}, in which we choose $u_o = 0$, and $a = 2$. Such choice
yields the surface
\[
S\ =\ \left\{(x,y,t)\in \mathbb H^1 \mid t^2\ =\ \frac{|z|^2}{2} - 1\right\}\ .
\]
We restrict attention to that portion of $S$ that lies above the plane $\{(x,y,t) \in \mathbb H^1\mid t=0\}$, for which the corresponding  horizontal Gauss map is given by \eqref{nuXcat}. Denoting with $\Om = \{z\in \mathbb R^2 \mid |z|>\sqrt 2\}$, we fix a point $z \in \Om$ and compute the corresponding seed curve $\gamma = \gamma_z$. The latter solves the system
\begin{equation}\label{seedcat1}
\begin{cases}
\gamma_1'(s)\ =\ -\ \sqrt{2}\ \left(\frac{\gamma_1(s) + \gamma_2(s) \sqrt{\frac{|\gamma(s)|^2}{2} - 1}}{|\gamma(s)|^2}\right)\ ,\quad\quad\quad \gamma_1(0)\ =\ x\ ,
\\
\gamma_2'(s)\ =\  -\ \sqrt{2}\ \left(\frac{\gamma_2(s) - \gamma_1(s) \sqrt{\frac{|\gamma(s)|^2}{2} - 1}}{|\gamma(s)|^2}\right)\ ,\quad\quad\quad \gamma_2(0)\ =\ y\ .
\end{cases}
\end{equation}
One recognizes that if $\gamma$ is a solution to \eqref{seedcat1}, then
\[
\frac{d}{ds}\ \left(\frac{|\gamma(s)|^2}{2}\right)\ =\ -\ \sqrt 2\ ,
\]
and therefore
\begin{equation}\label{seedcat2}
|\gamma(s)|\ =\ \sqrt{|z|^2 - 2 \sqrt 2 s}\ ,\quad\quad\quad -\infty\ <\ s\ \leq \frac{|z|^2}{2 \sqrt 2}\ .
\end{equation}
This shows that the seed curve of $S$ is similar to a spiral of Archimedes, and we have
\[
\underset{s\to - \infty}{\lim}\ |\gamma(s)|\ =\ +\ \infty\ ,\quad\quad\quad \underset{s\to \frac{|z|^2}{2 \sqrt 2}}{\lim}\ |\gamma(s)|\ =\ 0\ .
\]
Restricting $s$ to the range $(-\infty,0]$, we see that $|\gamma(s)|^2 \geq 2$, and therefore it is legitimate to substitute \eqref{seedcat2} into \eqref{seedcat1}, and the system reduces to the linear one
\begin{equation}\label{seedcat3}
\begin{cases}
\gamma_1'(s)\ =\ -\ \sqrt{2}\ \left(\frac{\gamma_1(s) + \gamma_2(s) \sqrt{\frac{|z|^2}{2} - \sqrt 2 s - 1}}{|z|^2 - 2 \sqrt 2 s}\right)\ ,\quad\quad\quad \gamma_1(0)\ =\ x\ ,
\\
\gamma_2'(s)\ =\  -\ \sqrt{2}\ \left(\frac{\gamma_2(s) - \gamma_1(s) \sqrt{\frac{|z|^2}{2} - \sqrt 2 s - 1}}{|z|^2 - 2 \sqrt 2 s}\right)\ ,\quad\quad\quad \gamma_2(0)\ =\ y\ .
\end{cases}
\end{equation}
In figure \ref{eg5.3} we again show $\nuX$, a seed curve in black and
some rules in grey for one of these surfaces, $t^2=|z|^2-4$.  Note that the vector field is not defined inside
the circle of radius $\sqrt{2}$.  We also note that this picture
corresponds to the seed data for the upper half of the surface, the
graph $t=\sqrt{|z|^2-4}$.  There is a corresponding picture for
the lower half resulting in a similar seed curve for the bottom portion of this surface.  We note that, in this case, we need two seed curves to describe the surface completely
 as the surface has two sheets when written as a graph over the $xy$-plane.  As with Example \ref{E:gencurve}, this justifies the necessity of the notion of generalized seed curve
that we introduce in Section \ref{S:characterization}.
}
\begin{figure}
\epsfig{figure=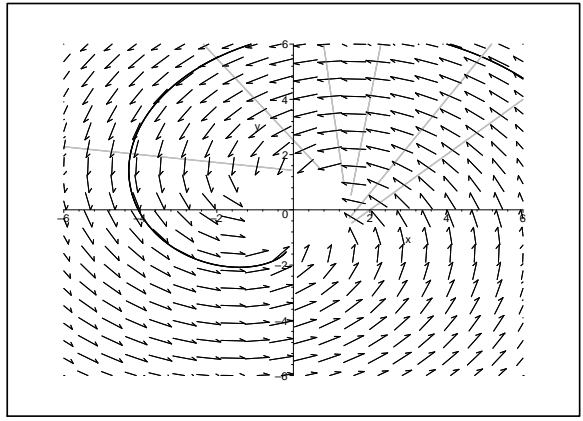,angle=-90,width=3in}
\caption{A seed curve and projections of rules for  $\{t=\sqrt{\frac{|z|^2}{2}-1}\}$}\label{eg5.3}
\end{figure}
\end{Example}

\medskip

\begin{Example}\label{E:ceseed}
\emph{We next analyze the surface $S$ in the Example \ref{E:ce}. To compute a seed curve we will use the parameterization \eqref{hor} of $S$,
\begin{equation}\label{hor2}
t\ =\ h(x,y)\ \overset{def}{=}\ -\  \tanh^{-1}\left(\tan^{-1}\left(\frac{y}{x}\right)\right)\ ,
\end{equation}
which is valid over the open set $\Om = \{(x,y) \in \mathbb R^2 \mid x \not= 0\}$. Recalling the corresponding expression \eqref{nuXce} of the horizontal Gauss map, we see that for every $z\in \Om$ the seed curve starting at $z$ must satisfy the system
\begin{equation}\label{ceseed1}
\begin{cases}
\gamma_1'(s)\ =\ -\ \frac{\left(\frac{1}{|\gamma(s)|^2 ( 1 - \alpha^2)} + \frac{1}{2}\right)\ \gamma_2(s)}{\sqrt{\frac{1}{|\gamma(s)|^2 (1 - \alpha^2)^2} + \frac{|\gamma(s)|^2}{4} + \frac{1}{1 - \alpha^2}}}\ ,\quad\quad\quad \gamma_1(0)\ =\ x\ ,
\\
\gamma_2'(s)\ =\  \quad \frac{\left(\frac{1}{|\gamma(s)|^2 ( 1 - \alpha^2)} + \frac{1}{2}\right)\ \gamma_1(s)}{\sqrt{\frac{1}{|\gamma(s)|^2 (1 - \alpha^2)^2} + \frac{|\gamma(s)|^2}{4} + \frac{1}{1 - \alpha^2}}}\ ,\quad\quad\quad \gamma_2(0)\ =\ y\ .
\end{cases}
\end{equation}
As in the previous example we easily recognize that if $\gamma(s)$ solves \eqref{ceseed1}, then it must be $|\gamma(s)| \equiv const = |z|$. With this information in hands, \eqref{ceseed1} reduces to
\begin{equation}\label{ceseed2}
\begin{cases}
\gamma_1'(s)\ =\ -\ \frac{\left(\frac{1}{|z|^2 ( 1 - \alpha^2)} + \frac{1}{2}\right)\ \gamma_2(s)}{\sqrt{\frac{1}{|z|^2 (1 - \alpha^2)^2} + \frac{|z|^2}{4} + \frac{1}{1 - \alpha^2}}}\ ,\quad\quad\quad \gamma_1(0)\ =\ x\ ,
\\
\gamma_2'(s)\ =\  \quad \frac{\left(\frac{1}{|z|^2 ( 1 - \alpha^2)} + \frac{1}{2}\right)\ \gamma_1(s)}{\sqrt{\frac{1}{|z|^2 (1 - \alpha^2)^2} + \frac{|z|^2}{4} + \frac{1}{1 - \alpha^2}}}\ ,\quad\quad\quad \gamma_2(0)\ =\ y\ .
\end{cases}
\end{equation}
It is easily recognized that
\[
\frac{\frac{1}{|z|^2 ( 1 - \alpha^2)} + \frac{1}{2}}{\sqrt{\frac{1}{|z|^2 (1 - \alpha^2)^2} + \frac{|z|^2}{4} + \frac{1}{1 - \alpha^2}}}\ \equiv \ \frac{1}{|z|}\ ,
\]
and therefore \eqref{ceseed2} further reduces to
\begin{equation}\label{ceseed3}
\begin{cases}
\gamma_1'(s)\ =\ -\ \frac{1}{|z|}\ \gamma_2(s)\ ,\quad\quad\quad\quad \gamma_1(0)\ =\ x\ ,
\\
\gamma_2'(s)\ =\  \quad \frac{1}{|z|}\ \gamma_1(s)\ ,\quad\quad\quad\quad \gamma_2(0)\ =\ y\ .
\end{cases}
\end{equation}
Incidentally, although the function $\alpha = \alpha(x,y) = \tan^{-1}(y/x)$ no longer plays a role in \eqref{ceseed3}, we note that its expression in the new coordinates $(s,r)$ is given by
\begin{align*}
\alpha\ & =\ \tan^{-1}\left(\frac{\gamma_2(s)}{\gamma_1(s)}\right)\ =\ \tan^{-1}\left(\frac{y + x \tan(s/|z|)}{x - y \tan(s/|z|)}\right)
\\
& =\ \tan^{-1}\left(\frac{\tan(\arg(z)) + \tan(s/|z|)}{1 - \tan(\arg(z)) \tan(s/|z|)}\right)
\\
& =\ \tan^{-1} \tan\left(\arg(z) + \frac{s}{|z|})\right)\ =\ \arg(z) + \frac{s}{|z|}\ .
\end{align*}
The solution to \eqref{ceseed3}, the seed curve $\gamma = \gamma_z$, is given by the circle
\[
\begin{cases}
\gamma_1(s)\ =\ x\ \cos\left(\frac{s}{|z|}\right)\ -\ y\ \sin\left(\frac{s}{|z|}\right)\ ,
\\
\gamma_2(s)\ =\ y\ \cos\left(\frac{s}{|z|}\right)\ +\ x\ \sin\left(\frac{s}{|z|}\right)\ .
\end{cases}
\]
From \eqref{kappa} we obtain, as in Example \ref{E:comp},
\begin{equation}\label{kappaex}
\kappa(s)\ \equiv\ -\ \frac{1}{|z|}\ ,
\end{equation}
with singular locus given by the straight line
\[
\mathcal C_\gamma\ =\ \left\{(s,r)\in \mathbb R^2 \mid r = - |z|\right\}\ ,
\]
Outside of $\mathcal C_\gamma$ the map
\begin{align}\label{point}
F(s,r)\ & =\ \left(1 + \frac{r}{|z|}\right)\ \gamma(s)
\\
& =\ \left(1 + \frac{r}{|z|}\right)\ \left(x\  \cos\left(\frac{s}{|z|}\right) - y\sin\left(\frac{s}{|z|}\right)\ ,\ y\ \cos\left(\frac{s}{|z|}\right) + x\sin\left(\frac{s}{|z|}\right)\right)\ ,
\notag\end{align}
defines a local diffeomorphism, similarly to Example \ref{E:comp}. Finally, we note that, in the coordinates $(s,r)$, the function $h$ in \eqref{hor2}
is given by
\[
h(s,r)\ =\ -\ \tanh^{-1}(\alpha(s,r))\ =\ -\ \tanh^{-1}\left(\arg(z) + \frac{s}{|z|}\right)\ ,
\]
thus, in particular,
\[
h_0(s)\ =\ h(s,0)\ =\ -\ \tanh^{-1}\left(\arg(z) + \frac{s}{|z|}\right)\ ,
\]
for
\[ -\ (1 + \arg(z))\ |z|\ <\ s\ <\ (1 - \arg(z))\ |z|\  .
\]
We now use \eqref{remchar} to verify, in conformity with what we found in the discussion of Example \ref{E:ce}, that the characteristic locus of $S$ is empty. Since
\[h_0'(s)\ =\ -\ \frac{1}{|z|\left[ 1 - \left(\arg(z) + \frac{s}{|z|}\right)^2\right]}\ ,
\]
\eqref{remchar} presently reads
\[
(r + |z|)^2\ +\ \frac{2}{ 1 - \left(\arg(z) + \frac{s}{|z|}\right)^2}\ =\ 0\ .
\]
It is clear that the latter equation has no solution in the region
\[
\{(s,r)\in \mathbb R^2\mid - (1 + \arg(z)) |z| < s < (1 - \arg(z)) |z|\ ,\ - \infty < r < \infty \}\ ,
\]
and therefore $\Sigma = \varnothing$. For instance, if $z = (1,0)$, then \eqref{kappaex} gives $\kappa(s) \equiv - 1$, and we find $h(s,r) = - \tanh^{-1}(s)$, with
\[
h_0'(s)\ =\ -\ \frac{1}{1 - s^2}\ .
\]
In figure \ref{eg5.4} we show a seed curve (in black) and some rules
(in grey) for this example.  We note that all the rules intersect at
the origin, reflecting the fact that, in this case, the singular locus maps to a point in the $xy$-plane, the origin $(0,0)$, see equation \eqref{point}.  By analyzing the equation of the surface, we see
that, near the origin in the $xy$-plane, the surface ceases to be a
graph over the $xy$-plane.  We also note that there are portions of this
plane where $\nuX$ is not defined because the height function,
$h_0(s)$, tends to $\pm \infty$ for $s\to \mp 1$.}

\begin{figure}
\epsfig{figure=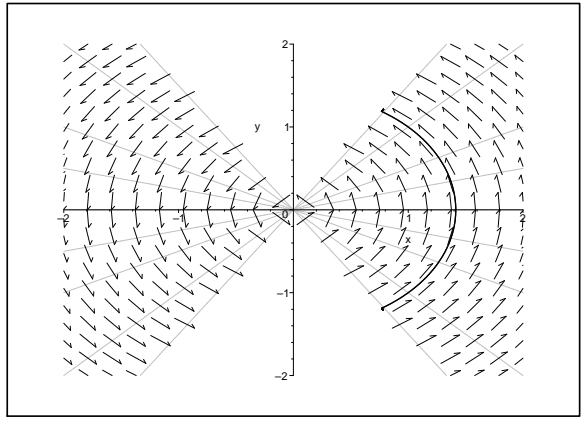,angle=-90,width=3in}
\caption{A seed curve and projections of rules for $\{y=-x\tan(\tanh t)\}$}\label{eg5.4}
\end{figure}

\end{Example}

\medskip

We next present two examples which demonstrate some of the
restrictions on regularity of $H$-minimal surfaces. We show how to
construct a $C^1$ $H$-minimal surface which is not $C^2$, and whose characteristic locus is not smooth.

\medskip

\begin{Example}\label{E:optreg}
\emph{Let $\psi:[-1,1] \subset \R \ra \R$ be a strictly positive, continuous, nowhere differentiable function, and
let $\Psi(s)=\int_{-1}^s \psi(t) \; dt$.  We define a curve in the plane by letting
\[
\gamma(s)\ =\ \left ( \int_{-1}^s \cos(\Psi(t))\;dt,\int_{-1}^s \sin(\Psi(t))\;dt
\right )\ .
\]
The function $\gamma$ gives a $C^2$ curve parameterized by arc-length which we use as a
seed curve for an $H$-minimal surface, $S$.  Note that
\[
\gamma'(s)\ =\ (\cos(\Psi(s)),\sin(\Psi(s)))\ ,
\]
and
\[
\gamma''(s)\ =\ (-\psi(s)\sin(\Psi(s)),\psi(s)\cos(\Psi(s)))\ ,
\]
and so \eqref{kappa} gives
\[
\kappa(s)\ =\ -\ \psi(s)\ .
\]
We point out that, since
$\gamma \in C^2$, if we choose $h_0(s)\in C^1$, the computations in
the second half of the proof of Theorem \ref{rep} apply and this seed
curve and height function pair define an H-minimal graph.
Moreover, using \eqref{remchar}, after picking $h_0(s)$, we can explicitly write down the
characteristic locus.  For example, if we take
\[
h_0'(s)\ =\ -\ \frac{1}{2}<\gamma'(s),\gamma(s)^\perp>\ +\ 1\ ,
\]
then \eqref{remchar} gives for the characteristic locus
\[
r\ =\ -\ \frac{1}{\psi(s)}\ \pm\ \frac{\sqrt{1+2\psi(s)}}{\psi(s)}\ .
\]
As $\psi$ is nowhere differentiable, we see that the characteristic locus
is also nowhere differentiable.  Moreover, by construction, the patch
of $S$ given by the equations \eqref{para}, \eqref{para2}, defined on the open set
\[
(s,r)\ \in\ (-1,1) \times \left (-\infty,
  -\frac{1}{\psi(s)}\right )\,
\]
is a $C^1$ graph over the $xy$-plane
containing one branch of the characteristic locus, $r = - \frac{1}{\psi(s)} - \frac{\sqrt{1+2\psi(s)}}{\psi(s)}$}.
\end{Example}

\medskip

\begin{Example}\label{E:optreg2}
\emph{As a more concrete situation, we will consider the case where
$\psi(s)=|s|$ in Example \ref{E:optreg}.  Of course, such $\psi$ is differentiable everywhere except at
$s=0$, but we will still produce a characteristic curve of the corresponding $H$-minimal surface $S$ which is not
differentiable at one point.  Using the construction above, we have that in the $(s,r)$ plane the characteristic locus is
given by the curve
\[
r\ =\ -\ \frac{1}{|s|}\ \pm\ \frac{\sqrt{1+2|s|}}{|s|}\ .
\]
Taking the branch
\[
r\ =\ -\ \frac{1 - \sqrt{1 + 2|s|}}{|s|}\ ,
 \]
\begin{figure}
\epsfig{figure=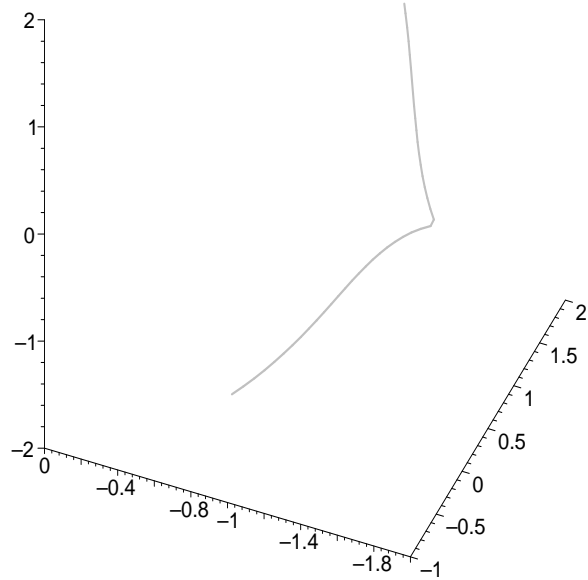,width=3in,height=3in}
\caption{The characteristic locus for Example 6.6}\label{5.6-2}
\end{figure}
we see that the characteristic locus is not differentiable at $s=0$.
Indeed, letting $c(s)=-\frac{1-\sqrt{1+2|s|}}{|s|}$, we have $c(s) = 1 + o(1)$ as $s\to 0$, and therefore if we set $c(0) = 1$, we find
\[
\lim_{s \ra 0^+}\frac{c(s)+1}{s}\ =\ -\ \frac{1}{2}\ ,
\]
and
\[
\lim_{s \ra 0^-}\frac{c(s)+1}{s}\ =\ \frac{1}{2}\ .
\]
The characteristic set of the $H$-minimal surface $S\subset \mathbb H^1$ is given by the image on $S$, $\tilde c(s)$, via the parameterization \eqref{para}, \eqref{para2}, of the curve $r = c(s)$,
\[
\tilde{c}(s)\ =\ \left (\gamma_1(s)+c(s)\gamma_2'(s),\gamma_2(s)-c(s)\gamma_1'(s),
  h_0(s) - \frac{c(s)}{2}<\gamma(s),\gamma'(s)>\right)\ .
\]
Direct computation shows}

\[\lim_{h \ra 0^+}\frac{\tilde{c}(h)+\tilde{c}(0)}{h}-\lim_{h \ra
  0^-}\frac{\tilde{c}(h)+\tilde{c}(0)}{h}= \left (-\sin\left (\frac{1}{2}\right),
  \cos\left(\frac{1}{2}\right),
  \frac{\sqrt{\pi}}{2}\sin^2\left(\frac{1}{2}\right)C_0+\frac{\sqrt{\pi}}{2}\cos^2\left(\frac{1}{2}\right)C_0 \right)\]
where \[C_0=\int_0^\frac{1}{\sqrt{\pi}}\cos\left(\frac{\pi}{2}t^2\right)
\; dt\ .
\]
\emph{Figure \ref{5.6-2} shows a picture of the characteristic locus - one
can see the point of non-differentiability.  Figure \ref{5.6-4}
shows a portion of the surface with the image of the seed
curve in black and the characteristic locus in grey.}

\begin{figure}
\epsfig{figure=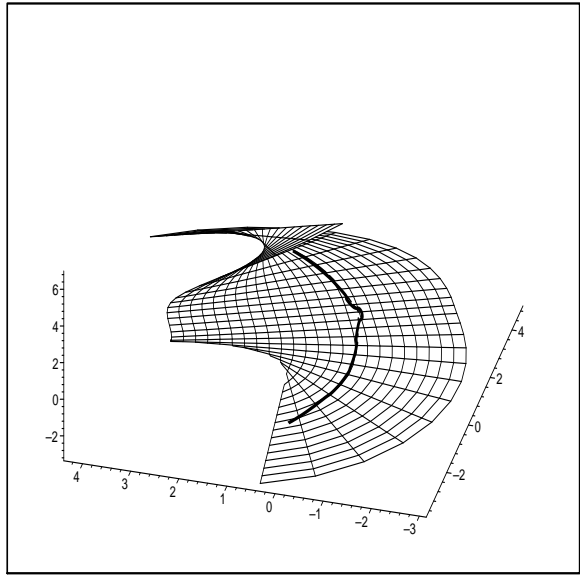,width=3in,height=3in}
\caption{A picture of Example 6.6 with characteristic locus denoted by
  the heavy black line.}\label{5.6-4}
\end{figure}

\end{Example}

\begin{Example}  \label{E:cyl}
\emph{In most of the previous examples, we examined surfaces first discussed
in Section \ref{S:examples} as graphs over a portion of a coordinate plane, in terms of their seed curves and height
functions.  While this brings a new perspective in the analysis of known $H$-minimal
surfaces, here we wish to emphasize the usefulness of Theorem \ref{rep} in
the creation of new examples.  Specifically, by simply picking a seed
curve and a height function, the representation given in Theorem
\ref{rep} yields an $H$-minimal surface.  To illustrate this, we will
address the following question:  do there exist $H$-minimal surfaces
that are topological cylinders with piecewise constant horizontal
Gauss map?  Here, we mean that $\nuX$, thought of as a vector field on
the plane, is piecewise constant.  Of course, Example \ref{E:catenoid}
gives an $H$-minimal surface which is a topological cylinder, but its
unit horizontal Gauss map is not piecewise constant in this sense.
\\
To answer this question, we recall that in the proof of Theorem
\ref{rep}, we showed that given an $H$-minimal graph over the $xy$-plane, $S$, with seed
curve $\gamma(s)$ and height function $h_0(s)$, the unit horizontal
Gauss map at a point $(s,r)$ is given by the vector $\gamma'(s)$.
Thus, if we wish to create a surface with constant unit horizontal
Gauss map, we need to pick $\gamma$ to be a straight line.  So, for any portion
of the desired surface which is a graph over the $xy$-plane, we must have
that the seed curve is a straight line.
\\
Now, to form a topological cylinder, we will attempt to glue two such
$H$-minimal graphs, $S_1$ and $S_2$, together by picking appropriate initial height
functions.  Experimentation yields the following solution:  for $S_1$ we take $\gamma^1(s)=(s,0)$ and $h_0^1(s)=\sqrt{1-s^2}$, while for $S_2$ we choose
$\gamma^2(s)=(s,0)$ and $h_0^2(s)=-\sqrt{1-s^2}$.  In other words, we
pick straight lines for seed curves and lift them using appropriate
height functions so that when glued together they form a circle.
Applying the representation formula \eqref{para} in Theorem \ref{rep}, we have
that $S_1$ is given by
\[\left (s,-r,\sqrt{1-s^2}+\frac{rs}{2} \right ) \;\; \text{for $s \in
  [-1,1]\ ,\quad\quad  r \in \R$}\ ,
\]
and $S_2$ is given by
\[\left (s,-r,-\sqrt{1-s^2}+\frac{rs}{2} \right ) \;\; \text{for $s \in
  [-1,1]\ ,\quad\quad r \in \R$}\ .
\]
We may now use \eqref{para2} to algebraically simplify the latter two equations to find that the
surface $S=S_1 \cup S_2$ is implicitly given by
\[ \left (t-\frac{xy}{2}\right)^2=1-x^2\ .
\]
We notice that $S$ is contained in the slab $\{(x,y,t)\in \mathbb H^1 \mid -1 \leq x \leq 1\}$, and that its characteristic locus is given by the
curve $\Sigma = \{(x,y,t)\in S \mid x = y (t - \frac{xy}{2})\}$.
If we let
$S^- = \{(x,y,t)\in S \mid x < y (t - \frac{xy}{2})\}$,  $S^+ = \{(x,y,t)\in S \mid x > y (t - \frac{xy}{2})\}$, then
we have $\nuX = (-1,0)$ on $S^-$, $\nuX = (1,0)$ on $S^+$. }
\end{Example}

\vskip 0.6in


\section{\textbf{Characteristic points on $H$-minimal graphs over the
  $xy$-plane}} \label{S:charpoints}

\vskip 0.2in

We turn next to understanding how characteristic points can arise on
$H$-minimal graphs. From the proof of Theorem \ref{rep},
we see that for a $C^2$ $H$-minimal surface which is a graph over a portion of the
$xy$-plane, the characteristic locus corresponds
to the equation
\begin{equation}\label{cl0}
h_0'(s)\ -\ r\ -\ \frac{1}{2} <\gamma'(s),\gamma(s)^\perp>\ +\ \frac{r^2}{2}\ \kappa(s)\ =\ 0\ .
\end{equation}

This formula, while explicit, is rather cumbersome.  In this section, we will provide a slightly different perspective from which
to view the characteristic locus which gives a different and
admittedly only slightly less cumbersome formula, see Theorem \ref{T:alph} and Remark \ref{R:rems}.  However, this
second formula will be much more useful in determining the
characteristic locus.

Under the same assumptions outlined at the beginning of Section \ref{S:representation}, we note that the (non-unit) Riemannian normal to the surface $S$ can be written
as
\begin{equation}\label{N}
\boldsymbol N\ =\ W\ \nuX\ +\   T\ =\ W \; \pb \; X_1\ +\ W \; \qb
\; X_2\ +\ T\ ,
\end{equation}
see formula \eqref{run}. We recall that the angle function
$W=\sqrt{p^2+q^2}$ is independent of $t$, as we assume that $S$ is
a graph over the $xy$-plane. The characteristic locus corresponds
to those points where $W(x,y) = 0$.  The formula \eqref{cl0} gives
one representation for such locus. For a different derivation, we
use the Fr\"{o}benius condition to impose a restriction on $W$.
First, we observe that the $2$-plane field perpendicular to
$\boldsymbol N$ is spanned by
\[
V_1\ =\ X_1\ -\ W\ \pb\ T\ ,\quad\quad\quad V_2\ =\ X_2\ -\ W\ \qb\ T\ .
\]

To verify the assumption in Fr\"{o}benius' theorem, we must impose the condition that
$[V_1,V_2] = A V_1 + B V_2$ for some $\{A,B\}$. Now,
\begin{equation}\label{frobenius}
[V_1,V_2]\ =\ -\ \big\{X_1(W \qb)\ -\ X_2(W \pb)\ -\ 1\big\}\ T\ .
\end{equation}

On the other hand, we have
\begin{equation}\label{lincomb}
A\ V_1\ +\ B\ V_2\ =\ A X_1\ +\ B X_2\ -\ (W \bar p A + W \bar q B)\ T\ .
\end{equation}

Equating \eqref{frobenius} and \eqref{lincomb} we find that we must have
$A=B=0$, and therefore,
\begin{equation}\label{nonpar1}
X_1(W \qb)\ -\ X_2(W \pb)\ -\ 1\ =\ 0\ .
\end{equation}

The fact that $W, \pb, \qb$ are functions of $x$ and $y$ only allows us to simplify equation (\ref{nonpar1}) as follows
\begin{equation}
\qb\ W_x\ -\ \pb\ W_y \ +\ W\ (\qb_x-\pb_y)\ -\ 1\ =\ 0\ ,
\end{equation}
or,
\begin{equation}\label{nonpar2}
<\nabla W, \nuX^\perp> \ =\ 1\ -\ W(\qb_x - \pb_y)
\end{equation}

In other words, the directional derivative of $W$ in the
direction of $\nuX^\perp$ is $1 - W(\qb_x-\pb_y)$.
Since $\nuX = \pb X_1 + \qb X_2$, we see from \eqref{goodnuX} that
\[
(\pb,\qb)\ =\  (\gamma_1'(s),\gamma_2'(s))\ .
\]
Using this information, and the equations \eqref{ps_xy}, we obtain
\begin{equation}\label{p&q}
\begin{split}
\qb_x\ -\ \pb_y\ &=\ \gamma_2''(s) \frac{\partial s}{\partial x}\ -\ \gamma_1''(s)
\frac{\partial s}{\partial y}\\
&=\
\frac{\gamma_2''(s)\gamma_1'(s)\ -\ \gamma_1''(s)\gamma_2'(s)}{1 - r \kappa(s)}\\
&=\ -\ \frac{\kappa(s)}{1-r\kappa(s)}\ ,
\end{split}
\end{equation}
where in the latter equality we have used \eqref{kappa}.
Substituting \eqref{p&q} in \eqref{nonpar2}, we conclude
\begin{equation}\label{nablaW}
<\nabla W ,\nuXp>\ =\ 1\ +\ W
\frac{\kappa(s)}{1-r\kappa(s)}\ .
\end{equation}

We now recall the local diffeomorphism introduced in \eqref{newpar}. In terms of the coordinates $(s,r)$ we see that
\[
W_r(s,r)\ =\ <\nabla W(s,r),\nuX^\perp(\gamma(s))>\ .
\]

If we use this in combination with \eqref{nablaW}, we find
\begin{equation}\label{deqW}
W_r(s,r)\ =\ 1\ +\ \frac{\kappa(s)}{1-r\kappa(s)}\ W(s,r)\ .
\end{equation}

Solving this differential equation, one obtains the following result.

\medskip

\begin{Thm}\label{T:alph} If the Riemannian normal to a $C^2$ $H$-minimal
  graph over the $xy$-plane  is given by \eqref{N}, then
\[
W(s,r)\ =\
\frac{W_0(s) + r - \frac{r^2}{2} \kappa(s)}{1 - r \kappa(s)}\ ,
\]
where $W_0(s)=W(s,0)$. Moreover, for $\kappa(s) \not= 0$ the characteristic locus
of the $H$-minimal surface is given by
\[
\Sigma\ =\ \left\{(s,r)\in \mathbb R^2 \mid W_0(s)\ +\ r\ -\ \frac{r^2}{2} \kappa(s)\ =\ 0\right\}\ .
\]
\end{Thm}

\medskip

\begin{Rmk}\label{R:rems}  Several comments are in order:
\begin{enumerate}
\item The functions $W_0(s)$ and $h_0(s)$ are related by
\[
W_0(s)\ =\ -\ h_0'(s)\ +\ \frac{1}{2} <\gamma'(s), \gamma(s)^\perp>\ ,
\]
see equation \eqref{cl0}.
\item If $\kappa(s_0)=0$ for some $s_0$, then the characteristic locus of the surface $S$ {\em must}
contain the point corresponding via the local diffeomorphism $F(s,r)$ to the point $(s_0,r_0)$ in the $sr$-plane, where
\[
r_0\ =\ -\ W_0(s_0)\ .
\]
\item If $W_0(s_0)\kappa(s_0) > - \frac{1}{2}$, the
  characteristic locus has two components
\[
r\ =\ \frac{1}{\kappa(s_0)}\ \pm\ \frac{\sqrt{1 + 2 W_0(s_0)\kappa(s_0)}}{\kappa(s_0)}\ ,
\]
one on each side of the point $(s_0,r_0)$, with
\[
r_0\ =\ \frac{1}{\kappa(s_0)}\ .
\]
\item If $W_0(s_0)\kappa(s_0) < - \frac{1}{2}$, then the surface has
  no characteristic points along the rule passing through $\gamma(s_0)$.
\item If $W_0(s_0) = - \frac{1}{2 \kappa(s_0)}$, then there is a single
  characteristic point at $r = \frac{1}{\kappa(s_0)}$ along the rule
  passing through $\gamma(s_0)$.
\end{enumerate}
\end{Rmk}

\medskip

\begin{Thm}\label{T:optregchar}   Suppose $k\ge 2$.  If the characteristic locus contains more than a single
  point, then the characteristic locus of a $C^k$ $H$-minimal surface
  is a
  $C^{k-1}$ curve.
\end{Thm}

\pf  To determine the regularity of the characteristic locus, we
first suppose that $\Sigma$ contains more than one point.  If
$\nu$ is the Riemannian unit normal to the surface, we may write
the (closed) condition that defines characteristic points as
$<\nu,T>=1$ (the reader should keep in mind that we have endowed
$\mathbb H^1$ with a left-invariant Riemannian metric with respect
to which $\{X_1, X_2, T\}$ constitute an orthonormal basis). So,
in a sufficiently small neighborhood of a characteristic point $S$
must be a graph over the $xy$-plane. Deleting the characteristic
points from this neighborhood, we have an open non-characteristic
patch of the surface and may use Theorem \ref{rep} to parameterize
the patch using seed curves and height functions.   Using Theorem
\ref{T:alph}, and (1) of Remark \ref{R:rems}, we have that the
characteristic locus is given by
\begin{equation}\label{cl}
-h_0'(s) - \frac{1}{2}<\gamma'(s),\gamma^\perp(s)> +\ r -\frac{r^2}{2}\kappa(s)\ =\ 0\ ,
\end{equation}
where $\kappa(s)=<\gamma'',(\gamma')^\perp>$.  As in the remarks
after Theorem \ref{T:alph}, we can describe the characteristic
locus as $r=f(s)$.  Using the implicit function theorem, we know
that this curve is differentiable if equation (\ref{cl}) is
differentiable in $s$.  By Theorem \ref{rep}, we know that if
$S\in C^k$ then  $\gamma \in C^{k+1}, h_0 \in  C^k$, and we have
that the right hand side of (\ref{cl}) is $k-1$ times
differentiable in $s$ as desired.  $\qed$

\begin{Rmk}  The statement that a $C^2$ $H$-minimal surfaces with
  characteristic locus larger than a single point is $C^1$ is also shown in
 Theorem B in \cite{CHMY}, although using different techniques.
\end{Rmk}

\vskip 0.6in

\section{{\textbf{Allowable gluing of patches of complete embedded $C^2$ $H$-minimal surfaces}}}\label{S:glue}

\vskip 0.2in

 In the previous sections, we have discussed a local
parameterization of $H$-minimal surfaces near non-characteristic
points.  In this section, we wish to investigate how two such
patches are allowed to meet under the assumption that they belong
to the same complete embedded $C^2$ $H$-minimal surface
$\tilde{S}$. Consider a neighborhood $U$ of a non-characteristic
point determined by a seed curve, $\gamma(s)$, and height
function, $h_0(s)$, as in Theorem \ref{rep}.
\[
S\ =\ \left \{ \left( F(s,r),h_0(s)-\frac{r}{2} <\gamma(s),
\gamma'(s)>
    \right ) \mid s \in (s_0,s_1), r \in (r_0(s),r_1(s)) \right \}\ .
\]

 For a point $g_o$ in the boundary of
$U$, there are two possibilities:
\begin{enumerate}
\item $g_o$ is a limit point of a rule $\mathcal L(r)$ as in \eqref{rule}:  i.e., there is a fixed $s$ and
  $i \in \{0,1\}$ so that
\[ g_o = \lim_{r \ra r_i(s)} \left( F(s,r),h_0(s)-\frac{r}{2} <\gamma(s) ,
\gamma'(s)>
    \right )\ ;
    \]
\item $g_o$ is a limit point of the lift to $\tilde S$ of an integral curve of $\til$:
  i.e., there is a fixed $r$ and  and
  $i \in \{0,1\}$ so that
\[ g_o = \lim_{s \ra s_i} \left( F(s,r),h_0(s)-\frac{r}{2} <\gamma(s) ,
\gamma'(s)>
    \right )\ .
    \]
\end{enumerate}

Each of the above limits is assumed to be the appropriate sided
limit depending on whether $i=0$ or $i=1$.

We consider the possibilities for each case. We recall from
\eqref{N} that the (non-unit) Riemannian normal to $U$ is given by
\begin{equation}\label{rn}
\boldsymbol N\ =\ W(s,r)\ \nuX\ +\ T\ ,
\end{equation}
where
\begin{equation}\label{E:t1}
W(s,r)\ =\ \frac{W_0(s)-r+\frac{r^2}{2}\kappa(s)}{1-r\kappa(s)}\ .
\end{equation}

Next, we want to calculate the Euclidean normal $\boldsymbol n$ to
the surface as well. From the equation \eqref{rn} one immediately
sees that the vector fields
\[
V_1\ =\ X_1\ -\ W(s,r)\ \pb \; T
\]
and
 \[
V_2\ =\ X_2\ -\ W(s,r)\ \qb \ T
\]
are tangential to the surface. Using \eqref{coord}, in combination with \eqref{np}, we can re-write these vector fields in terms of the
standard Euclidean basis, obtaining
\[
V_1\ =\ \partial_x\ -\ \left(W(s,r)\ \pb\ +\
\frac{1}{2}\left(\gamma_2(s)- r\ \gamma_1'(s)\right)\right)\
\partial_t\ ,
\]
and
\[
V_2\ =\ \partial_y\ -\ \left(W(s,r)\ \qb
-\frac{1}{2}\left(\gamma_1(s) + r\ \gamma_2'(s)\right)\right)\
\partial_t\ .
\]

The latter two formulas give for the Euclidean normal $\boldsymbol
n = V_1 \times V_2$
\begin{equation}\label{enorm}
\boldsymbol n(s,r)\ =\ \left(W(s,r)\ \pb\ +\ \frac{\gamma_2(s)- r
\gamma_1'(s)}{2}\ ,\  W(s,r)\ \qb\ -\ \frac{(\gamma_1(s) + r
\gamma_2'(s))}{2}\ ,\ 1\right)\ .
\end{equation}

Now, cases (1) and (2) above are each subdivided in three
subcases. If $g_o$ is a limit point of a rule in $U$, we see three
possibilities:
\begin{enumerate}
\item[($A_1$)] $g_o$ is a non-characteristic point that lies in a
  neighborhood that can be written as a graph over the $xy$-plane.
\item[($B_1$)] $U$ and $r_i(s)$ satisfy one of Remark \ref{R:rems} (2), (3) and
  $g_o$ is a characteristic point.
\item[($C_1$)] $r_i(s)= \frac{1}{\kappa(s)}$ and $g_o$ is in the image of the
  singular locus of the parameterization given by $F$.  In this case,
  $g_o$ is either a characteristic point (by Remark \ref{R:rems} (5))
  or is a point where $\tilde{S}$ ceases to be a graph over the
  $xy$-plane as $W(s,r) \ra \pm \infty$ and hence $\boldsymbol n$ tends to
  a vector with no $\partial_t$ component.
\end{enumerate}

If $g_o$ is a limit point of the lift to $\tilde S$ of an
integral curve of $\til$, we see similar possibilities:

\begin{enumerate}
\item[($A_2$)] If $W(s,r)$ tends to a finite number as $s \ra s_i$, then by
  \eqref{rn} $g_o$ is a non-characteristic point of $\tilde{S}$ that lies in a
  neighborhood that can be written as a graph over the $xy$-plane.
\item[($B_2$)] If $W(s,r)$ tends to zero as $s \ra s_i$, then $g_o$ is a
  characteristic point.
\item[($C_2$)] If $W(s,r)$ tends to $\pm \infty$ as $s \ra s_i$ for a
  fixed $r \neq \lim_{s \ra s_i \frac{1}{\kappa(s)}}$, then, as before,
  we have that $\boldsymbol n$ tends to a vector with no $\partial_t$
  component and hence $g_o$ is a point where $\tilde{S}$ ceases to be
  a graph over the $xy$-plane.
\end{enumerate}

Thus, there are in total three types of possible boundary points:
\begin{enumerate}
\item[(A)] non-characteristic points that lie in surface neighborhoods which are
  graphs over the $xy$-plane;
\item[(B)] characteristic points;
\item[(C)] points where $\tilde{S}$ ceases to be a graph over the
$xy$-plane\ .
\end{enumerate}

We will examine each of these cases separately but wish to point
out that the same basic idea guides the treatment of all cases:
the assumption that the surface be $C^2$ forces the normal to be
$C^1$. By examining the behavior of the normal vector (or, in some
cases, of the horizontal Gauss map) near the boundary points, we
are able to discern the manner in which different patches of
surface are allowed to meet.

In the subsection \ref{ss1}, we will address cases ($A_1$) and
($A_2$).  In subsection \ref{ss2}, we analyze the behavior near
characteristic points, cases ($B_1$) and ($B_2$).  Finally, in
subsection \ref{SS:vert} we study the cases ($C_1$) and ($C_2$),
points at which the surface ceases to be a graph over the
$xy$-plane.

\subsection{Non-characteristic boundary points in graphical neighborhoods}\label{ss1}

Suppose we have a neighborhood $U$ of a non-characteristic point
that is parameterized by \eqref{para} and \eqref{para2} via
Theorem \ref{rep}, and suppose further that we are considering a
point in $\partial U$ which is also non-characteristic.  Then,
again by Theorem \ref{rep}, there is a neighborhood of the
boundary point, $U'$, parameterized by a different seed curve and
height function.

\medskip

\begin{Lem}\label{concat}
Let $\tilde{S}$ be a complete $C^2$ embedded H-minimal graph and
suppose $U, U'$ are neighborhoods of non-characteristic points $g,
g'$ respectively, and that $U$ is parameterized by
\[ \left (\gamma_1(s)+r
  \gamma_2'(s),\gamma_2(s)-r\gamma_1'(s),h_0(s)-\frac{r}{2}<\gamma(s),\gamma'(s)> \right) \ =\ \mathcal{L}_{\gamma(s)}(r)
  \] for $(s,r) \in O$, and $U'$ is parameterized by
\[\left (\rho_1(s)+r
  \rho_2'(s),\rho_2(s)-r\rho_1'(s),k_0(s)-\frac{r}{2}<\rho(s),
  \rho'(s)> \right)\ =\ \mathcal{L}_{\rho(s)}(r)
\]
for $(s,r) \in O'$. Further, suppose that $U \cap U' \neq
\varnothing$. Let $g_o \in U \cap U'$ be a non-characteristic
point and $g_o=\mathcal{L}_{\gamma(s)}(r) =
\mathcal{L}_{\rho(s')}(r')$ for some $(s,r)$ and $(s',r')$, then
$\mathcal{L}_{\gamma(s)}\cup \ \mathcal{L}_{\rho(s')}$ forms a
single straight line in $\mathbb H^1$.
\end{Lem}

\pf  Consider any non-characteristic point $g_o \in U \cap U'$.
There exist $(s_0,r_0),(s_0',r_0')$ so that
\[
\mathcal{L}_{\gamma(s_0)}(r_0)\ =\ \mathcal{L}_{\rho(s_0')}(r_0')\ .
\]

As $\tilde{S}$ is $C^2$, the horizontal Gauss map is $C^1$ away
from characteristic points.  As $\tilde{S}$ is embedded, we have
that the union of intersecting sub-surfaces of $\tilde{S}$ cannot
be transverse and must combine to form a $C^2$ sub-surface. Thus,
\eqref{wedge} implies that
\[ \lim_{r \ra r_0} \nuX(\mathcal{L}_{\gamma(s_0)}(r)) = \lim_{r \ra
  r_0'}\nuX(\mathcal{L}_{\rho(s_0')}(r))\ .\]

Since $\tilde{\mathcal{L}}'_{\gamma(s)}(r) = \tilp(\gamma(s))$ and
$\tilde{\mathcal{L}}'_{\rho(s)}(r) = \tilp(\rho(s))$ on their
respective domains of definition, we conclude that
$\mathcal{L}_{\gamma(s)}\cup \ \mathcal{L}_{\rho(s')}$ must form a
single straight line in $\mathbb H^1$.

 $\qed$

The proof of Lemma \ref{concat} shows that $U \cup U'$ is a
portion of $H$-minimal surface, defined as a graph over the
$xy$-plane, which is foliated by straight line segments.

\medskip

\subsection{Characteristic boundary points}\label{ss2}

Next, consider the case where the boundary point of the
neighborhood $U$ is a characteristic point.  Since, at a
characteristic point, the normal vector has a $\partial_t$
component, by the implicit function theorem there exists a
neighborhood of the characteristic point, $S$, that can be written
as a graph over the $xy$-plane. Using Theorem \ref{rep} we can
write $S = \cup S_i \cup \Sigma$. Here, each $S_i$ is a
neighborhood of a non-characteristic point given by Theorem
\ref{rep} parameterized using seed curve $\gamma^i$ and height
function $h_0^i$, while $\Sigma$ is the characteristic locus of
$S$.

Recalling the work of the previous section (in particular, Theorem
\ref{T:optregchar}), we see that there are two distinct cases,
first when $\Sigma$ is an isolated point and second when $\Sigma$
is a smooth curve. We first consider the case when $\Sigma$
consists of an isolated characteristic point.

\medskip

\begin{Lem}\label{L:isochar}  Suppose $S$, a complete embedded $C^2$ H-minimal
  surface, contains an isolated
  characteristic point $g_o$, then $S$ is a characteristic plane and is
  determined by a single seed curve which is a circle.
\end{Lem}

\pf Let $S_i$ be a portion of $S$ with $g_o \in \overline{S}_i$.
As the point is characteristic, the tangent plane at $g_o$ is the
horizontal plane passing through $g_o$, and any rule meeting this
point must therefore be completely contained in the plane. The
characteristic point, as a boundary point of one of the $S_i$,
arises either in case ($B_1$) or ($B_2$) described at the
beginning of the section.  If the point arises in case ($B_1$),
then there exists an $s'$ so that $g_o \in
\overline{\mathcal{L}}_{\gamma^i(s')}$.  If the point arises as in
case ($B_2$), then there exists a straight line containing $g_o$,
denoted again by $\mathcal{L}_{\gamma^i(s')}$, where $s'$ is in
the closure of the domain of $\gamma^i$.  This line is the limit
of lines $\mathcal{L}_{\gamma^i(s)}$ as $s$ approaches $s'$.  We
next observe that in the case of an isolated characteristic point,
$g_o$, all rules in $S_i$ must intersect $g_o$. To see this we
argue by contradiction, supposing that there are $s_j \ra s'$ so
that $g_o \nin \overline{\mathcal{L}}_{\gamma^i(s_j)}$. Since
$g_o$ is a characteristic
  point we have that, by equation \eqref{cl0}, there exists an $r_0$
  so that
\begin{equation}
(h_0^i)'(s')-\frac{1}{2}<(\gamma^i)'(s'),(\gamma^i)^\perp(s')>
-r_0+\frac{r_0^2}{2}\kappa^i(s')\ =\ 0\ .
\end{equation}

As $\gamma^i \in C^3$ and $h_0^i\in C^2$, by the implicit function
theorem we have that, for $\tilde{s}$ sufficiently close to $s'$,
there exists a $C^1$ function $r(\tilde{s})$ so that
\begin{equation}
(h_0^i)'(\tilde{s})-\frac{1}{2}<(\gamma^i)'(\tilde{s}),(\gamma^i)^\perp(\tilde{s})>
-r(\tilde{s})+\frac{r(\tilde{s})^2}{2}\kappa^i(\tilde{s})\ =\ 0\ .
\end{equation}

In other words, $\mathcal{L}_{\gamma^i(\tilde{s})}(r(\tilde{s}))$
is a characteristic point.  Applying this to the $s_j$ for $j$
sufficiently large, we see that either every
$\mathcal{L}_{\gamma^i(s_j)}$ contains $g_o$, or there exists a
sequence of characteristic points, $g_j \in
\mathcal{L}_{\gamma^i(s_j)}$ with $g_j \ra g_o$.  This is a
contradiction of the assumption that $g_o$ is an isolated
characteristic point.  Thus, all rules in $S_i$ contain $g_o$ and
we conclude that $S_i$ is a portion of a plane.  Now, by composing
with a left translation (which preserve $H$-minimality), we may
assume that the characteristic point is at the origin. But then,
every $S_i$ that has the origin in its closure is a subset of the
characteristic plane given by $t=0$, which is described in
Examples \ref{E:charplane} and \ref{E:comp}.  Let $S^0$ be the
union of all the $S_i$ that have the origin in their closure. By
the preceding discussion, $S^0$ is a subset of the plane and its
closure contains an open neighborhood of the origin.  From Example
\ref{E:charplane}, we have that if we choose $\gamma(0)=z=(x,y)$,
we obtain that
\[\gamma(s)=\left(x \cos\left ( \frac{s}{|z|}\right ) - y \sin \left (
  \frac{s}{|z|}\right ),y \cos\left ( \frac{s}{|z|}\right ) +x \sin \left (
  \frac{s}{|z|}\right ) \right )\]

Clearly, by picking $z$ sufficiently close to $0$, we may assume
that
  $\gamma \subset S^0$ and hence $\partial S^0 \cap \gamma =
  \varnothing$.  We claim that $\partial S^0$ contains only
  the isolated characteristic point in question.  To prove this, we
  suppose that $\partial S^0$ contains at least one other point, $g_1$.
  Since $\gamma \cap \partial S^0 = \varnothing$, there must exist a
  rule, $\mathcal{L}$, so that $g_1 \in \overline{\mathcal{L}}$.  By
  construction of the $S_i$, we know that: 1)
  either $g_1$ is another characteristic point; or, 2) the tangent
  plane to $S$ coincides with a vertical plane at $g_1$. In Example \ref{E:charplane} we computed the characteristic locus of a characteristic plane
  and we have seen that along any rule emanating from the isolated
  characteristic point at the origin, there are no other
  characteristic points. Thus, $g_1$ cannot be another characteristic
  point, and the possibility 1) is ruled out. Now, for $t=0$, the non-unit Riemannian normal to the point
  $(x,y,0)$ is $-\frac{y}{2}\; X_1 +\frac{x}{2} \; X_2 +T$ and hence, at any point, the tangent
  plane is not a vertical plane.  Hence, the possibility 2) for $g_1$ is ruled out as well and we conclude that $\partial S_i$ contains only the
  isolated characteristic point.  By construction, we
  have that $\tilde{S}=\overline{S^0}$ and therefore $\tilde{S}$ is a characteristic
  plane and is determined by the single seed curve given above.

$\qed$

\medskip

We next address the case when $\Sigma$ is a smooth curve.

\medskip

 \begin{Lem}\label{L:extchar}  Let $\tilde{S}$ be a complete embedded
   $C^2$ H-minimal surface.  Let $S$ be a connected component of
   $\tilde{S}$ which can be written as a graph over the $xy$-plane and
   is decomposed as $S = \cup S_i \cup \Sigma$ via Theorem \ref{rep}.
Suppose $g_o \in \Sigma$ is a non-isolated characteristic point of
type ($B_1$) and $S_j$ and
  $S_k$ are pieces of the decomposition of $S$ which have $g_o$ in their boundary. Let
  $\mathcal{L}_j$ be a rule in $S_j$ which terminates at $g_o$, then
  $\mathcal{L}'_j$ can be extended over $g_o$ and remains in $\tilde{S}$.
\end{Lem}

\pf Once again, the conclusion follows from the assumption that
the surface is $C^2$, and hence the Riemannian normal vector is
$C^1$. By the implicit function theorem, we can write a
neighborhood of the characteristic point as a graph over a portion
of the $xy$-plane. In other words, there exist $\Omega_0 \subset
\R^2$, and  $f_0: \Omega_0 \ra \R$, so that the neighborhood is
given by $(x,y,f_0(x,y))$ where $(x,y) \in \Omega_0$.  By the
assumption, the subsets $S_k$ and $S_j$ intersect this
neighborhood. Now, a non-unit Euclidean normal over this
neighborhood is given by
\[
 N\ =\ \left(-\frac{\partial f_0}{\partial x}(x,y), -\frac{\partial f_0}{\partial y}(x,y) ,1\right)\ .
 \]

The latter formula gives
\[
 N_x\ =\ \left(-\frac{\partial^2 f_0}{\partial
  x^2}(x,y),-\frac{\partial^2 f_0}{\partial y \partial x}(x,y),0\right)\ ,
 \]
\[
N_y\ =\  \left(-\frac{\partial^2 f_0}{\partial
  x \partial y}(x,y),-\frac{\partial^2 f_0}{\partial y^2}(x,y),0\right)\ .
 \]

First, suppose that $S_k$ contains a rule, $\mathcal{L}_k$, so
that $g_o \in \overline{\mathcal{L}}_k$.  We consider a Riemannian
normal vector along the two rules, $\mathcal{L}_j, \mathcal{L}_k$.
Denoting the horizontal Gauss map of $S_i$ by $\nuX^i=\pb^i
X_1+\qb^i X_2$, $i = j , k$, in $(s,r)$ coordinates, we obtain
from \eqref{run}
\[
N_i\ =\ W^i(s,r)\ \nuX^i\ +\ T\ .
\]

By \eqref{enorm}, we have that re-calculating the normal along the
rules with respect to the Euclidean metric yields:
\begin{equation}\label{Ni}
N_i(s,r)\ =\ \left(W^i(s,r) \pb^i + \frac{\gamma_2^i(s)- r
(\gamma_1^i)'(s)}{2},
W^i(s,r)\qb^i-\frac{\gamma_1^i(s)+r(\gamma_2^i)'(s)}{2}\ ,\
1\right)\ .
\end{equation}

Note that the patch of surface is a planar one if and only if
$N_i(s,r)$ is a constant vector, and $\partial N_i/\partial
r(s,r)=(0,0,0)$. Since we assume that the surface is not planar,
we may suppose that $\Sigma$ does not consist of a single isolated
point - Example \ref{E:charplane} shows that all planes that are
graphs over the $xy$-plane have a single characteristic point.
Recalling that $\pb^i=(\gamma_1^i)'(s)$ and
$\qb^i=(\gamma_2^i)'(s)$, we have
\[
N_i'(s,r)\ =\ \left(W^i_r(s,r)\pb^i-\frac{\pb^i}{2},W^i_r(s,r)
\qb^i-\frac{\qb^i}{2},0\right)\ .
\]

By formula \eqref{deqW},
\[
W^i_r(s,r)\ =\ 1\ +\ \frac{\kappa^i(s)}{1-r\kappa^i(s)} W^i(s,r)\
.
\]

Since at characteristic points $W^i(s,r)=0$, we conclude that at
such points $W^i_r = 1$. Next, we restrict our attention to
$\mathcal{L}_j$ and $\mathcal{L}_k$.  We know the following facts:

\begin{enumerate}
\item there are two $r$ values, $r_j,r_k$ so that the appropriate one
  sided limits along the rules approach $g_o$ as $r$ approaches these
  values.  As we do not, a priori, know if the limits are left or
  right handed limits, we will denote them as follows:
\[ \lim_{r \ra r_j^\sharp} \mathcal{L}_j(r) = g_o\ ,\quad\quad \text{and}\quad\quad    \lim_{r \ra r_k^\sharp} \mathcal{L}_k(r) = g_o\ ;
\]
\item $(\pb^j,\qb^j)$ and $(\pb^k,\qb^k)$ are constant along
  $\tilde{\mathcal{L}}_j$ and $\tilde{\mathcal{L}}_k$ respectively.
\item By Remark \ref{R:rems}, we know that there are at most two
  characteristic points along any given rule.  Thus, for $r$
  sufficiently close to $r_i$, there are no other characteristic
  points along $\mathcal{L}_i$ for $i=j,k$
   .
\end{enumerate}
Using these facts, we have
\begin{equation*}
\lim_{r \ra r_j^\sharp}N_j'(r) =
\left(\frac{\pb^j}{2},\frac{\qb^j}{2},0\right)\ ,\quad\quad\quad
\lim_{r \ra r_k^\sharp}N_k'(r) =
\left(\frac{\pb^k}{2},\frac{\qb^k}{2},0\right)\ .
\end{equation*}

Now, both $N_j'$ and $N_k'$ are (one-sided) directional
derivatives of $N$ in the directions of the derivatives of $\mathcal{L}_j$ and $\mathcal{L}_k$ respectively. Precisely, we have for $i=j,k$,

\begin{equation}
N'_i= \begin{pmatrix}
  -(f_0)_{xx} & -(f_0)_{yx} &0\\ -(f_0)_{xy} & -(f_0)_{yy}^0\\ 0&0&0 \end{pmatrix} \begin{pmatrix} \qb^i \\ -\pb^i\\0 \end{pmatrix}
\end{equation}

Re-writing these equations, we have
\begin{equation}\label{c2eqs}
\begin{split}
\frac{\pb^j}{2} &= -\qb^j (f_0)_{xx}+\pb^j (f_0)_{yx} \\
\frac{\pb^k}{2} &= -\qb^k (f_0)_{xx}+\pb^k (f_0)_{yx} \\
\frac{\qb^j}{2} &= -\qb^j (f_0)_{xy}+\pb^j (f_0)_{yy} \\
\frac{\qb^k}{2} &= -\qb^k (f_0)_{xy}+\pb^k (f_0)_{yy}\ .
\end{split}
\end{equation}

As we may compose with left-translations and rotations without
affecting minimality, we may assume without a loss of generality
that $\pb_j,\pb_k,\qb_j,\qb_k$ are all non-zero.  Indeed, we may
left-translate the surface so that $g_o$ is at the origin.  Then,
rotation of the $H$-minimal surface has the effect of rotating the
horizontal Gauss map.  If $(\pb_j,\qb_j)$ and/or $(\pb_k,\qb_k)$
had a zero component, composition with a rotation would place the
vectors in general position.  Solving the first and second
equations for $(f_0)_{yx}$, and the third and fourth equations for
$(f_0)_{xy}$, we have
\begin{equation*}
\begin{split}
\frac{1}{2} + \frac{\qb^j}{\pb^j} (f_0)_{xx}&=(f_0)_{yx}\ ,\\
\frac{1}{2} + \frac{\qb^k}{\pb^k} (f_0)_{xx}&=(f_0)_{yx}\ ,\\
-\frac{1}{2} + \frac{\pb^j}{\qb^j} (f_0)_{yy} &=(f_0)_{xy}\ ,\\
-\frac{1}{2} + \frac{\pb^k}{\qb^k} (f_0)_{yy} &=(f_0)_{xy}\ .\\
\end{split}
\end{equation*}

The first pair tells us that
\begin{equation}\label{eqs3}
\frac{\qb^j}{\pb^j}=\frac{\qb^k}{\pb^k}\ ,
\end{equation}
or $(f_0)_{xx}=0$.  The second pair says either (\ref{eqs3}) is
true, or $(f_0)_{yy}=0$.  If equation (\ref{eqs3}) is true, then
the two rules, $\mathcal{L}_j,\mathcal{L}_k$, point in the same
direction and, since they both contain $g_o$, they are therefore
parts of the same line. The fact that they are parameterized in
opposite directions implies $\nuX^j=-\nuX^k$.  If, on the other
hand, $(f_0)_{xx}=0=(f_0)_{yy}$, the equations are inconsistent,
implying that $(f_0)_{xy}$ is both $\frac{1}{2}$ and
$-\frac{1}{2}$. Now, consider the case where $g_o \in \partial
S_k$ arises from case ($B_2$).  Then, $g_o= \lim_{s \ra s'}
(F(s,r),h(s,r))$ (where $r$ is fixed).  We can construct a line,
$\mathcal{L}_k= \lim_{s \ra s'} \mathcal{L}_{\gamma^k(s)}$.  As
$\tilde{S}$ is metrically complete , $\mathcal{L}_k$ is a line that lies in the
boundary of $S_k$. Next, consider the normal along $\mathcal{L}_k$.  Repeating the analyis at the beginning of this proof along the line $\mathcal{L}_k$, i.e. equations \eqref{c2eqs} and \eqref{eqs3},  yields that $\mathcal{L}_j \cup \mathcal{L}_k$ forms a
single line and hence $\mathcal{L}_j$ can be extended over $g_o$.
This completes the proof.

$\qed$

\subsection{Points where $\tilde{S}$ ceases to be a graph}\label{SS:vert}

Last, we consider the case where a boundary point of $U$ is a
point where $\tilde{S}$ ceases to be a graph over the $xy$-plane.
At such a point, the Euclidean normal vector must have no
$\partial_t$ component.  Recalling the discussion at the beginning
of this section we have two possibilities:
\begin{enumerate}
\item[($C_1$)] The boundary of $U$ contains the image of the singular locus, $r=
\frac{1}{\kappa(s)}$\ ;
\item[($C_2$)] The seed curve, $\gamma$, hits the edge of $S$ and $W_0(s)=-(h_0)'(s)$
  tends to $\pm \infty$ as we approach these points.
\end{enumerate}

To deal with the first case, we prove a lemma showing that  we can extend the rules over the singular locus.  Again,
this follows from the assumption that the surface is $C^2$.

\medskip

\begin{Lem}\label{L:singloc} Let $g_o \in \partial U$ be a point which is
  non-characteristic and of type ($C_1$) (i.e., in the image
  of the singular locus, $r = \frac{1}{\kappa(s)}$).  If $\mathcal{L}_{\gamma(s)}$
  is a rule in $U$ that contains $g_o$ in its closure, then $\mathcal{L}_{\gamma(s)}$ can be
  extended over $g_o$ so as to remain in $\tilde{S}$.
\end{Lem}

\pf  Consider a neighborhood $S_0 \subset \tilde{S}$ of $g_o$.  At the outset of this section, we fixed $U$ as a portion of $\tilde{S}$ that could be parameterized in terms of a single seed curve and height function as in theorem 4.6.  In particular, this implies that $U$ can be written as a graph over the xy-plane and thus $S_0 \cap U$ can also be written as a graph over the $xy$-plane. Then, either
$S_0 \minus U$ has a component, $U'$, with $g_o$ as a boundary
point, that can be written as a graph over the $xy$-plane, or $S_0
\minus U$ has no such component. In the second case, by Lemma
\ref{trivial}, $S_0 \minus U$ is contained entirely in a vertical
plane with normal $\nuX(g_o)$. Consider a nearby rule,
$\mathcal{L}_{\gamma(s+\epsilon)}$, where $\epsilon$ is
sufficiently small.  If the limit points of
$\mathcal{L}_{\gamma(s+\epsilon)}$ do not contain the image of
$r=\frac{1}{\kappa(s+\epsilon)}$, then by either Lemma
\ref{concat} or Lemmas \ref{L:isochar} and \ref{L:extchar}, we can
extend $\mathcal{L}_{\gamma(s+\epsilon)}$ until they do.  Let
$g_1$ be such a limit point. Repeating the same analysis, we have
that $\mathcal{L}_{\gamma(s+\epsilon)}$ must lie inside a vertical
plane with normal given by $\nuX(g_1)$.  By the hypothesis that
$g_o$ occurs at $r = \frac{1}{\kappa(s)}$, we know that $\kappa(s)
\neq 0$.  Thus, by the continuity of $\kappa$, for $\epsilon$
sufficiently small, $\kappa(s+\epsilon) \neq 0$, and thus
$\gamma'$ is non-constant in a neighborhood of $s$.  So, again for
$\epsilon$ sufficiently small, we may assume that $\nuX(g_o) \neq
\nuX(g_1)$, and hence these portions of vertical planes must
intersect.  But, for $\tilde{S}$ to be a $C^2$ embedded surface,
we must have that $\nuX(g_1) = \nuX(g_o)$.  This is a
contradiction.

Having dealt with the case where $S_0 \minus U$ has no graphical
component, we may assume that there is another neighborhood,
$U^*$, parameterized using seed curve $\rho$ and height function
$f_0$, which is a graph over the $xy$-plane that has $g_o$ as a
boundary point. Again, there are two possibilities, corresponding
to ($C_1$) and ($C_2$) above.  In case ($C_1$), where $g_o$ is in
the image of the singular locus of $U^*$, we consider two rules
$\mathcal{L} \subset U, \mathcal{L}^* \subset U^*$ which hit
$g_o$.
 We again compute the normals along these rules.  Recall
 that the unit Riemannian normal along $\mathcal{L}$ is given by
\[ \frac{W(s,r)}{\sqrt{1+W(s,r)^2}}\ \nuX\ +\ \frac{1}{\sqrt{1+W(s,r)^2}}\ T\ ,
\]
and that as $r \ra \frac{1}{\kappa(s)}$, $W(s,r) \ra \pm \infty$.
Similarly, the unit Riemannian normal along $\mathcal{L}^*$ is
given by
\[ \frac{{W^*}(s,r)}{\sqrt{1 + {W^*}(s,r)^2}}\ \nuX^*\ +\ \frac{1}{\sqrt{1+ {W^*}(s,r)^2}}\ T\ ,
\]
with ${W^*}(s,r) \ra \pm \infty$ as $r \ra
\frac{1}{{\kappa^*}(s)}$. As $r \ra \frac{1}{\kappa(s)}$ we have
that the normal along $\mathcal{L}$ tends to $\nuX$ and that as $r
\ra \frac{1}{{\kappa^*}(s)}$ we have that the normal along
${\mathcal{L}^*}$ tends to $\boldsymbol{\nu}^*_H$.  Since the
surface is $C^2$, these limits must coincide with each other and
the Riemannian normal at the point.  Hence, we must have
$\nuX=\boldsymbol{\nu}^*_H$, and thus we may extend
$\mathcal{L}_{\gamma(s)}$ over the image of the singular locus.

If we are in case ($C_2$), the lifted seed curve of $U^*$,
$(\rho(s),f_0(s))$, contains $g_o$ as a limit point.  Let $\nu_1$
be the limit of the horizontal Gauss map along
$\mathcal{L}_{\gamma(s)}$ as $r \ra r_i$. In other words:
\[\nu_1=\lim_{r \ra r_i} (\mathcal{L}'_{\gamma(s)}(r))^\perp= \gamma'(s)\]

As noted before, $f_0'(s)$ tends to $\pm \infty$ in this limit, so
that the unit Riemannian normal tends to $\nu_2=\rho'(s)$ as $s
\ra s_j$.  Again, as the normal is $C^1$, we have that
$\nu_1=\nu_2$.  Now, consider a sequence of rules in $U^*$,
${\mathcal{L}^*}_{\rho(s_k)}$ where $s_k \ra s_j$ as $k \ra
\infty$.  Then, since $\tilde{S}$ is complete, $\lim_{k \ra
\infty} {\mathcal{L}^*}_{\rho(s_k)}$ is a segment of a line
containing $g_o$ and perpendicular to $\nu_2$. As $\nu_1=\nu_2$,
this line segment extends $\mathcal{L}_{\gamma(s)}$.  $\qed$

\medskip

We next turn to points of type ($C_2$).  Consider $U$, an open set
which can be written as a graph over the $xy$-plane, with $g_o$ in
its boundary.  Consider a neighborhood $S_0$ of $g_o$ in the
$H$-minimal surface $S$.  Then, $S_0 \minus U$ may have a
component, $U'$ which can be written as a graph over the
$xy$-plane, or it may not have such a component.  We address these
two possibilities in the next lemma.

\medskip

\begin{Lem}\label{L:C2} Let $g_o \in \partial U$ be a point of type
  ($C_2$).  Then, one of three
  possibilities occur:
\begin{enumerate}
\item There exists another neighborhood,
  $U' \subset S$, containing $g_o$ as a boundary point, which can be
  written as a graph over the $xy$-plane and that $g_o$ is a point of
  type ($C_1$) for $U'$.  Then, the rule $\mathcal{L} \subset U'$
  whose closure contains $g_o$ can be extended over $g_o$ and remain
  in $\tilde{S}$.
\item There exists another neighborhood,
  $U' \subset S$, containing $g_o$ as a boundary point, which can be
  written as a graph over the $xy$-plane and that $g_o$ is a point of
  type ($C_2$) for $U'$.  Then, there exists a
  vertical plane $V$ and a rule, $\mathcal{L} \subset V \cap S$, both
  containing $g_o$, so that $S$ is tangent to $V$ along $\mathcal{L}$.
\item There is no other neighborhood in $S$ which contains $g_o$ as a
  boundary point can be written as a graph over the $xy$-plane.  Then,
  there exists a vertical plane, $V_0$, and a rule $\mathcal{L}
  \subset S$ with  $g_o \in \mathcal{L} \subset V_0$ so that $S$ is
  tangent to $V_0$ along $\mathcal{L}$.
\end{enumerate}
\end{Lem}
\pf We may assume that $U$ is parameterized by a
seed curve $\gamma$ and height function $h_0$:
\[ U = \{F(s,r),h(s,r)| s \in (a,b), r \in (c,d)\}\]
and that  there exist parameter values $s_0$, $r_0 \neq \lim_{s
\ra s_0} \frac{1}{\kappa(s)}$ such that
\[ g_o = \lim_{s \ra s_0} (F(s,r_0),h(s,r_0))\ .\]

Let \[\mathcal{L}=\lim_{s \ra s_0} \mathcal{L}_{\gamma(s)}\]

We next examine the possible behavior of $S$ near $g_o$.  First,
suppose there exists another neighborhood,
   $U' \subset S$, which contains $g_o$ as a boundary point, and can be
   written as a graph over the $xy$-plane.  Given the set $U'$, if the point $g_o$ is of type ($C_1$), then by
 reversing the roles of $U$ and $U'$, we may apply Lemma
 \ref{L:singloc}, and conclude that the rule in $U'$ may be extended
 over the point $g_o$.  We may now assume that $g_o$ is a point of type
   ($C_2$) for $U'$, or that no such $U'$ exists.  Condition ($C_2$) assumes that $W(s,r_0) \ra \pm \infty$ as we approach
this point.  Examining \eqref{E:t1}, and recalling that
$W_0(s)=-h_0'(s)-\frac{1}{2} <\gamma(s),(\gamma')^\perp(s)>$, this
implies that, for fixed $r=r_0$, $h_0'(s) \ra \pm \infty$ as $s
\ra s_0$.  Recall that the unit Riemannian normal to the surface
is given by
\[ \frac{W}{\sqrt{1+(W)^2}} \; \nuX + \frac{1}{\sqrt{1+(W)^2}}
\; T\ ,\]
 and that
\[W(s,r)= \frac{W_0(s) +r - \frac{r^2}{2} \kappa(s)}{1-r\kappa(s)}\ ,\]
we have that as $s$ approaches $s_0$ for any fixed $r \neq \lim_{s \ra
  s_0} \frac{1}{\kappa(s)}$, the Riemannian normal tends to the
horizontal Gauss map, $\nuX$.  Since $\nuX$ is constant in $r$, we see
that the Riemannian normal tends to a fixed vector $\lim_{s \ra s_0}
\nuX$ as $s \ra s_0$ along $\mathcal{L}$.  Let $V_0$ be the vertical
  plane passing through $g_o$ with normal given by $\nuX(g_o)$.  Then,
  by construction, $\mathcal{L} \subset V_0$ and $S$ is tangent to
  $V_0$ along $\mathcal{L}$.

  $\qed$

We now provide an application of these gluing operations, showing
that a neighborhood of a non-characteristic point on a complete
embedded $C^2$ $H$-minimal graph parameterized by a single seed
curve and height function can be extended so that the same seed
curve and height function determine a substantial portion of the
surface.

\medskip

\begin{Lem}[Extension lemma]\label{extension}
Let $\tilde{S}$ be a complete $C^2$ embedded $H$-minimal graph and
suppose $S$ is a neighborhood of a point $g \in \tilde{S}$ which
is parameterized by
\begin{equation}\label{extlem1}
\left (\gamma_1(s)+r
\gamma_2'(s),\gamma_2(s)-r\gamma_1'(s),h_0(s)-\frac{r}{2}<\gamma(s),
\gamma'(s)> \right)
\end{equation}
for $(s,r) \in O = \{ (s,r) | s \in(s_0,s_1), r \in (-r_0(s),r_1(s))
\}$.  Then, $S$ and the representation of $S$ given by \eqref{extlem1}
can be extended to a domain $\tilde{O} = \{ (s,r) | s \in (s_0,s_1), r
\in (-\infty,\infty)\}$ and $S'$, the surface given by \eqref{extlem1}
with $(s,r) \in \tilde{O}$ is a subset of $\tilde{S}$.
\end{Lem}

\pf Consider, for a fixed $s \in (s_0,s_1)$, the line
\[ \mathcal{L}_{\gamma(s)}(r) = \left (F(s,r),h_0(s) - \frac{r}{2}
<\gamma(s), \gamma'(s)> \right )\ .
\] for $r \in
(-r_0(s),r_1(s))$. Consider first the case where $\kappa(s) >0$
and $r_1(s) < \frac{1}{\kappa(s)}$.  Then
  \eqref{enorm} shows that as $r \ra r_1(s)$, the normal to the
  surface has a $\frac{\partial}{\partial t}$ component, and hence
  $\tilde{S}$ is a graph over the $xy$-plane in a neighborhood of the
  point $\lim_{r \ra r_1(s)^-} \mathcal{L}_{\gamma(s)}(r)$.  By
  Theorem \ref{rep}, the representation of
  $\tilde{S}$ in a neighborhood of this limit point is determined by
  a seed curve and height function.  If the limiting point is
  non-characteristic, then Lemma \ref{concat} allows us to extend
  the line $\mathcal{L}_{\gamma(s)}(r)$ over a neighborhood of this
  limit point.  Suppose instead that the limiting point is characteristic.  Then, as the
  normal at this point has a nonzero $\frac{\partial}{\partial t}$ component,
  a neighborhood of the characteristic point can be written as a graph
  over the $xy$-plane.  Thus, in addition to $S$, there is at least one
  other patch of $\tilde{S}$, $S_0$, parameterized as in Theorem
  \ref{rep}, so that the limit point is in $\overline{S}_0$.  Lemma
  \ref{L:extchar} shows us that the conclusion of Lemma \ref{concat} is
  still true in this case, allowing us to extend the line in $S$ over the
  characteristic point in $S_0$. Continuing in this manner, we see that we can extend
  the line $\mathcal{L}_{\gamma(s)}(r)$ at least up until $r =
  \frac{1}{\kappa(s)}$, where the $\frac{\partial}{\partial t}$
    component of $\vec{n}$ vanishes.  At such a point, if one of the
    first two components of $\vec{n}$ is nonzero, the surface ceases
    to be a graph over the $xy$-plane.  Lastly, Lemma \ref{L:singloc} shows
    that if we extend the line $\mathcal{L}_{\gamma(s)}(r)$ beyond $r =
    \frac{1}{\kappa(s)}$, the line remains in $\tilde{S}$.  Similarly, there is no
    obstruction to extending the line in the other direction using Lemma \ref{concat}.  Thus, in this case, we may replace the
    domain of $\mathcal{L}_{\gamma(s)}$ with the interval $\left (-\infty,
    \infty \right )$.

    $\qed$

\medskip

The following is a technical lemma needed for Theorem
\ref{T:graphs}.

\medskip

\begin{Lem}\label{L:intersect}
Let $\tilde{S}$ be a complete embedded $C^2$ $H$-minimal surface
and suppose $S_1,S_2,S_3,S_4$ are open subsurfaces of $\tilde{S}$
which are graphs over regions of the $xy$-plane such that $S_1
\cup S_2 \cup S_3 \cup S_4$ forms a connected set, and for which
one has for each $j$, $\cup_{i \neq j} S_i \neq \cup_{i=1}^4 S_i$.
Further, suppose that each $S_j$ is parameterized via a seed curve
$\gamma^j$ and height function $h_0^j$ using equations
\eqref{para} and \eqref{para2} with $s \in (s_0^j,s_1^j)$ and $r
\in \R$.  Then, $S_1$ can intersect at most two of
$\{S_2,S_3,S_4\}$.
\end{Lem}

\pf  We argue by contradiction, and assume that there exist points
$g_2 \in S_2$, $g_3 \in S_3$ and $g_4 \in S_4$ so that $g_i \in
S_1$ as well.  Thus, there exist $s_i,r_i$ so that
\[\mathcal{L}_{\gamma^1(s_i)}(r_i)=g_i\]
for $i=\{2,3,4\}$. By Lemma \ref{extension},
\[\mathcal{L}_{\gamma^1(s_i)}(r) \subset S_1 \cap S_i\] for $r \in
I^1(s_i)$.  By reordering the indices, we may assume that $s_2\le
s_3\le s_4$.  By construction, $S_3$ contains either
$\{\gamma^1(s)|s \in [s_3, s_1^1)\}$ or $\{\gamma^1(s)|s \in
(s_0^1,s_3]\}$.  In the first case, $S_3$ contains $g_4$, and
depending on the length of the intervals $(s_0^3,s_1^3)$ and
$(s_0^4,s_1^4)$, either $S_3 \subset S_1 \cup S_4$ or $S_4 \subset
S_1 \cup S_3$, violating the assumption that for each $j$,
$\cup_{i \neq j} S_i \neq \cup_{i=1}^4 S_i$. The other cases give
similar violations of this assumption.  Hence, $S_1$ can intersect
at most two of the other $S_j$.  $\qed$

\begin{Thm}\label{T:graphs}  Let $\tilde{S}$ be a complete embedded
  $C^2$ H-minimal surface.  For a fixed $g_0 \in \tilde{S}$, let $S
  \subset \tilde{S}$ be the largest connected component of
  $\tilde{S}$ containing $g_0$ that can be written as a graph over the
  $xy$-plane.  Then $S$ can be parameterized by \eqref{para} and \eqref{para2} using a single
  seed curve and height function with $s \in (s_0,s_1)$ and
\begin{equation*}
r \in I(s) =
\begin{cases}
\left (-\infty ,\frac{1}{\kappa(s)}\right) \;\;\; \kappa(s)>0\\
\left (\frac{1}{\kappa(s)},\infty\right) \;\;\; \kappa(s)<0
\end{cases}
\end{equation*}

Moreover, if we extend $S$ to $S'$ by allowing for $ r \in\R$,
then $S' \subset \tilde{S}$.
\end{Thm}

\pf  Covering $S$ with open neighborhoods generated by Theorem
\ref{rep} and extending them via the extension Lemma
\ref{extension}, we can write $S = \cup S_j$ where each $S_j$ is
parameterized by a seed curve $\gamma^j(s)$ and height function
$h_0^j(s)$ with $s\in(s_0^j,s_1^j)$.  We fix $S_0$ to be the
neighborhood containing $g_0$ and assume, by picking a new basepoint for the seed curve if necessary, that $(\gamma^0(0),h^0_0(0))=g_0$.
 By Lemma
\ref{L:intersect}, we may assume that we can order the $S_j$ so
that \[S= \cup_{j=-\infty}^\infty S_j\]
 with the property that $S_j \cap S_{j+1} \neq \varnothing$ but $S_i
 \cap S_j = \varnothing$ if $|i-j| > 1$.

Since all of the $S_i$ are open, there exists $\epsilon >0$ so
that for $s \in (s_1^0-\epsilon, s_1^0)$, $\gamma^0(s) \in S_1$.
Thus, for $s \in (s_1^0-\epsilon, s_1^0)$, there exists $\sigma
\in (s_0^1,s_1^1)$, $r \in I(\sigma)$ so that $F^1(\sigma, r)=
\gamma^0(s)$.  By remark \ref{R:intcurve}, we know that for fixed
$r$, the image of $F(\sigma, r)$ coincides with the image of an
integral curve of $\nuX$ away from characteristic points.  As
$\gamma^0$ is such an integral curve, we can use $F^1(\sigma, r)$
to extend $\gamma^0$ as follows.  First, let $\gamma^{0,1}(s)$ be
a re-parametrization of the curve $F^1(\sigma, r)$ by arclength. A
priori, $\gamma^{0,j}(s)$ has domain that may not match up with
the domain of $\gamma^0$ on the overlap between $S_0$ and $S_1$ -
i.e. the point $\gamma^0(s)$ may be given by $\gamma^{0,1}(s')$
where $s \neq s'$.  But, using a simple re-parametrization of
$\gamma^j$, we may assume that the domain of $\gamma^{0,1}$ is
$(\tilde{s}_0^1,\tilde{s}_1^1)$ with $\tilde{s}_1^1 \in
(s_0^0,s_1^0)$ and for every $\sigma \in
(\tilde{s}_0^1,\tilde{s}_1^1) \cap (s_0^0,s_1^0)$,
$\gamma^0(\sigma)=\gamma^{0,1}(\sigma)$.  We note that
$\gamma^{0,1}$ may contain characteristic points but is still a
well defined curve that, away from the characteristic points is an
integral curve of $\nuX$.  We observe that if $\tilde{S} \in C^k$,
$\gamma^{0,1} \in C^{k+1}$ even over characteristic points.  To
see this, we note that
\[ \frac{\partial}{\partial s} F^1(s,r) = ((\gamma^1_1)'(s)+r
(\gamma^1_2)''(s), (\gamma_2^1)'(s)-r(\gamma_1^1)''(s))\] and the
length of this vector is $1-r \kappa^1(s)$.  Thus, the unit
tangent vector to this curve coincides, up to sign, with
$(\gamma^{0,1})'$ due to the re-parametrization and is given by
\begin{equation*}
\begin{split}
(\pm \gamma^{0,1})'(s) &=
\left ( \frac{(\gamma^1_1)'(s)+r (\gamma^1_2)''(s)}{1-r\kappa^1(s)},
  \frac{ (\gamma_2^1)'(s)-r(\gamma_1^1)''(s)}{1-r\kappa^1(s)} \right)\\
&= \left ( (\gamma_1^1)'(s)+ r
  \frac{(\gamma_2^1)''(s)+(\gamma_1^1)'(s)\kappa^1(s)}{1-r\kappa^1(s)},(\gamma_2^1)'(s)- r
  \frac{(\gamma_1^1)''(s)-(\gamma_2^1)'(s)\kappa^1(s)}{1-r\kappa^1(s)}
  \right )
\end{split}
\end{equation*}

Now, since
$\kappa^1(s)=(\gamma_1^1)''(s)(\gamma_2^1)'(s)-(\gamma_2^1)''(s)(\gamma_1^1)'(s)$
and $\frac{\partial}{\partial s} < (\gamma^1)'(s),(\gamma^1)'(s)> =0$,
we have
\begin{equation*}
\begin{split}
(\gamma_2^1)''(s)+(\gamma_1^1)'(s)\kappa^1(s) &=
(1-(\gamma_1^1)'(s)^2)(\gamma_2^1)'(s)(\gamma_2^1)''(s)+(\gamma_1^1)'(s)(\gamma_2^1)'(s)^2(\gamma_1^1)''(s)
\\
&= (\gamma_2^1)'(s)^2 ((\gamma_1^1)'(s)(\gamma_1^1)''(s)
+(\gamma_2^1)'(s)(\gamma_2^1)''(s))\\
&=0
\end{split}
\end{equation*}
Similarly, $(\gamma_1^1)''(s)-(\gamma_2^1)'(s)\kappa^1(s) =0$ and so
\[\pm (\gamma^{0,1})'(s)= (\gamma^1)'(s)\]
Since $\gamma^1 \in C^{k+1}$, so is $\gamma^{0,1}$.  We note that the
signs of these two vectors must match (i.e. $(\gamma^{0,1})'(s)=
(\gamma^1)'(s)$), otherwise $\nuX$ would not be constant along the
rules.

Second, we extend $\gamma^0$ as follows:
\begin{equation*}
\gamma(s) =
\begin{cases}
\gamma^0(s) \;\;\; s \in (s_0^0,s_1^0)\\
\gamma^{0,1}(s) \;\;\; s \in (\tilde{s}_0^1,\tilde{s}_1^1)
\end{cases}
\end{equation*}

Similarly, we can form a function
$h_0^{0,1}(s)=h_0^1(s)-\frac{r}{2}\gamma^1(s) \cdot (\gamma^1)'(s)$
and extend $h_0^0$ as well:
\begin{equation*}
h_0(s) =
\begin{cases}
h_0^0(s) \;\;\; s \in (s_0^0,s_1^0)\\
h_0^{0,1}(s) \;\;\; s \in (\tilde{s}_0^1,\tilde{s}_1^1)
\end{cases}
\end{equation*}
We note that the extended curve now parameterizes $S_0 \cup S_1$ using
equations \eqref{para} and \eqref{para2}.  Repeating this construction
for the endpoints $s_1^j$ and the analogous extension for the
endpoints $s_0^j$, we see that we may extend $\gamma$ and $h_0$
infinitely many times so that using these two functions in
\eqref{para} and \eqref{para2} parameterizes all of $S$.  $\qed$

\begin{Def}\label{D:extgraph}  If $S$ is a subset of a complete
  embedded $C^2$ H-minimal surface that is a graph over the $xy$-plane,
  the set $S'$, described in the previous theorem and obtained by
  extending all the rules foliating $S$, is called the {\bf extended
    graph} associated to $S$.
\end{Def}

\vskip 0.6in


\section{\textbf{A characterization of complete embedded $C^2$ $H$-minimal surfaces}}\label{S:characterization}

\vskip 0.2in

With the results of the last section, we can now extend our
description of complete embedded $C^2$ $H$-minimal surfaces using a
generalization of the notion of seed curves.  At this point, we can use Theorem \ref{T:graphs} to say the following:
given a complete embedded $C^2$ $H$-minimal surface $S$, let $S_0 = S_1 \cup
S_2 \cup ...$ be the union of subsets
of $S$, each of which can be written as extended graphs over the
$xy$-plane.  Thus, we can find seed curves $\gamma^i$ and height
functions $h_0^i$ so that
\[S_i = \left \{ F^i(s,r),h_0^i(s)-\frac{r}{2}\gamma^i \cdot
  (\gamma^i)'(s) | s \in (s_0^i,s_1^i), r \in \R \right \}\]
where $F^i(s,r) = (\gamma_1^i(s)+r
(\gamma_2^i)'(s),\gamma_2^i(s)-r(\gamma_1^i)'(s))$.  This gives us a
decomposition of $\tilde{S}$:
\[\tilde{S} = \cup S_i \bigcup V\]
where $V$ is the portion of $\tilde{S}$ where the tangent plane
becomes vertical and that does not arise from the image of the
singular locus for one of the $\gamma^i$.  Indeed, the boundary points
of $S_i$ must fall into category ($C_2$) in section \ref{S:glue}.  By
lemma \ref{L:C2}, the surface $S$ is tangent to a vetrical plane
at these points.  To finish our description
of $\tilde{S}$, we need to understand how these pieces can fit
together.



To allow a succinct summary and to make precise the way the various
pieces fit together,  we introduce the following definitions.

\medskip

\begin{Def}
If no portion of $S$ can be written as a graph over the $xy$-plane, we say that $S$ has \emph{trivial seed
  curve}.  Otherwise, we say that $S$ has \emph{non-trivial seed curves}.
\end{Def}

\medskip

\begin{Def}\label{D:genseed}
  If, in the representation $S=V\cup S_1 \cup S_2 \cup \dots$, the surface $S$ has non-trivial seed
  curves, then we call the collection
  $\{(\gamma_1^i(s),\gamma_2^i(s),h_0^i(s)), s \in (a_i,b_i)\}$, where
  $\lim_{s \ra b_i^-}\gamma^i(s)=\lim_{s \ra a_{i+1}^+} \gamma^{i+1}(s)$,
  a {\em generalized seed curve} for $S$.
\end{Def}

\medskip

Putting the results of the previous sections together, we obtain the following theorem.

\medskip

\begin{Thm}\label{T:mainrep}  Let $S$ be a $C^2$ complete embedded connected $H$-minimal
surface,
  then either $S$ is a vertical plane, or $S$ is determined by a
  generalized seed curve.
\end{Thm}
\pf  As we did at the beginning of the section, decompose $S=V \cup
S_1 \cup S_2 \cup \dots$.  Lemma \ref{L:singloc} shows us that
the $S_i$ can be extended and joined
over the image of the singular
locus.  Thus, we may assume that for each $S_i$ with seed curve
$\gamma^i(s)|_{s \in (a_i,b_i)}$, $S_i$ contains the union of all rules
emanating from $\gamma^i$ and that $V$ contains only places where the
surface is tangent to a vertical plane that are not in the image of
the singular locus for one of the $\gamma^i$.  The proof of lemma \ref{trivial} shows that if $V$ has
nonempty interior then $V$ must composed of subsets of a vertical planes.  Let
$V=V_1 \cup V_2 \cup ...$ be a decomposition of $V$ into portions of
distinct vertical planes.

By lemma
\ref{L:C2}, we have that for each $S_i$, there exists a vertical plane
$V_i$ and a limit point $p_i$ of $\gamma^i$ in the boundary of
$\overline{S_i}$  so that $S$ is tangent to $V_i$ at $p_i$ and at $p_i$, $W_0^i$ must tend to $\pm
\infty$.  By construction, $V_i$ contains a component of $V$ which, by
abuse of notation, we will also denote $V_i$.  We recall a fact from the proof of lemma \ref{L:C2}:  Since the unit Riemannian normal to the surface is
given by
\[ \frac{W^i}{\sqrt{1+(W^i)^2}} \; \nuX + \frac{1}{\sqrt{1+(W^i)^2}}
\; T\]
and that
\[W^i(s,r)= \frac{W_0^i(s) +r - \frac{r^2}{2} \kappa^i(s)}{1-r\kappa^i(s)}\]
we have that as we approach $p_i$, the Riemannian normal tends to the unit
horizontal Gauss map, $\nuX$, and hence is tangent to a vertical plane
with the same normal as the surface at $p_i$.  We define a set-valued
function $v$ to denote the component of
$V$ to which $S_i$ is tangent at $p_i$ by letting $v(p_i)=V_i$.  This
set-valued function $v$ provides
the extra data needed to illustrate how $S_i$ fits together with the
rest of the pieces of $S$.  Specifically, if $S_i$ is determined by
$\Gamma^i(s)=(\gamma_1^i(s),\gamma_2^i(s),h_0^i(s)), s \in (a_i,b_i)$
and
\[p^i_1 = \lim_{s \ra a_i^+} \Gamma^i(s), \; p^i_2= \lim_{s \ra b_i^-}
    \Gamma^i(s)\]
then $\{v(p^i_1),v(p^i_2)\}$ indicate how $S_i$ meets $V$.  To
summarize, we see
that the surface is then determined by a list of curves
$\{\Gamma_i(s)=(\gamma^i(s),h_0^i(s)), s \in
(a_i,b_i)\}$, one corresponding to each $S_i$ coupled with an associated list
designating vertical planes, $\{v(p^i_1),v(p^i_2)\}$. This data
essentially contains the same information as a generalized seed curve
- it provides a complete description of the H-minimal surface.  To
convert this data to the data described in the definition of a
generalized seed curve, we need to ensure that we may choose the
$\gamma^i$ so that
\begin{equation}\label{eq:gen}
 \lim_{s \ra b_i^-} \gamma^i(s) = \lim_{s \ra a_{i+1}^+} \gamma^{i+1}
\end{equation}

To achieve this, we first re-order the $S_i$ (and therefore the
$\gamma^i$ as well) so that any two of the $S_i$ that meet the same
component of $V$ are listed sequentially.  In other words, we reorder
to ensure that $v(p^i_2)=v(p^{i+1}_1)$ for all $i$.  Second, we will
show that we may change the base-points of the seed curves,
$\{\gamma^i\}$, so that the condition in equation \eqref{eq:gen} is
satisfied.  Assume first that equation \eqref{eq:gen} is not satisfied and
consider $\gamma^{i+1}$.  As $p^i_2$ and $p^{i+1}_1$ lie in the same
vertical plane, their projections to the $xy$-plane, denoted
$\tilde{p}^i_2$ and $\tilde{p}^{i+1}_1$ respectively, lie on the same
line, $L$, and are given by
\begin{equation}\label{eq:seedadj}
 \tilde{p}^i_2 = \lim_{s \ra b_i^-} \gamma^i(s),\; \;
 \tilde{p}^{i+1}_1=\lim_{s \ra a_{i+1}^+} \gamma^{i+1}
\end{equation}
Now, using the representation given in theorem \ref{rep} (and extended
by theorem \ref{T:graphs}), we see that
$L$ is the limit of rules in $S_{i+1}$, i.e.
\[L(r)=\lim_{s \ra a_{i+1}^+} \mathcal{L}_{\gamma^{i+1}(s)}(r)\]
Thus, there exists an $r_0$ so that $\tilde{p}^i_2 = L(r_0)$.  Again,
using the representation of theorem \ref{rep}, we know that the
projection of $S_{i+1}$ to the xy-plane is parameterized by a function
\[F_{i+1}(s,r)=(\gamma^{i+1}_1(s)+r(\gamma^{i+1}_2)'(s),\gamma^{i+1}_2(s)-r(\gamma^{i+1}_1)'(s))
\]
Moreover, $F_{i+1}(s,r_0)$ is a re-parametrization of a seed curve
with base-point at $F_{i+1}(s_0,r_0)$ for some $s_0 \in
(a_{i+1},b_{i+1})$ and
\[\lim_{s \ra a_{i+1}^+} F(s,r_0) = \tilde{p}^i_2\]
So, by repicking $\gamma^{i+1}$ with initial basepoint $F_{i+1}(s_0,r_0)$ for some $s_0 \in
(a_{i+1},b_{i+1})$, $\gamma^i$ and $\gamma^{i+1}$ now satisfy equation
\eqref{eq:gen}.  By repeating this construction for all $i$, we have
that $\{\Gamma^i(s)=(\gamma_1^i(s),\gamma_2^i(s),h_0^i(s)), s \in
(a_i,b_i)\}$ satisfies the definition of a generalized seed curve.
Moreover, we now have that $p^i_2$ and $p^{i+1}_1$ are connected by a vertical
line for all $i$ and so taking all the $\Gamma^i$ patched together with the
vertical segments describes a single curve in $\mathbb{H}^1$.  By
construction, this curve characterizes the H-minimal surface $S$.  To
finish the proof, we note that if $\{\Gamma^i\}$ is an empty list,
then $S$ has trivial seed curve and must be a vertical plane.  $\qed$
\medskip

\begin{Rmk}\label{R:genreg}  We remark that the gluing of graphical pieces to vertical
  pieces that can occur in the previous theorem in the presence of
  regularity assumptions imposes more conditions on the $\Gamma^i$.
  In particular, the condition that $W^i$ tends to $\pm \infty$
  ensures that the resulting glued surface is at least $C^1$.  If $S$
  is a $C^k$ surface, then the normal vector would need be $C^{k-1}$.
  For example, consider the case where a component of $V$ contained a portion of a vertical
  plane that had non empty interior that was glued an $S_i$.  As the normal vector to any portion of a vertical plane is constant,
  we see that a necessary and sufficient condition for the glued surface to be $C^k$
  is that $k-1$ derivatives of the normal vector to $S_i$ tend to zero
  as we tend towards the portion of the boundary of $S_i$ that meets
  the vertical plane.

We note that if one is constructing H-minimal surfaces by hand via
  picking seed curves, height functions and vertical planes, this condition can be quite
  complicated to write down due to the multitude of possibilities
  allowed by a generalized seed curve.
\end{Rmk}

\medskip

\begin{Example}\label{E:gencurve}
\emph{In this example, we will construct a variety of surfaces that are best described in using a
generalized seed curve, and some that can only be described in terms
of a generalized seed curve.  We will define two pieces of the surface,
$S_1$ and $S_2$, that are graphs
over the xy-plane, both of which have seed curve a straight line.
To do so, let
\[\gamma^1(s)=(s,0), \; h_0^1(s)= s^\frac{1}{3}\;\; \; \text{for $s \in (0, \infty)$}\] and let
\[ \gamma^2(s)=(s,0), \; h_0^2(s)=s^\frac{1}{3} \;\;\; \text{for $s \in (-\infty,0)$}\]
From these choices, we have that $S_1$ is given by:
\[\left (s,-r,s^\frac{1}{3}+\frac{sr}{2} \right )\;\; \text{for $s >0$}\]
or, equivalently, by $\left(t-\frac{xy}{2}\right)^3-x=0$.  Similarly, $S_2$ is given as
  the portion of $\left(t-\frac{xy}{2}\right)^3-x=0$ with $x<0$.  These two
  pieces are joined by a single line, $V$, given by
  $(0,r,0)$ where the tangent plane to the glued surface $S_1 \cup V
  \cup S_2$ is equal to the vertical plane $x=0$.
This surface $S=S_1 \cup V \cup S_2$ is characterized by the
generalized seed curve:
\[\{(s,0, s^\frac{1}{3}); \; s \in (0,
  \infty), (s,0,s^\frac{1}{3}); \;  s \in (-\infty,0)\}\]
As we can write this surface as $\left(t-\frac{xy}{2}\right)^3-x=0$, it is
  obviously real analytic.  We note that simply by changing the height
  functions, we get a family of similar H-minimal surfaces.  For
  example, the surface $\left(t-\frac{xy}{2}\right)^n-x=0$ is an
  H-minimal surface for any nonzero integer $n$ and is determined by a generalized
  seed curve as follows.  If $n$ is even, then the generalized seed
  curve is
\[\{(s,0, s^\frac{1}{n}); \; s \in (0,
  \infty), (s,0,-s^\frac{1}{n}); \;  s \in (0,\infty)\}\]
whereas if $n$ is odd, the generalized seed curve is
\[\{(s,0, s^\frac{1}{n}); \; s \in (0,
  \infty), (s,0,s^\frac{1}{n}); \;  s \in (-\infty,0)\}\]
We note that when $n$ is even, the surface
$\left(t-\frac{xy}{2}\right)^n-x=0$ can only be described by a
generalized seed curve.}
\end{Example}

\medskip

\begin{Def}\label{D:constcurv}
An $H$-minimal surface has \emph{constant curvature} if either it has
trivial seed curve, or if the signed curvature $\kappa(s)$ of each
seed curve in its generalized seed curve is constant.
\end{Def}

\medskip

We pause to emphasize two important points.  First, that the
assumption of constant curvature does not imply that the curve in
$\mathbb{H}^1$ defined by the generalized seed curve, $\Gamma = \{(\gamma^i(s),h_0^i(s))\}$, has constant
curvature, merely that the $\gamma^i$ do.  Indeed, Example
\ref{E:gencurve} explicitly shows this:  each of its seed curves are
lines (constant curvature zero) but the lifted curves, $(\gamma^i,h_0^i)$ do not have
constant curvature.   Second, we point out that
the individual seed curves, $\gamma^i$, for an $H$-minimal
surface of constant curvature may have different signed curvatures for
different $i$.  We outline a method to construct examples where this can happen.  We
will construct an $H$-minimal surface which has three pieces, $S_1$ and
$S_2$, which are graphs over the $xy$-plane and a single
line $\mathcal{L}$ so that $\overline{S_1}\cap\overline{S_2}
=\mathcal{L}$.  Moreover, along $\mathcal{L}$, the surface is tangent
to a vertical plane.   To define $S_1$,
we pick the seed curve, $\gamma^1$, to be a straight line and choose the initial height
function, $h_0^1$, so that $(h_0^1)'$ tends to $\infty$ as $s\ra
s_0^+$.  To define $S_2$, pick the seed curve, $\gamma^2$, to be a
circle with $\lim_{s \ra s_0^-}\gamma^2(s)=\lim_{s \ra s_0^+}
\gamma^1(s)$ and choose
the initial height function, $h_0^2$, so that $(h_0^2)'$ tends to
$-\infty$ as $s\ra s_0^-$ and $\lim_{s \ra s_0^-}h_0^2(s)=\lim_{s \ra
  s_0^+} h_0^1(s)$.  This type of example is very similar to Example
\ref{E:gencurve} where instead of picking one seed curve to be a circle
and one to be a line, we have that both seed curves are lines.  We
note that in light of Remark \ref{R:genreg} this type of gluing may
lower the  overall regularity of the surface, unless further
restrictions are placed on $h_0^i$.

\vskip 0.6in


\section{\textbf{Complete connected $C^2$ $H$-minimal graphs have constant
  curvature}}\label{S:main}

\vskip 0.2in

We are now in the position to prove our main theorem of Bernstein type.

\medskip

\begin{Thm}\label{T:main} Let $S$ be a connected $C^2$ $H$-minimal surface.
  Suppose further that $S$ is a graph over some plane, $P$.  Then, $S$ has constant
  curvature.
\end{Thm}

\medskip

We begin with a necessary condition for $S$ to be a graph over some plane $P$.

\medskip

\begin{Lem}\label{const} Suppose $S$ is a complete connected $C^2$
  $H$-minimal surface with non-trivial seed curve which is a graph over
  some plane $P$.
  Moreover, assume that $\gamma(0)=0$, $\gamma'(0)=(0,1)$, and $\gamma_1'(s)$ is not identically zero in a neighborhood of $s=0$.  Then, for
  $s$ in a sufficiently small neighborhood of $s=0$,
\[
\frac{\gamma_1'(s)}{<\gamma(s),\gamma'(s)>}\ =\ \kappa(0)\ .
\]
\end{Lem}

\pf Consider two lines in $\R^3$,
$l_0=x_0+t\vec{v}_0$ and $l_1=x_1+t\vec{v}_1$ where $\vec{v}_0$ and
$\vec{v}_1$ are distinct unit vectors and $x_0 \neq x_1$.  We claim there exists
a plane $P$ in $\R^3$ so that the projections of $\{l_1,l_2\}$ to $P$
are non-intersecting lines
if and only if the unit normal to $P$, $\vec{n}$ is a linear
combination of $v_0$ and $v_1$. To prove this, consider the projection of these lines to $P$
\[
\overline{l}_0\ =\ x_0-<\vec{n},x_0>\vec{n} + t \vec{w}_0\ ,
\]
and
\[
\overline{l}_1\ =\ x_1-<\vec{n},x_1>\vec{n} + t \vec{w}_1\ ,
\]
where
\[
\vec{w}_i\ =\ \vec{v}_i-(<\vec{n},\vec{v}_i>)\vec{n}\ .
\]

These lines do not intersect if and only if $\vec{w}_0 \cross
\vec{w}_1 = \vec{0}$. Letting $\vec{n}_0=\vec{v}_0 \cross
\vec{v}_1$, and writing $\vec{n} = a
\vec{v}_0+b\vec{v}_1+c\vec{n}_0$, we have
\begin{equation*}
\begin{split}
\vec{0} = \vec{w}_0 \cross \vec{w}_1 &= \vec{n}_0 - (<\vec{n},\vec{v}_1>)\vec{v}_0 \cross \vec{n} - (<\vec{n},\vec{v}_0>) \vec{n}
\cross \vec{v}_1 \\
&= \vec{n}_0-b(b\vec{n}_0+c(<\vec{v}_0, \vec{v}_1>
\vec{v}_0-\vec{v}_1))-a(a\vec{n}_0+c(\vec{v}_0-<\vec{v}_0, \vec{v}_1>
\vec{v}_1)) \\
&= \vec{n}_0 - b^2 \vec{n}_0 - a^2 \vec{n}_0 +c\vec{v}_0(-b<\vec{v}_0,
\vec{v}_1>-a)+c\vec{v}_1(a<\vec{v}_0, \vec{v}_1>+b) \\
&= c(c \vec{n}_0-\vec{v}_0(b<\vec{v}_0,
\vec{v}_1>+a)+\vec{v}_1(a<\vec{v}_0, \vec{v}_1>+b))\\
\end{split}
\end{equation*}
As $\{\vec{v}_0,\vec{v}_1,\vec{n}_0\}$ form a linearly independent set
under the assumption that $\vec{v}_0$ and $\vec{v}_1$ are not
parallel, we have that all of the coefficients of the last equation
must be zero.  In particular, $c=0$ and $\vec{n}$ is a linear combination of $\vec{v}_0$
and $\vec{v}_1$.

We apply this to two rules of $S$, the rule through
$(\gamma(0),h_0(0))$,
$l_0(r)=(\gamma(0)+r\nuXp,h_0(0)-\frac{r}{2}<\gamma(0),\gamma'(0)>)$,
and the rule through $(\gamma(s),h_0(s))$,
$l_s(r)=(\gamma(s)+r\nuXp,h_0(s)-\frac{r}{2}<\gamma(s),\gamma'(s)>)$.  If the $H$-minimal surface is the graph over some
plane, $P$, then the projection of these lines to $P$ must be
nonintersecting lines.  Since $\gamma(s)$ is not a line near $s=0$ (via the
assumption $\gamma_1'(s)$ is not identically zero in a neighborhood of
$s=0$), we know that
$\kappa(s)$ is not identically zero.  Thus, for all $s$
in a neighborhood of $s=0$,
$l_0$ and $l_s$
are not parallel. We have $l_0'(r)=(1,0,0)$
and $l_s'(r)=(\gamma_2'(s),-\gamma_1'(s),-\frac{1}{2}<\gamma(s),\gamma'(s)>)$.  The plane determined by these vectors has normal
\[
\vec{n}_s = \left(0, \frac{1}{2}<\gamma(s),\gamma'(s)>,-\gamma_1'(s)\right)
\]

Thus, if there is a single plane over which $S$ is a graph, it is
necessary that $\vec{n}_s$ point in the same direction for all
$s$, i.e., that the unit normal vector is constant.  After
normalizing $\vec{n}_s$, the third component must be constant,
i.e. there is a $C$ so that
\[
\frac{-\gamma_1'(s)}{\sqrt{\frac{1}{4}<\gamma(s),\gamma'(s)>^2+\gamma_1'(s)^2}}\ =\ C\ .
\]
By a simple
computation, this is equivalent to
\[
\frac{\gamma_1'(s)}{<\gamma(s),\gamma'(s)>}\ =\ C'\ .
\]
Since $\frac{d}{ds}<\gamma(s),\gamma'(s)>=1+<\gamma(s),\gamma''(s)>$
tends to $1$ as $s$ tends to zero, we have that there exists
an $\epsilon$ so that $<\gamma(s),\gamma'(s)>$ is nonzero on $(-\epsilon, 0) \cup
(0,\epsilon)$.  Thus, on this set,
the quotient above is continuous.  Hence, $\frac{\gamma_1'(s)}{<\gamma(s),
  \gamma'(s)>} = C'$ for all $s$ in the appropriate open set.
Taking the limit as $s \ra 0$, we see that $C'=\kappa(0)$.  $\qed$

\begin{Cor}  Under the assumptions of the previous lemma, $\gamma$ is
  a circle for a neighborhood of $s=0$.
\end{Cor}
\pf By the previous lemma, we have
\[
\frac{\gamma_1'(s)}{<\gamma(s),\gamma'(s)>}\ =\ \kappa(0)\ .
\]
Algebraically simplifying, we have
\[
\gamma_1'(s)(\kappa(0)\gamma_1(s)-1)+\gamma_2'(s)\kappa(0)\gamma_2(s)
= 0\ ,
\]
which shows that $(\kappa(0)\gamma_1(s)-1,\kappa(0)\gamma_2(s))$ is a constant
multiple of $\nuX^\perp$.  In other words, there exists a constant
$\mu$ so that

\begin{equation}\label{l1}
\kappa(0)\gamma_1(s)-1= \mu \gamma_2'(s)
\end{equation}
and
\begin{equation}\label{l2}
\kappa(0)\gamma_2(s)= \mu \gamma_1'(s)
\end{equation}

Differentiating \eqref{l1}, and substituting the resulting equation in \eqref{l2},
we obtain
\[
\gamma_2(s)\ =\ -\ \frac{\mu^2}{\kappa(0)^2}\gamma_2''(s)\ .
\]

Solving this equation subject to the relation between $\gamma_1$
and $\gamma_2$ and the assumption that $\gamma(0)=(0,0),
\gamma'(0)=(0,1)$, finally yields
\[
\gamma(s)\ =\ \left (-\frac{1}{\kappa(0)}\cos(-\kappa(0)s)+
\frac{1}{\kappa(0)}, -\frac{1}{\kappa(0)} \sin (-\kappa(0) s)\right)\ .
\]

This is a circle centered at $\left(\frac{1}{\kappa(0)},0\right)$.

$\qed$

\medskip

\noindent {\em Proof of Theorem \ref{T:main}:}  By theorem \ref{T:mainrep},
$S$ either has trivial seed curve or is determined by a generalized
seed curve $\Gamma = \{(\gamma^i(s),h_0^i(s))\}$.  Consider any seed
curve $\gamma \in \Gamma$.  By composing with
a left translation, which preserves $H$-minimality, we may assume that
$\gamma(0)=(0,0)$ and $\gamma'(0)=(a,b)$.  We note that if
$\gamma'$ is constant in a neighborhood of $s=0$ then
$\gamma$ is a line in this neighborhood and if $<\gamma, \gamma'>$ is
identically zero in  a neighborhood of $s=0$, then $\gamma$ must be a
circle in  this neighborhood.  In both of these cases, $\gamma$ has
constant curvature in the specified neighborhood.  As this is true for
and $\gamma \in \Gamma$, if $S$ has non-trivial seed curve and is a graph
over some plane, then it must have constant curvature.  $\qed$

\vskip 0.6in


\section{\textbf{Classification of $H$-minimal graphs over the $xy$-plane}}

\vskip 0.2in

We end the paper with a proof of Theorem D of the introduction.

\medskip

\noindent
{\em Proof of Theorem D:}  Suppose $S$ is a $C^2$ complete connected $H$-minimal
surface which is a graph over the $xy$-plane, i.e., $S$ can be
parameterized by $(x,y,f(x,y))$.  Then by theorem \ref{T:graphs}, we
know that we can rewrite $S$ in terms of a seed curve, $\gamma(s)$.
By theorem \ref{T:main}, we have that $\gamma(s)$ is either a circle or a line.  By left
translation, we may assume that we have either
\begin{enumerate}
\item
\[
\gamma(s)\ =\ (as,bs)\ ,
\]
with $a^2+b^2=1$, or
\item
\[
\gamma(s)\ =\ \left(R\sin\left(\frac{s}{R}\right),R\cos\left(\frac{s}{R}\right)-R\right)\ ,
\]
where $R$ is the radius of the circle.
\end{enumerate}

In the first case, then $S$ is given by
\begin{equation}\label{ceq1}
\left(as+br,bs-ar,h_0(s)-\frac{rs}{2}\right)\ ,
\end{equation}
Letting $x=as+br$ and $y=bs-ar$, we solve the equations to find $s$
and $r$ in terms of $x$ and $y$:
\begin{equation*}
\begin{split}
s &= ax+by \\
r &= bx-ay\\
\end{split}
\end{equation*}
Substituting these into equation \eqref{ceq1} yields that the surface
is given by:
\[ (x,y,h_0(ax+by)-\frac{1}{2}ab(x^2-y^2)-\frac{1}{2}(b^2-a^2)xy\]
Left translating this representation by the point $(x_0,y_0,t_0)$
yields the desired equation in the statement of theorem D.

In the second case, we have $S$ given as
\[ \left((R-r)\sin\left(\frac{s}{R}\right),(R-r)\cos\left(\frac{s}{R}\right)-R, h_0(s)-\frac{rR}{2}\sin\left(\frac{s}{R}\right)\right)\ .
\]

From Lemma \ref{L:sl}, we know that this parameterization may
cease to be a diffeomorphism when $r=R$.  The points corresponding
to $\left (s,R\right)$ are
\[\left(0,-R,h_0(s)-\frac{R^2}{2}\sin\left(\frac{s}{R}\right)\right)\ .
\]
Thus, for $S$ to remain a graph over the $xy$-plane, we must have
\[
h_0(s)\ =\ \frac{R^2}{2}\sin\left(\frac{s}{R}\right)\ +\ C\ ,
\]
where $C$ is some constant. In this case, we have
\[
h_0(s)-\frac{rR}{2}\cos\left(\frac{s}{R}\right)=(R-r)\left(\frac{R}{2}\sin\left(\frac{s}{R}\right)\right)+C\ .
\]
Letting
$x=(R-r)\sin\left(\frac{s}{R}\right)$, $y=(R-r)\cos\left(\frac{s}{R}\right)$, $t=(R-r)\left(\frac{R}{2}\sin\left(\frac{s}{R}\right)\right)+C$ we
see that $S$ is given by
\[\left(x,y,\frac{R}{2}x+C\right)\]
as claimed.

$\qed$


\end{document}